%% file: paper.tex
\documentclass{article}

\RequirePackage{bm}
\RequirePackage{mathtools,amssymb}
\RequirePackage{graphicx,xcolor, soul}
\RequirePackage{hyperref}
\RequirePackage{a4wide}

\newcommand{\salto}[1]{\ensuremath{[\![#1]\!]}}
\newcommand{\media}[1]{\ensuremath{\lbrace #1 \rbrace}}
\newcommand{\normal}{\bm{n}}
\newcommand{\norm}[1]{{\left\Vert {#1} \right\Vert}}

\DeclareMathOperator{\tr}{\operatorname{tr}}

\graphicspath{{fig/}{fig/examples/ex1/},{fig/examples/ex2/}}

\newtheorem{problem}{Problem}

\begin{document}

\title{A multi-layer reactive transport model for fractured porous media}

\author{%
  Luca Formaggia$^1$,
  Alessio Fumagalli$^1$,
  and
  Anna Scotti$^1$%\affil{1,}\corrauth
}
\date{\small Dipartimento di Matematica, Politecnico di Milano, piazza Leonardo da Vinci 32, 20133 Milano, Italy}

%% \shortauthors is used in copyright information in the end of the paper
%\shortauthors{L. Formaggia, A. Fumagalli, A. Scotti}
%
%\address{%
%  \addr{\affilnum{1}}{Dipartimento di Matematica, Politecnico di Milano, piazza Leonardo da Vinci 32, 20133 Milano, Italy}
%}
%
%% corresponding author
%\corraddr{anna.scotti@polimi.it
%}
%

\maketitle

\begin{abstract}%232 parole
    An accurate modeling of reactive flows in fractured porous media is a key
    ingredient to obtain reliable numerical simulations of several industrial and environmental applications.
    For some values of the physical parameters we can observe the formation of a
    narrow region or layer around the fractures where chemical reactions are focused. Here the transported solute may precipitate and form a salt, or vice-versa. This phenomenon has been observed and reported in real outcrops. By changing its physical properties this layer might substantially alter the global flow response of the system and thus the actual transport  of solute: the problem is thus non-linear and fully coupled. The aim of this work is to
    propose a new mathematical model for reactive flow in fractured porous
    media, by approximating both the fracture and these surrounding layers via a reduced model. In particular, our main goal is to describe the layer thickness evolution with a new  mathematical model, and compare it to a fully resolved equidimensional model for validation.
    As concerns numerical approximation we extend an operator splitting scheme in time to solve sequentially, at each time step,
    each physical process thus avoiding the need for a
    non-linear monolithic solver, which might be challenging due to the
    non-smoothness of the reaction rate.
    We consider  bi- and tridimensional numerical test cases to asses the accuracy and
    benefit of the proposed model in realistic scenarios.
\end{abstract}
%\keywords{
% fractured porous media, reactive flow in porous media,
%    mixed-dimensional modeling, reduced order modeling, operator splitting}

%\input{luca}
\section{Introduction}
The study of reactive flows in porous media is a challenging problem in a large
variety of applications, from geothermal energy to $CO_2$ sequestration up to the
study of flow in tissues or that of the degradation of monuments and cultural
heritage sites. In many cases the porous material presents networks of fractures
that may greatly affect the flow field.
%Fractures are often filled with normally span a wide range of space
These fractures could be responsible for the fast transport of reactants and
heat and thus, in the proximity of fractures, it is possible to observe strong
geochemical effects such as mineral precipitation, dissolution of transformation
that can significantly alter the structure of the porous matrix. Depending on
the relative speed of reaction and transport, namely depending on Damk\"oler
number, we can observe different patterns: a diffused effect on a large part of
the domain, or steeper concentration profiles leading to mineral precipitation
focusing in thin layers around the fractures.

This work presents a mathematical model for this phenomenon based on a
geometrical model reduction that allows to represent thin, heterogeneous
portions of the domains, such as fractures, as lower dimensional manifolds
immersed in the rock matrix.  The proposed model does indeed follow an important
line of research of flow in fractured porous media where fractures are modeled
as one-codimensional manifolds (typically planar) immersed in porous media.
These models, often indicated as hybrid, or mixed-dimensional, describe the
evolution of flow and related fields inside the fracture using a dimensionally
reduced set of equations, and coupling conditions with the surrounding porous
media. With no pretence of being exhaustive, we give a brief overview of
literature related to the techniques used in this framework. A first
hybrid-dimensional model for the coupling of Darcy's flow in porous media and a
single immersed fracture has been presented in~\cite{Martin2005}, and later
extended to networks of fractures by several authors, among
which~\cite{Formaggia2012,Schwenck2015,Scialo2017}. In all those works
single-phase flow was considered, while in~\cite{Fumagalli2012d,Brenner2018, Aghili2019} the
authors deal with two-phase flow formulation. To treat this class of problems, a
large variety of numerical schemes have been exploited. The literature on the
subject is vary vast, we give here only a few suggestions for the interested
readers.  Discretization methods for this type of problems are broadly
subdivided into conforming and non-conforming. In a conforming method the
computational grid used for the porous media is conformal to that used in the
fractures, which means that the elements of the grid used to discretize the
fractures coincide geometrically with facets of the mesh used for the porous
medium. In this setting, many numerical schemes have been proposed, from
classical finite volume approaches, like in~\cite{Stefansson2018a}, to mimetic
finite differencing~\cite{Antonietti}, gradient schemes~\cite{Brenner2015},
discontinuous Galerkin~\cite{Antonietti2019} and hybrid-high order
schemes~\cite{Chave2018}, just to mention some recent works. We recall also some
literature concerning non-conforming methods, which can be again subdivided into
two subsets. The first concerns  the so called geometrically non-matching
discretizations, where the grid used in the fracture is completely independent
to that of the porous media. Among this type of techniques we mention the
embedded discrete fracture network (e-DFM)~\cite{Fumagalli2015,Tene2017} and
approaches based on the use of eXtended finite elements~\cite{Flemisch2008}. In
the second set we have techniques where the fracture is still geometrically
conforming with the porous media grid, but the computational grid can be
different on the two sides. In this class we mention the framework presented
in~\cite{Boon2018,Nordbotten2018} where a mortaring-type technique is used to
connect the solution on domains of different dimensions. See also
\cite{Flemisch2016a,Berre2020a} for a comparison of some of these models.
\begin{figure}[tb]
    \centering
    \includegraphics[width=0.7\textwidth]{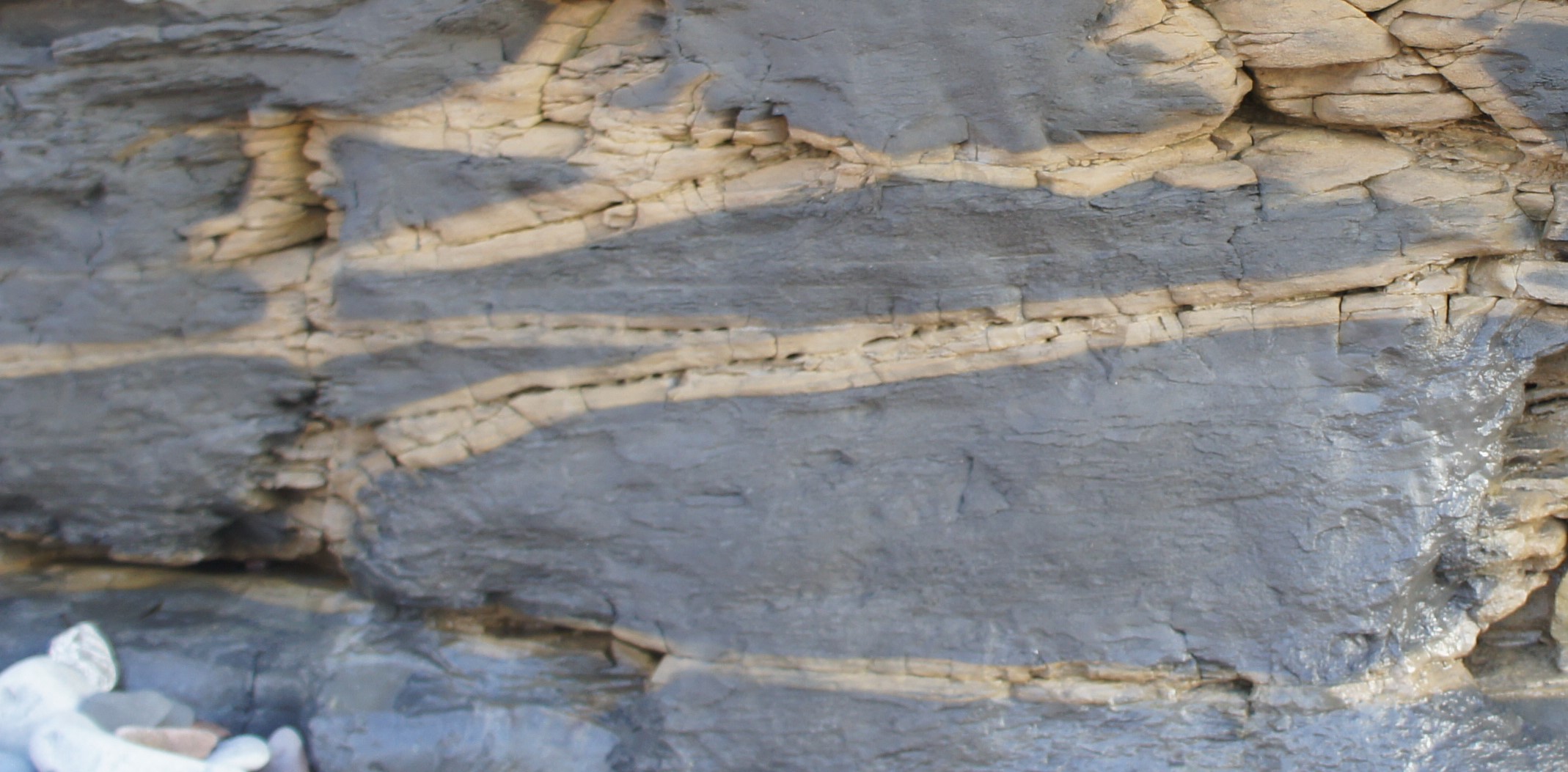}
    \caption{An example of fractures where it is evident the presence of a thin
    layer of altered material in the vicinity of the fracture. The alteration is
    due to geochemical processes driven by the fluid carried by the
    fractures.}\label{fig:alessio}
\end{figure}

In this paper we extend the model presented in~\cite{Fumagalli2020,Fumagalli2020e}, where the
authors developed a model for flow in fractured media accounting for
dissolution-precipitation processes that may alter the flow behavior in of
both fractures and rock matrix.  In ~\cite{Fumagalli2020,Fumagalli2020e} the fracture is represented by an
immersed one-codimensional manifold and special interface conditions were
devised for the diffusion-transport-reaction problem. However, it is known, see
Figure~\ref{fig:alessio}, that the geochemical processes may heavily affect a
very thin layer around the fracture.
Simulating the processes in that region is crucial, but since it is part of the
rock matrix we would need a very fine grid
resolution to obtain an accurate approximation. Consequently, in this work we
consider a model where also those layers are described with a one-codimensional
representation. Thus, the proposed hybrid-dimensional model comprises three
embedded structures, one for the fracture and two for the damage zone
surrounding the fracture at the two sides. We consider the simple reaction model
for solute/mineral reactions illustrated in~\cite{Fumagalli2020e}, in particular
we will consider a single mobile species dissolved in water, representing one of
the two ions in a salt precipitation reaction, and track its transport solving
the single phase Darcy problem and a suitable advection-diffusion-reaction PDE.
Simultaneously, we will keep track of the corresponding precipitate
concentration in the domain.
  %If these reactive regions are thin enough it can become challenging to capture this phenomenon with the computational grid since we want to avoid excessively small grid elements for to reasons: on one hand, because we want to keep the system size under control, in particular since the problem is time dependent; on the other hand, the dimensional reduction of fractures is also based on the assumption that mesh elements are larger than fracture aperture.
The model is similar to the one proposed in
\cite{Faille2014a,Faille2014,Fumagalli2020e,Murad,Fumagalli2019a} to model fault cores and
their surrounding damage zones. It couples three lower
dimensional domains among them and with the surrounding porous matrix by means
of multi-dimensional conservation operators and suitable interface conditions.
This procedure is applied to the Darcy problem and to the evolution equation for
the solute concentration. It could be easily extended to the heat equation to
obtain a more complete physical description of the problem. Another original
contributions of this work consists in the fact that the thickness of the
reactive layers is not fixed a priori, but computed at each time based on the local Darcy velocity and solute concentration.  To this aim, we have derived a
simplified problem on the direction normal to the fracture that provides an
idealized, but useful estimate of the area affected by precipitation.

The numerical discretization is based on a sequential operator splitting
strategy for the decoupling of the equations, and on mixed finite elements for a
good spatial approximation of the fluxes. The model is implemented in the open source
library PorePy, a simulation tool for fractured and
deformable porous media written in Python, see \cite{Keilegavlen2019}.
Some numerical tests are presented, with the aim of verifying the applicability
of the proposed reduced model and its limits, for both two and three-dimensional
settings.

The paper is structured as follows: in Section \ref{sec:model} the single and
multi-layer mathematical model is introduced and described in details. We
introduce also the model to describe the evolution of the layer surrounding the
fracture. Section \ref{sec:numerical_discretization} defines the numerical
discretization, in space and time. In particular, a splitting scheme in time is
detailed to allow the solution of each physical process sequentially. In Section
\ref{sec:results} we present the numerical test cases for the comparison between
the new model and the one already present in literature. Finally, Section \ref{sec:conclusion}
is devoted to the conclusions.

\section{Mathematical model}\label{sec:model}

Let us start by illustrating the governing equations before performing
dimensional reduction of the fracture region.  We will consider here a simple
setting with a single fracture, and depict the domains in two dimensions for
simplicity, even if the presentation is given in a general setting and
three-dimensional results will be presented in Section \ref{sec:results}.

\begin{figure}[htb]
    \centering
    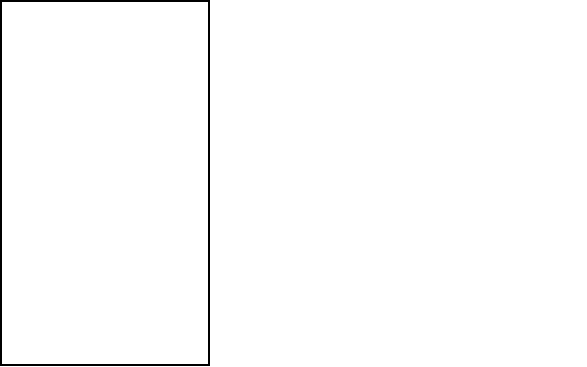
    \caption{Equi-dimensional representation of the rock matrix $\Omega$, the
    fracture $\gamma$ and surrounding layers $\mu$.}%
    \label{fig:domain_equi}
\end{figure}

Let $\Theta\subset \mathbb{R}^d$, with $d=2$ or $3$, be the domain filled by
porous material, where we can identify three parts, as depicted in Figure
\ref{fig:domain_equi}: the porous matrix $\Omega$, occupying the larger part of
the domain; the fracture $\gamma$, characterized by a small thickness, called
aperture, and a disconnected subdomain $\mu$, formed by two layers ($\mu^-$ and
$\mu^+$) adjacent the fracture at both sides. The domain $\Omega$ is split in
two disjoint parts $\Omega^+$ and $\Omega^-$ by the two sides of the layer.
Clearly, $\overline{\Theta} = \overline{\Omega} \cup \overline{\gamma} \cup
\overline{\mu}$ and  $\Omega$, $\gamma$, and $\mu$ have mutually disjoint
interior. In the following, barred quantities are given boundary data.

We assume that $\Theta$ is filled by a single phase fluid, water, with
constant density, and that average fluid velocity $\bm{q}_\Theta$ and pressure
$p_\Theta$
can be obtained as the solution of the Darcy's problem
\begin{align}\label{eq:darcy}
    \begin{aligned}
    &\begin{aligned}
        &k_\Theta^{-1} \bm{q}_\Theta +  \nabla p_\Theta =\bm{0}\\
        &\partial_t \phi_\Theta + \nabla \cdot \bm{q}_\Theta + f_\Theta =0
    \end{aligned}
    \quad&&\text{in } \Theta \times (0, T)\\
    &p_\Theta = \overline{p_\Theta} && \text{on } \partial \Theta_p\times(0, T)\\
    &\bm{q}_\Theta \cdot \bm{\nu} = \overline{q_\Theta} && \text{on } \partial
    \Theta_q\times(0, T)
    \end{aligned}
\end{align}
where $\phi_\Theta$ denotes the porosity (variable in space and time), $k_\Theta =
k_\Theta(\phi_\Theta)$ is the
intrinsic permeability tensor (already divided by fluid viscosity), which can
depend on porosity and may show large variations among the three different subdomains,
and $f_\Theta$ is a volumetric forcing term. The boundary $\partial \Theta$ is
subdivided into two disjoint subsets $\partial \Theta_p$ and $\partial \Theta_u$
such that $\overline{\partial\Theta}=\overline{\partial\Theta_p} \cup
\overline{\partial\Theta_q}$. We assume that $\partial\Theta_p\neq \emptyset$.
The $\overline{p_\Theta}$ and $\overline{q_\Theta}$ are given boundary
conditions. Note that, even if the subdomains are characterized by
different physical parameters, we have continuity of pressure and flux at the
interface between $\Omega$ and $\mu$, and $\mu$ and $\gamma$. Finally we denote with
$T>0$ the final simulation time.

The Darcy's problem is coupled with a simple chemical system with two species, \cite{Radu2010}: a solute $U$, whose
concentration is denoted by $u$, and a precipitate $W$, whose concentration is
denoted by $w$. The solute $U$ can represent the anion and cation in
a salt precipitation model. Thanks to the usual assumption of electrical
equilibrium, the concentrations of these two species are equal. The solute is
transported by water, therefore its evolution is governed by an advection-diffusion-reaction
equation for $u$, while that of the  precipitate can be described
by an ordinary differential equation for $w$ at each point in $\Theta$. We have then
\begin{subequations}\label{eq:problem_solute_precipitate}
\begin{align}
    \begin{aligned}
    &
    \begin{aligned}\label{eq:solute}
        &\bm{\chi}_\Theta - \bm{q}_\Theta u_\Theta + \phi_\Theta D_\Theta \nabla
        u_\Theta = \bm{0}\\
        &\partial_t (\phi_\Theta u_\Theta) + \nabla \cdot \bm{\chi}_\Theta +
        \phi_\Theta r_w(u_\Theta, w_\Theta; \theta_\Theta) =
        0
    \end{aligned}
    \quad&&\text{in } \Theta \times (0, T)\\
    &u_\Theta = \overline{u_\Theta} &&\text{on } \partial \Theta_u \times(0,
    T)\\
    &\bm{\chi}_\Theta\cdot \bm{\nu} = \overline{\chi_\Theta} && \text{on } \partial
    \Theta_\chi \times (0, T)\\
    &u_\Theta(t=0) = u_{\Theta, 0} && \text{in } \Theta \times \{0\}
    \end{aligned},
\end{align}
and
\begin{align}
    \begin{aligned}
        &\partial_t (\phi_\Theta w_\Theta) - \phi_\Theta r_w(u_\Theta, w_\Theta;
        \theta_\Theta) = 0 \label{eq:precipitate}
        \quad &&\text{in } \Theta \times (0, T)\\
        &w_\Theta (t=0) = w_{\Theta, 0} && \text{in } \Theta \times \{0\}
    \end{aligned},
\end{align}
\end{subequations}
Here the problem is presented in mixed form and $\bm{\chi}_\Theta$ is the total flux
accounting for advection and diffusion. $D_\Theta$ is the diffusion coefficient and
$r_w$ the reaction rate, whose expression depends on the type of reaction
considered. In the following we will use a linear (oversimplified) model where
\begin{gather}\label{eq:chimica1}
    r_w(u, w; \theta) =\lambda(\theta)u,
\end{gather}
as well as a more complex model, taken from \cite{Bouillard}, and used in~\cite{Agosti2016,Fumagalli2020e},
\begin{gather}\label{eq:chimica2}
    r_w(u, w; \theta) =\lambda(\theta)
    \left\{\max[r(u)-1,0]+H(w)\min[r(u)-1,0]\right\}.
\end{gather}
In both cases, the reaction rate depends on $\lambda$ (which can be a constant or
depend on the local temperature according to Arrhenius law) and on the reactant
concentration. While in \eqref{eq:chimica1} the transformation of $U$ into $W$
proceeds in a single direction until $u=0$, the more realistic equation
\eqref{eq:chimica2} could describe a reaction that proceeds in both directions
depending on the solute concentration compared to the equilibrium one (taken equal to
one in this a-dimensional setting). It also accounts for the fact that mineral
dissolution must stop when $w=0$, hence the dependence on the Heaviside function
$H(w)=\max(0,w)$, \cite{Knabner1995}.

Finally, the porosity $\phi_\Theta$ can change in time as the result of
mineral precipitation with the following law
\begin{gather}\label{phi_w}
    \begin{aligned}
        &\partial_t \phi_\Theta + \eta_\Theta \phi_\Theta \partial_t w_\Theta =
        0 \quad&& \text{in } \Theta\times (0,T)\\
        &\phi_\Theta(t=0) = \phi_{\Theta, 0} &&\text{in }\Theta \times
        \{0\}
    \end{aligned},
\end{gather}
with $\eta_\Theta$ being a positive parameter. See \cite{Noorden2009} for an in-depth discussion of the microscale phenomenon at the basis of \eqref{phi_w}.

%definizione di Damk\"oler number
The transport-reaction process can be characterized by means of the Damkh\"oler
number, which can be interpreted as the ratio between the characteristic times
of transport and reaction~\cite{bearmodeling}. If the dominant transport
mechanism is advection we can define the first Damk\"ohler number as
\begin{gather*}
    Da_I=\dfrac{\lambda \phi L}{\norm{\bm{q}}}
\end{gather*}
where $L$ is the characteristic length of the phenomenon. Conversely, if diffusion
is prevalent, one should consider the second Damk\"oler number
\begin{gather*}
    Da_{II}=\dfrac{\lambda  L^2}{D}.
\end{gather*}
In both cases, a large Damk\"ohler number means that reaction is fast compared to
transport and will result in a precipitation (or dissolution) concentrated in
space. In this work, we treat situations arising for high Damk\"ohler number.
\begin{problem}[Equi-dimensional problem]\label{pb:equi}
    The problem of reactive transport in the porous media $\Theta\times(0, T)$ gives
    $(\bm{q}_\Theta, p_\Theta, \bm{\chi}_\Theta, u_\Theta, w_\Theta,
    \phi_\Theta)$ by solving the coupled equations \eqref{eq:darcy},
    \eqref{eq:problem_solute_precipitate}, \eqref{phi_w}.
\end{problem}

\subsection{Standard fracture-matrix flow and transport model}

\begin{figure}[htb]
    \centering
    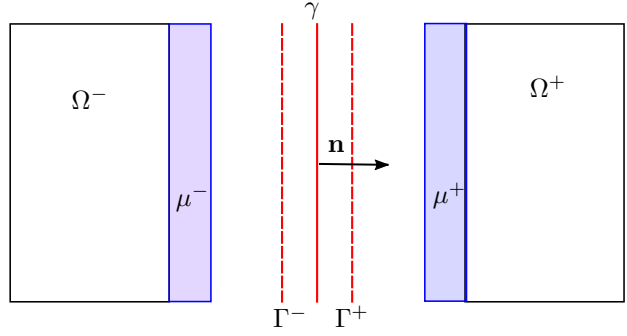
    \caption{Mixed-dimensional representation of the rock matrix $\Omega$, the equi-dimensional layers $\mu$  and the lower-dimensional fracture $\gamma$.}%
    \label{fig:domain_mono}
\end{figure}

We are interested in the effect of mineral precipitation on fractured porous
media.
In the standard setting, like the one illustrated in~\cite{Fumagalli2020}, the portion $\mu$ is
still considered as part of the $d$-dimensional domain, while fractures are modeled as
lower dimensional entities, since they are characterized by a small aperture compared to
the other characteristic lengths. We indicate with $\overline{\Psi} = \overline{\Omega}
\cup \overline{\mu}$ and $\overline{\Psi^\pm} = \overline{\Omega^\pm} \cup
\overline{\mu^\pm}$.
A sketch of the domain is shown
in Figure \ref{fig:domain_mono}, where $\gamma$ indicates now, with a slight
abuse of notation, the center line of the fracture, with aperture $\epsilon_\gamma$. While, with $\Gamma$ we indicate $\partial\mu\cap\gamma$, i.e.
the portion of the boundary of the porous matrix that coincides geometrically with $\gamma$. Indeed, $\Gamma$ is formed by two parts, $\Gamma^+$ and $\Gamma^-$, corresponding to
the $+$ and $-$ parts of the porous matrix, on the right and left side of the fracture in Figure~\ref{fig:domain_mono}, respectively.
We remark that in the figure, $\Gamma$ is drawn separately form $\gamma$, but in
fact $\Gamma$ and $\gamma$ coincide geometrically.

To make the notation more compact in the hybrid-dimensional setting, from now on
we use the following convention.  When no subscript is present a scalar and
vector field is understood as the compound variable of fields defined in the different
hybrid-dimensional domains.  For instance,
$\bm{q}=(\bm{q}_\Psi,\bm{q}_\gamma)$ represents the fluxes in the rock matrix
$\Psi$ and in the fracture $\gamma$, each indicated with the corresponding
subscript. Analogously for $p = (p_\Psi, p_\gamma)$. Moreover, in the following, for a given field $f$ we indicate with
$\tr_{\beta} f$ the trace of $f$ on $\beta$. In particular, $\tr_{\Gamma^-}$ and
$\tr_{\Gamma^+}$ indicate the trace operators on the two parts of $\Gamma$.

We can define the jump operator for a scalar function $p$ and the normal component of a vector function $\bm{q}$, as
\begin{gather*}
    \salto{p}_\gamma=\tr_{\Gamma^+}p_\Psi -\tr_{\Gamma^-}p_\Psi
    \quad \text{and} \quad
    \salto{\bm{q}\cdot \bm{n}}_\gamma=\tr_{\Gamma^+}(\bm{q}_\Psi \cdot\bm{n})-\tr_{\Gamma^-} (\bm{q}_\Psi\cdot\bm{n}),
\end{gather*}
where $\bm{n}$ is the normal to $\gamma$ pointing towards the $+$ side.
We can also define the average operators,
\begin{gather*}
    \media{p}_\gamma=\frac{1}{2}(\tr_{\Gamma^+}p +\tr_{\Gamma^-}p)
    \quad \text{and} \quad
    \media{\bm{q}\cdot \bm{n}}_\gamma=\frac{1}{2}(\tr_{\Gamma^+}(\bm{q}_\Psi\cdot\bm{n})+\tr_{\Gamma^-} (\bm{q}_\Psi\cdot\bm{n})).
\end{gather*}

In this framework, the governing equations should be formulated for the variables $\bm{q}_\Psi$ and $p_\Psi$
in the porous matrix domain $\Psi$, and for the flux $\bm{q}_\gamma$, and the pressure $p_\gamma$ in the
fracture $\gamma$.
We note that $\bm{q}_\gamma$ is aligned along $\gamma$, i.e. $\bm{q}_\gamma\cdot\bm{n}=0$, and we can define on $\gamma$
a mixed-dimensional divergence $\nabla_\gamma\cdot$ as
\begin{gather*}
    \nabla_{\gamma} \cdot \bm{q} = {\nabla} \cdot
    \bm{q}_{\gamma} - \salto{\bm{q}\cdot \bm{n}}_\gamma,
\end{gather*}
where ${\nabla} \cdot\bm{q}_{\gamma}$ is the standard divergence on the
tangent space of $\gamma$ and
the jump term accounts for the exchange between fracture and porous matrix.  We
note that in a two-dimensional setting like the one depicted in
Figure~\ref{fig:domain_mono}, ${\nabla} \cdot\bm{q}_{\gamma}=\partial_y
\bm{q}_{\gamma}$, where $y$ is, in general, the intrinsic
coordinate of $\gamma$. More details on those operators may be found in the
cited literature. In this case the boundary of $\Psi$ is divided in the
following three non-intersecting subsets $\overline{\partial \Psi} =
\overline{\partial \Psi_u} \cup \overline{\partial \Psi_q} \cup
\overline{\Gamma}$, with the similar division also for the boundary in the
solute equation.

The resulting mixed dimensional set of equation is, in the domain $\Psi$,
\begin{subequations}\label{eq:model}
\begin{align}
    \begin{aligned}
    &
    \begin{aligned}
        &k^{-1}_\Psi \bm{q}_\Psi + \nabla p_\Psi = \bm{0}\\
        &\partial_t \phi_\Psi + \nabla \cdot \bm{q}_\Psi + f_\Psi = 0
    \end{aligned}
    \quad&&\text{in } \Psi\times(0, T)\\
    &p_\Psi = \overline{p_\Psi} &&\text{on } \partial
    \Psi_{p}\times(0, T)\\
    &\bm{q}_\Psi \cdot \bm{\normal}
    = \overline{q_\Psi} &&\text{on } \partial \Psi_{u}\times(0, T)
    \end{aligned}
\end{align}
and also in the fracture $\gamma$
\begin{align}
    \begin{aligned}
    &
    \begin{aligned}
        &\epsilon_\gamma^{-1} k_\gamma^{-1} \bm{q}_\gamma + \nabla p_\gamma = \bm{0}\\
        &\partial_t \epsilon_\gamma + \nabla_\gamma \cdot \bm{q} + f_\gamma = 0
    \end{aligned}
    \quad&&\text{in } \gamma\times(0, T)\\
    &p_\gamma = \overline{p_\gamma} &&\text{on } \partial
    \gamma_{p}\times(0, T)\\
    &\bm{q}_\gamma \cdot \bm{\normal}
    = \overline{q_\gamma} &&\text{on } \partial \gamma_{u}\times(0, T).
    \end{aligned}
\end{align}
Note that the mixed-dimensional divergence couples the equations in the porous matrix with those in the fracture.
%Note that the
%subdomains $\mu$ and $\Omega$ may have different properties, such as porosity
%and permeability, but flux and pressure are continuous in $\Psi$.
Equations~\label{eq:mode} are complemented with the following interface conditions
on $\Gamma$,
\begin{align}\label{eq:model_cc}
    \begin{aligned}
     &k_\gamma^{-1}\epsilon_\gamma  \media{\bm{q}\cdot \bm{n}}_\gamma -
    \salto{p}_\gamma = 0\\
    &\dfrac{k_\gamma^{-1}\epsilon_\gamma}{4}  \salto{\bm{q}\cdot \bm{n}}_\gamma + p_{\gamma} -
    \media{p}_\gamma = 0
    \end{aligned}
    \quad \text{in } \Gamma \times (0, T)
\end{align}
\end{subequations}
where we have assumed an isotropic permeability $k_\gamma$ in the fracture, i.e. permeability is the same in the tangential and normal direction.
The first condition \eqref{eq:model_cc} states that the net flux of $\bm{q}_{\Psi}$ through $\gamma$ is
proportional to the jump of pressure across the fracture, while the
second states that the flux exchange between porous matrix and fracture is
proportional to the difference between the pressure in the fracture and the
average pressure in the surrounding porous medium. We may note that the second
relation is a particular case of that proposed
in~\cite{Formaggia2018,Martin2005}, where a family of conditions have been
proposed depending on a modeling parameter.

Accordingly, the advection-diffusion-reaction problem can be written in
mixed-form in the rock matrix as
\begin{subequations}\label{eq:modelADR}
\begin{gather}
    \begin{aligned}
        &
        \begin{aligned}
            & \bm{\chi}_\Psi -\bm{q}_\Psi u_\Psi + \phi_\Psi D_\Psi\nabla u_\Psi = \bm{0}\\
            &\partial_t (\phi_\Psi u_\Psi) + \nabla \cdot \bm{\chi}_\Psi + r_\Psi = 0
        \end{aligned}
        \quad &&\text{in } \Psi\times(0, T)\\
        &u_\Psi = \overline{u_\Psi} &&\text{on }\partial \Psi_u \times (0, T)\\
        &\bm{\chi}_\Psi \cdot \bm{\nu} = \overline{\chi_\Psi} && \text{on }
        \partial \Psi_\chi \times (0, T)\\
        &u_\Psi(t=0) = u_{\Psi, 0} && \Psi \times \{0\}
    \end{aligned}
\end{gather}
and in the fracture $\gamma$ as
\begin{gather}
    \begin{aligned}
        &
        \begin{aligned}
            &\bm{\chi}_\gamma - \bm{q}_\gamma u_\gamma + \epsilon_\gamma  D_\gamma \nabla u_\gamma = \bm{0}\\
            &\partial_t (\epsilon_\gamma u_\gamma) + \nabla_\gamma \cdot \bm{\chi} + r_\gamma = 0
        \end{aligned}
        &&\text{in } \gamma\times(0, T)\\
        &u_\gamma = \overline{u_\gamma} && \text{on } \partial\gamma_u \times
        (0, T)\\
        &\bm{\chi}\cdot \bm{\nu} = \overline{\chi_\gamma} && \text{on }
        \partial\gamma_\chi \times (0, T)\\
        &u_\gamma(t=0) = u_{\gamma, 0} && \text{in } \gamma \times \{0\}
    \end{aligned}
\end{gather}
with the same definition for the mixed dimensional divergence operator and a
similar interface conditions
\begin{align}\label{eq:modelADR_cc}
    \begin{aligned}
        &D_\gamma^{-1} \epsilon_\gamma \left(\media{\bm{\chi}\cdot\bm{n}}_\gamma - \media{\bm{q}\cdot\bm{n}}_\gamma \right) -\salto{u}_\gamma=0\\
        &\dfrac{D_\gamma^{-1}\epsilon_\gamma}{4}  \left(\salto{\bm{\chi}\cdot \bm{n}}_\gamma - \salto{\bm{q}\cdot \bm{n}}_\gamma\right) + u_{\gamma} -
    \media{u}_\gamma = 0
    \end{aligned}
    \quad \text{in } \Gamma \times (0, T).
\end{align}
\end{subequations}
The porosity $\phi_\Psi$ evolves in time according to
\eqref{phi_w}, and fracture aperture can vary due to mineral
precipitation with the following law:
\begin{gather}\label{eq:reduced_aperture}
    \begin{aligned}
        &\partial_t \epsilon_\gamma  + \eta_\gamma \epsilon_\gamma \partial_t
        w_\gamma = 0&&
        \text{in } \gamma \times (0, T)\\
        &\epsilon_\gamma(t=0) = \epsilon_{\gamma, 0} && \text{in } \gamma \times \{0\}.
    \end{aligned}
\end{gather}
\begin{problem}[Fracture mixed-dimensional problem]\label{pb:fractured}
    The problem of reactive transport in the fractured porous media gives in
    $\Psi\times(0, T)$ the fields
    $(\bm{q}_\Psi, p_\Psi, \bm{\chi}_\Psi, u_\Psi, w_\Psi,
    \phi_\Psi)$ and in
    $\gamma\times(0, T)$ the fields
    $(\bm{q}_\gamma, p_\gamma, \bm{\chi}_\gamma, u_\gamma, w_\gamma,
    \epsilon_\gamma)$
    by solving the coupled equations \eqref{eq:model},
    \eqref{eq:modelADR}, \eqref{phi_w} for $\phi_\Psi$, and
    \eqref{eq:reduced_aperture}.
\end{problem}

\subsection{Multi-layer flow and transport model}
In the previous section we have revised a mixed-dimensional model where only the fracture is treated
as a lower dimensional interface. However, if we assume that fractures play a major role on fluid flow and
solute transport, we can identify cases in which the Damk\"oler number is high, and consequently
the precipitation (or dissolution) of minerals is concentrated in a thin region close to the fracture. This
occurs, for instance, if solute is injected in clean water through a fracture,
the fracture is more permeable than the surrounding domain and reaction is
significantly faster than transport. In this case, as shown in
\cite{Fumagalli2020e}, the solute profile decays rapidly in a thin
region near the fracture. It is then difficult to capture the
phenomenon numerically without resorting to a very fine grid in the porous region near the fracture
where most geochemical reactions occurs, which we call reactive layer.

To reduce the computational cost,  we propose here a three
layers model where also the reactive layers $\mu$ surrounding the fracture are
represented as lower dimensional domains, of thickness $\epsilon_\mu$, suitably
coupled with the fracture on one side and the porous matrix on the other side.
The derivation of such multi-layer model is similar to the one presented
in~\cite{Faille2014a,Faille2014,Fumagalli2020e,Murad,Fumagalli2019a}, where its
introduction was motivated by the modelling of faults and their surrounding
damage zone.

\begin{figure}[htb]
    \centering
    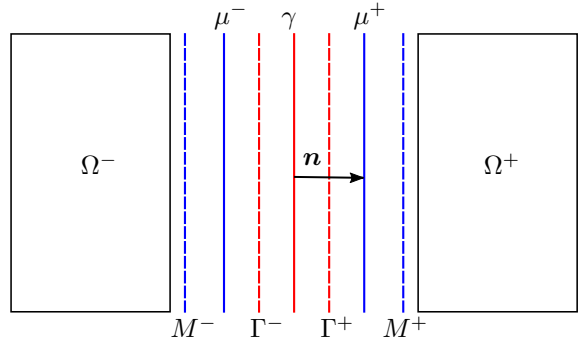
    \caption{Mixed-dimensional representation of the rock matrix $\Omega$,
    damage zone $\mu$,
    and fault $\gamma$.}%
    \label{fig:domain}
\end{figure}

To keep the notation simple, we preserve the same notation used in the previously described model, even if the domains are geometrically different,
since $\mu$ is now formed by two lower dimensional reactive layers $\mu^-$ and $\mu^+$, located at each side of the fracture $\gamma$. Moreover,
we let $M=\{M^-,M^+\}$ denote the interface between $\Omega$ and $\mu$, while $\Gamma=\{\Gamma^-,\Gamma^+\}$ is now the
interface between $\mu$ and $\gamma$, see Figure \ref{fig:domain}. Note that
even if $\gamma$, $\mu$, $M$, $\Gamma$ are geometrically superimposed, they
play a different role in the model: lower dimensional domains and interfaces,
respectively.

In addition to $\bm{q}_\Omega$, $p_\Omega$,
$\bm{q}_\gamma$, and $p_\gamma$, we define the flux $\bm{q}_\mu$,
and the average pressure $p_\mu$ in $\mu$. Similarly,
$u_\mu$, $w_\mu$ will denote the concentrations in $\mu$,  and $\bm{\chi}_\mu$ the relative flux. We follow also here the convention
that fields without a subscript identify the collection of quantities in the different domains.

While now $M$ can be identified as the part of the boundary of $\Omega$ that coincides with the model of the fracture,
here $\Gamma=\{\Gamma^-,\Gamma^+\}$ are fictitious additional interfaces needed to define the coupling, and on which we define
the normal fluxes $q_\Gamma$ and $\chi_\Gamma$, both scalar.

We also need to revise the definition of jump and average operators. In particular,
\begin{gather*}
    \salto{p}_{\mu^-} = p_{\gamma} - \tr_{M^-} p_\Omega\quad \text{and}\quad
    \salto{p}_{\mu^+}=\tr_{M^+} p_\Omega - p_{\gamma},\\
    \media{p}_{\mu^-} = \frac{1}{2}(p_{\mu^-} + \tr_{M^-} p_\Omega) \quad \text{and} \quad
    \media{p}_{\mu^+} =\frac{1}{2}(\tr_{M^+} p_\Omega - p_{\mu^+}),
\end{gather*}
and
\begin{gather*}
    \salto{\bm{q}\cdot\bm{n}}_{\mu^-} = q_{\Gamma^-} - \tr_{M^-}(\bm{q}_\Omega\cdot\bm{n})
    \quad \text{and} \quad
    \salto{\bm{q}\cdot\bm{n}}_{\mu^+}=\tr_{M^+} (\bm{q}_\Omega\cdot\bm{n}) - q_{\Gamma^+},\\
    \media{\bm{q}\cdot\bm{n}}_{\mu^-} = \frac{1}{2}(q_{\Gamma^-} - \tr_{M^-}(\bm{q}_\Omega\cdot\bm{n}))
    \quad \text{and}\quad
    \media{\bm{q}\cdot\bm{n}}_{\mu^+} =\frac{1}{2}(\tr_{M^+} (\bm{q}_\Omega\cdot\bm{n}) - q_{\Gamma^+}),
\end{gather*}
depending on whether we are considering $\mu^-$ or $\mu^+$ of $\mu$, respectively. While,
\begin{gather*}
    \salto{p}_\gamma = p_{\mu^+} - p_{\mu^-}
    \quad \text{and} \quad
    \media{p}_\gamma = \frac{1}{2}(p_{\mu^-} + p_{\mu^+}),\\
    \media{\bm{q}\cdot\normal}_\gamma = q_{\Gamma^+} - q_{\Gamma^-}
    \quad\text{and}\quad
    \media{\bm{q}\cdot\normal}_\gamma = \frac{1}{2}(q_{\Gamma^-} + q_{\Gamma^+}).
\end{gather*}
Analogous definitions for hold for $u$, $w$ and $\bm{\chi}$.

We are now in the position to define the mixed dimensional divergence operators in this new
setting: given a vector field $\bm{q}$
we have
\begin{gather}\label{eq:gen_divergence}
    \nabla_{\mu} \cdot \bm{q} = \nabla \cdot
    \bm{q}_{\mu} - \salto{\bm{q}\cdot\normal}_\mu\quad \text{and}\quad
    \nabla_\gamma \cdot \bm{q} = \nabla \cdot \bm{q}_\gamma -\salto{\bm{q}\cdot\normal}_\gamma,
\end{gather}
where, following the convention,  $\bm{q}_{\mu}$ and  $\bm{q}_{\gamma}$ are the
components of $\bm{q}$ in the corresponding lower dimensional domains, while
$\nabla_\gamma\cdot$  and $\nabla_\mu\cdot$ the divergence operator acting on
the corresponding domain.

We now write the differential problem representing the new mixed-dimensional
model, where we also impose boundary conditions for the  flux and for the
pressure on portions of the boundaries of $\Omega$, $\mu$ and $\gamma$,
indicated by the subscript $u$ and $p$, respectively.  Note that
$\overline{\partial\Omega}=\overline{\partial\Omega_u}\cup\overline{\partial\Omega_p}\cup
\overline{M}$, with a similar division also for the boundary in the solute
equation. We also assume, for
well-posedness, that $\partial\Omega_p$ is not empty.

In the porous matrix we have
\begin{subequations}\label{eq:model3}
\begin{gather}
    \begin{aligned}
        &
        \begin{aligned}
            &k_\Omega^{-1} \bm{q}_\Omega + \nabla p_\Omega = \bm{0}\\
            &{\partial_t \phi_\Omega} + \nabla \cdot \bm{q}_\Omega + f_\Omega = 0
        \end{aligned}
        &&\text{in } \Omega\times(0, T)\\
        &p_\Omega = \overline{p_\Omega} &&\text{on } \partial
        \Omega_{p}\times(0, T)\\
        &\bm{q}_\Omega \cdot \bm{\normal}
        = \overline{q_\Omega} &&\text{on } \partial \Omega_{u}\times(0, T)
    \end{aligned},
\end{gather}
while for the layer $\mu$ we have
\begin{gather}\label{2c}%
    \begin{aligned}
        &
        \begin{aligned}
            &\epsilon_\mu^{-1}k_\mu^{-1} \bm{q}_{\mu} + \nabla p_{\mu} = \bm{0}\\
            &\partial_t (\epsilon_\mu\phi_\mu)+\nabla_{\mu} \cdot \bm{q} + f_{\mu} = 0
        \end{aligned}
        && \text{in } \mu\times(0, T)\\
        &p_{\mu} = \overline{p_\mu} &&\text{on } \partial \mu_p\times(0, T)\\
        &\bm{u}_\mu \cdot \bm{\nu} = \overline{u_\mu} &&\text{on } \partial
        \mu_{u}\times(0, T)
    \end{aligned},
\end{gather}
and, finally, for in the fracture we have
\begin{gather}\label{2d}%
    \begin{aligned}
        &
        \begin{aligned}
            &\epsilon_\gamma^{-1} k_\gamma^{-1 }\bm{q}_\gamma + \nabla p_\gamma = \bm{0}\\
            &{\partial_t \epsilon_\gamma}+ \nabla_\gamma \cdot \bm{q} + f_\gamma = 0
        \end{aligned}
        &&\text{in } \gamma\times(0, T)\\
        &p_\gamma = \overline{p_\gamma} &&\text{on } \partial
        \gamma_p\times(0, T)\\
        &\bm{u}_\gamma \cdot \bm{\nu} = \overline{u_\gamma} &&\text{on } \partial
        \gamma_{u}\times(0, T)
    \end{aligned}.
\end{gather}
Note that the fracture is considered filled just by fluid, and that the flow velocity is sufficiently small to model it using lubrication theory,
which gives an equation akin to Darcy's with a ``porosity'' equal to $1$. Fracture aperture can change as an effect of
precipitation.  Moreover, we have the following interface conditions on $M$ and
$\Gamma$, respectively,
\begin{gather}\label{interface}
    \begin{aligned}
    &
    \begin{aligned}
        &\epsilon_\mu k_\mu^{-1} \media{\bm{q}\cdot \bm{n}}_\mu - \salto{p}_\mu = 0\\
        &\dfrac{\epsilon_\mu k_\mu^{-1}}{4} \salto{\bm{q}\cdot \bm{n}}_\mu + p_{\mu} - \media{p}_\mu = 0
    \end{aligned}
    \quad &\text{on } M \times (0, T)\\
    &
    \begin{aligned}
        &\epsilon_\gamma k_\gamma^{-1} \media{\bm{q}\cdot \bm{n}}_\gamma - \salto{p}_\gamma = 0\\
        &\dfrac{\epsilon_\gamma k_\gamma^{-1}}{4} \salto{\bm{q}\cdot \bm{n}}_\gamma + p_{\gamma} - \media{p}_\gamma = 0
    \end{aligned}
    & \text{on } \Gamma \times (0, T).
    \end{aligned}
\end{gather}
\end{subequations}
Similarly, the transport
and reaction problem in the multi-layer domain becomes, first for the rock matrix
\begin{subequations}\label{eq:modelADR3}
\begin{gather}
    \begin{aligned}
        &
         \begin{aligned}
            & \bm{\chi}_\Omega -\bm{q}_\Omega u_\Omega + \phi_\Omega D_\Omega\nabla u_\Omega = \bm{0}\\
            &\partial_t (\phi_\Omega u_\Omega) + \nabla \cdot \bm{\chi}_\Omega + r_\Omega = 0
        \end{aligned}
        \quad &&\text{in } \Omega\times(0, T)\\
        &u_\Omega = \overline{u_\Omega} &&\text{on } \partial
        \Omega_{u}\times(0, T)\\
        &\bm{\chi}_\Omega \cdot \bm{\nu}
        = \overline{\chi_\Omega} &&\text{on } \partial \Omega_{\chi}\times(0,
        T)\\
        &u_\Omega(t=0) = u_{\Omega, 0}&& \text{in } \Omega \times \{0\}
    \end{aligned},
\end{gather}
while for the layer $\mu$ we have
\begin{gather}\label{2cADR}
    \begin{aligned}
        &
        \begin{aligned}
            & \bm{\chi}_{\mu}-\bm{q}_\mu u_\mu + \epsilon_\mu \phi_\mu D_\mu\nabla u_\mu = \bm{0}\\
            &{\partial_t (\epsilon_\mu\phi_\mu u_\mu)} +\nabla_{\mu} \cdot \bm{\chi} + r_{\mu} = 0
        \end{aligned}
        \quad&& \text{in } \mu\times(0,T)\\
        &u_{\mu} = \overline{u_\mu} &&\text{on } \partial \mu_u\times(0, T)\\
        &\bm{\chi}_\mu \cdot \bm{\nu} = \overline{\chi_\mu} &&\text{on } \partial
        \mu_{\chi}\times(0, T)\\
        &u_\mu(t=0) = u_{\mu, 0}&& \text{in } \mu \times \{0\}
    \end{aligned},
\end{gather}
and finally for the fracture $\gamma$
\begin{gather}\label{2dADR}
    \begin{aligned}
        &
        \begin{aligned}
            &\bm{\chi}_\gamma - \bm{q}_\gamma u_\gamma + \epsilon_\gamma  D_\gamma \nabla u_\gamma = \bm{0}\\
            &\partial_t (\epsilon_\gamma u_\gamma) + \nabla_\gamma \cdot \bm{\chi} + r_\gamma = 0
        \end{aligned}
        \quad&&\text{in } \gamma\times(0, T)\\
        &u_\gamma = \overline{u_\gamma} &&\text{on } \partial
        \gamma_u\times(0, T)\\
        &\bm{\chi}_\gamma \cdot \bm{\nu} = \overline{\chi_\gamma} &&\text{on } \partial
        \gamma_{\chi}\times(0, T)\\
        &u_\gamma(t=0) = u_{\gamma, 0}&& \text{in } \gamma \times \{0\}.
    \end{aligned}
\end{gather}
\end{subequations}
The interface conditions on $\Gamma$ and $M$ similar to \eqref{interface} to
couple concentrations and fluxes in the subdomains,
\begin{gather}\label{interfaceADR}
    \begin{aligned}
    &
    \begin{aligned}
        &\epsilon_\mu D_\mu^{-1} \left(\media{\bm{\chi}\cdot \bm{n}}_\mu-\media{\bm{q}\cdot \bm{n}}_\mu \right)- \salto{p}_\mu = 0\\
        &\dfrac{\epsilon_\mu D_\mu^{-1}}{4} \left(\salto{\bm{\chi}\cdot \bm{n}}_\mu-\salto{\bm{q}\cdot \bm{n}}_\mu \right)+ p_{\mu} - \media{p}_\mu = 0
    \end{aligned}
    \quad &\text{on } M \times (0, T)\\
    &
    \begin{aligned}
        &\epsilon_\gamma D_\gamma^{-1}\left( \media{\bm{\chi}\cdot \bm{n}}_\gamma-\media{\bm{q}\cdot \bm{n}}_\gamma \right)- \salto{p}_\gamma = 0\\
        &\dfrac{\epsilon_\gamma D_\gamma^{-1}}{4} \salto{\bm{q}\cdot \bm{n}}_\gamma + p_{\gamma} - \media{p}_\gamma = 0.
    \end{aligned}
    &\text{on } \Gamma \times (0, T).
    \end{aligned}
\end{gather}
In this multi-layer model the
porosities $\phi_\Omega$ and $\phi_\mu$ depend on the corresponding values of
precipitate concentration according to \eqref{phi_w}, and fracture aperture
$\epsilon_\gamma$ follows \eqref{eq:reduced_aperture}. The only missing part is a
model for the evolution of the thickness $\epsilon_\mu$, which will be discusses
in the next section.
\begin{problem}[Multi-layer fractured mixed-dimensional problem]\label{pb:multi}
    The problem of reactive transport in the multi-layer fractured porous media gives in
    $\Omega \times(0, T)$ the fields
    $(\bm{q}_\Omega, p_\Omega, \bm{\chi}_\Omega, u_\Omega, w_\Omega,
    \phi_\Omega)$, in $\mu\times(0, T)$ the fields
    $(\bm{q}_\mu, p_\mu, \bm{\chi}_\mu, u_\mu, w_\mu, \phi_\mu, \epsilon_\mu)$,
    in
    $\gamma\times(0, T)$ the fields
    $(\bm{q}_\gamma, p_\gamma, \bm{\chi}_\gamma, u_\gamma, w_\gamma,
    \epsilon_\gamma)$, and in $\Gamma \times (0, T)$ the interface fluxes
    $(q_\Gamma, \chi_\Gamma)$
    by solving the coupled equations \eqref{eq:model3},
    \eqref{eq:modelADR3}, \eqref{phi_w} for $\phi_\Omega$ and $\phi_\mu$, and
    \eqref{eq:reduced_aperture} for $\epsilon_\gamma$. While for $\epsilon_\mu$ one of the model
    discussed in Subsection \ref{subsec:model_layer}.
\end{problem}

\subsection{A model for layer thickness}\label{subsec:model_layer}

\begin{figure}[tb]
    \centering
    \includegraphics[width=0.45\textwidth]{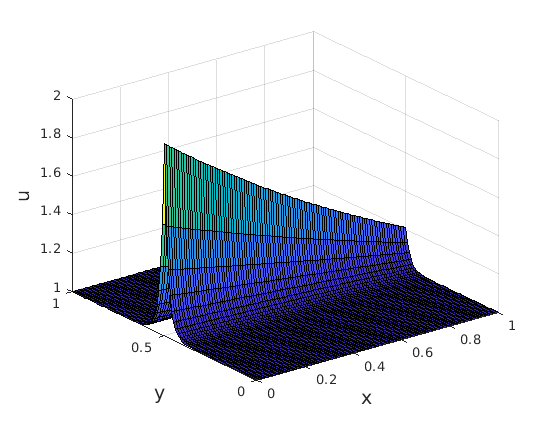}%
    \includegraphics[width=0.45\textwidth]{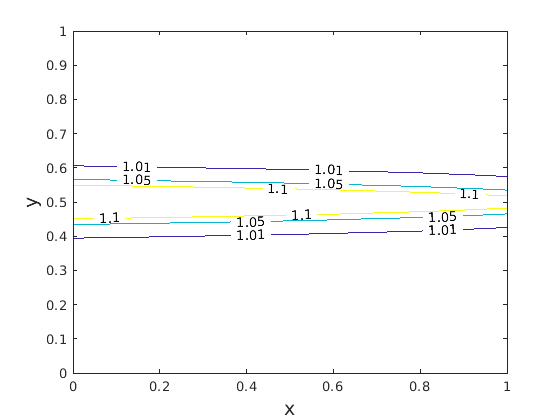}
    \caption{On the left, an example of solute concentration distribution around
    a fracture. On the right, the corresponding isolines.}
    \label{profilo}
\end{figure}

We want to obtain a model for the thickness of the layers $\mu$, i.e. we want to
model $\epsilon_\mu$ as a function of the physical parameters and the solution
itself, to compute values that can change in space and in time accounting for chemical reactions. We recall that we
assume that
there is a well-identifiable region, around the fracture, where dissolution or
precipitation take place, and that this region is ``thin'' if reaction is
sufficiently fast with respect to the transport mechanism of interest, advection
and/or diffusion. However, we cannot obtain this information from the solute and
precipitate distribution in the porous matrix due, in practice, to
insufficient grid resolution. For this reason we have resorted to
one-dimensional models that will allow us to compute analytical solutions for the evolution of the layer in
simplified settings. In particular we assume that
\begin{itemize}
 \item the transport of solute near the fracture can be approximated as
 one-dimensional in the direction normal to the fracture, for each section;
 \item the changes in porosity due to precipitation have a small impact on the
 advection field;
 \item solute is transported more easily in the fracture, thus the concentration
 of solute in the fracture can be considered as a boundary condition for its
 diffusion/advection in the neighboring layers;
 \item the Damk\"oler number is such that, from the solute profile we can, after
 fixing a cutoff concentration value, find a small thickness $\epsilon_\mu$ for
 each of the two layers $\mu^+$ and $\mu^-$ at each time $t$.
\end{itemize}
Consider Figure \ref{profilo}: starting from the solute concentration in the
fracture we obtain the concentration profile in the neighborhood. If, for
instance, we consider a precipitation model such that precipitation occurs where
$u>1$ then the region $\mu$ is encompassed by the corresponding concentration
isoline. If $\epsilon_\mu$ is small enough, it is reasonable to use the proposed mixed-dimensional model, by
collapsing $\mu$  into a lower dimensional domain, as explained previously.

\subsubsection{Pure advection, linear reaction}

In the simplified case of no diffusion, and with the advection field given by
the Darcy velocity normal to the fracture, we can obtain an analytical
expression for the solute concentration under the assumptions stated above. We
are assuming that flux is exiting the fracture, i.e. the normal Darcy velocity
$Q$ is positive and can be considered constant in time.

If we denote with $s$ the arc length in the direction normal to the fracture the
one-dimensional problem for solute concentration reads:
\begin{gather*}
    \begin{aligned}
        &\phi{\partial_t u} + Q{\partial_s u} = -\lambda \phi u \quad && \text{in  }
        (0, +\infty) \times (0, T)\\
        &u = u_\gamma && \text{on } 0\times (0, T)\\
        &u(t=0) = 0 && \text{in } (0, +\infty) \times \{0\}
    \end{aligned}
\end{gather*}
and has the exact solution  $u(s,t) = u_0\left(s-{Q}t/\phi\right)\exp(-\phi\lambda/Q s)$ where
\begin{gather*}
    u_0=\begin{cases}
            u_\gamma &  s=0\\
            0 & s>0
           \end{cases}
    \quad \text{and} \quad
    u=\begin{cases}
        u_\gamma\exp(-\phi\lambda/Q t) & s\leq {Q}t/\phi\\
        0 & s> {Q}t/\phi
    \end{cases}.
\end{gather*}
Note that, with this linear reaction term, we have precipitation whenever $u>0$,
however, in practice, we can choose a cut-off value, i.e. the layer is defined  by the condition $u(s,t)>\delta$.

Thus, we seek the point $\overline{s}=\epsilon_\mu$ where $u(\epsilon_\mu,t)=\delta$.  We obtain
\begin{gather}\label{growth}
    \epsilon_\mu= Q/ \phi \min(t,
    -\ln({\delta}/{u_\gamma})/\lambda) ,
\end{gather}
 i.e. the layer thickness grows
linearly in time until it reaches its steady state value. The time to reach the
steady state can be estimated as
$\overline{t}=\ln({\delta}/{u_\gamma})/{\lambda}$.

\subsubsection{A more realistic reaction term}

The linear decay term considered in the previous section is however too simple
for most diagenetic processes. For the case of mineral precipitation, under some
simplifying assumptions, one can consider the reaction term given in~\eqref{eq:chimica2}.
If we consider just the case of mineral precipitation, i.e. we assume that the solution
is supersaturated, its expression simplifies to
\begin{gather*}
     r_w(u) =-\lambda\left(u^2-1\right).
\end{gather*}
Under the assumptions stated in the previous section, we can estimate the
thickness of the layer where precipitation occurs by solving the following one
dimensional problem in the direction normal to the fracture
\begin{gather*}
    \begin{aligned}
        &
        \phi {\partial_t u} +Q{\partial_s u} = -\lambda\phi\left(u^2-1\right)
        \quad && \text{in  }
        (0, +\infty) \times (0, T)\\
        &u = u_\gamma && \text{on } 0\times (0, T)\\
        &u(t=0) = 0 && \text{in } (0, +\infty) \times \{0\}
    \end{aligned}.
\end{gather*}
At steady state, with ${\partial_t u} =0$, the problem above admits the exact solution
\begin{gather*}
    u(s)=\dfrac{C+\exp({2\lambda \phi s}/{Q})}{\exp({2\lambda \phi
    s}/{Q})-C}\quad\text{with}\quad
    C=\dfrac{u_\gamma-1}{u_\gamma+1}
\end{gather*}
to satisfy the boundary condition at
the interface with the fracture. Note that, with this reaction model,
precipitation only occurs where $u>1$ i.e. where the concentration is above
equilibrium. Therefore we consider a cutoff value $\overline{u}=1+\delta$ with
$\delta<1$ a small enough number, and seek the corresponding layer thickness
\begin{gather}\label{eq:thickness_layer_non_linear}
    \dfrac{C+\exp({2\lambda\phi \epsilon_\mu}/{Q})}{\exp({2\lambda \phi
    \epsilon_\mu}/{Q})-C}=1+\delta
%C+\exp(\frac{2\lambda \phi\epsilon_\mu}{Q})=(1+f)(\exp(\frac{2\lambda\phi \epsilon_\mu}{Q})-C)\\
%\exp(\frac{2\lambda \phi \epsilon_\mu}{Q}) = C\dfrac{f+2}{f}\\
    \quad \Rightarrow \quad
    \epsilon_\mu=\dfrac{Q}{2\lambda \phi}\log\left(C\dfrac{2+\delta}{\delta}\right).
\end{gather}
Once again the steady state layer thickness depends linearly on the ratio
$Q/\lambda$. In this case however it is more difficult to obtain an
expression for its growth in time: for this reason, in the results section, we
will just verify this estimate and defer the actual application of this model to
future work.

\section{Numerical approximation}\label{sec:numerical_discretization}

In this section we discuss the approximation strategies adopted to solve the
model presented in Problem \ref{pb:multi}, in particular the spatial
and temporal approximation schemes and the procedure to solve the resulting
coupled and non-linear system. In Subsection \ref{subsec:time_discretization} we
consider the temporal discretization of the problem along with the splitting
algorithm, which can be considered an extension of the one introduced and
studied in \cite{Fumagalli2020e}. In Subsection
\ref{subsec:spatial_discretization} we will briefly present the spatial
discretization adopted.

\subsection{Time discretization and splitting}\label{subsec:time_discretization}

The global physical problem, in \eqref{eq:model3}, involves several processes that are
coupled in a non-linear way. To overcome the need for a  monolithic non-linear
solver, and rely more on legacy simulation codes, suited for each single physical process,
we consider a splitting strategy in time, such that each equation can be solved
separately. However, we recall that the operator splitting approach usually
introduces an additional error in time.  Furthermore, since some of our physical
variables (i.e., porosity, solute, and precipitate) are very sensitive to volume
changes we also need to design the splitting strategy such that no mass or
volume is unexpectedly lost. Finally, since the reaction term for the solute may
be rather complex and highly non-linear, an additional operator splitting is
employed to separate the diffusive and advective part from the reaction in
equations \eqref{eq:modelADR3}. In this way, we can use ad-hoc numerical schemes to solve the
latter.

For these reasons, we extend the strategy developed in \cite{Fumagalli2020e} to our
needs, in particular incorporating the physical processes linked to the reactive layers $\mu$.
The extension is quite straightforward, however we recall the splitting algorithm
for reader's convenience. We divide the time interval in $N$ steps and we denote
with $t^n = n\Delta t$, with $\Delta t$ the time step assumed constant for simplicity. We set the initial condition as
\begin{gather*}
    \phi^0_\Omega = \phi_{\Omega, 0} \quad \epsilon^0_\gamma=\epsilon_{\gamma,
    0}\quad
    \phi^0_\mu = \phi_{\mu, 0} \quad \epsilon^0_\mu = \epsilon_{\mu, 0}
    \quad
    \theta^0_\Omega=\theta_{\Omega, 0} \quad \theta^0_\gamma=\theta_{\gamma,
    0}\quad
    \theta^0_\mu = \theta_{\mu, 0}\quad\\
    u^0_\Omega=u_{\Omega, 0} \quad u^0_\gamma = u_{\gamma, 0} \quad u^0_\mu =
    u_{\mu, 0}\quad
    w^{-1}_\Omega = w^0_\Omega=w_{\Omega, 0} \quad w^{-1}_\gamma =
    w^0_\gamma=w_{\gamma, 0}\quad w^{-1}_\mu =
    w^0_\mu =w_{\mu, 0}.
\end{gather*}

In each time step $(t^n, t^{n+1})$, we perform the following steps.
\begin{enumerate}
    \item To get a better estimate
    of the porosity as well as the fracture aperture computed in the
    Step \ref{step:porosity_2}, we extrapolate the concentration of the precipitate as in
        \cite{Giovanardi2015,Agosti2015}. We obtain
        \begin{gather*}
            w^*_\Omega=2w^n_\Omega-w^{n-1}_\Omega
            \quad \text{and} \quad
            w^*_\gamma=2w^n_\gamma-w^{n-1}_\gamma
            \quad \text{and} \quad
            w^*_\mu=2w^n_\mu-w^{n-1}_\mu
            .
        \end{gather*}
    \item\label{step:porosity_2} We then compute the porous media and layer
        porosity and fracture aperture, from \eqref{phi_w} for $\phi_\Omega$ and
        $\phi_\mu$ and \eqref{eq:reduced_aperture} for $\epsilon_\gamma$, by the following relations
        \begin{gather*}
            \phi_\Omega^*=\frac{\phi_\Omega^n}{1+ \eta_\Omega(w^*_\Omega-w^{n}_\Omega)}
            \quad \text{and} \quad
            \phi_\mu^*=\frac{\phi_\mu^n}{1+ \eta_\mu(w^*_\mu-w^{n}_\mu)}
            \quad \text{and} \quad
            \epsilon^*_\gamma=\frac{\epsilon^n_\gamma}{1+ \eta_\gamma(w^*_\gamma-w^{n}_\gamma)}.
        \end{gather*}
        Note that we do not compute an estimate of the thickness layer
        $\epsilon_\mu$ since the models presented in Subsection \ref{subsec:model_layer} are not differential.
    \item To prepare the computation of the pressure and Darcy velocity, we update
    the  permeability of the porous media $k_\Omega(\phi_\Omega^*)$ as well as fracture and layer permeabilities
    $k_\gamma(\epsilon^*_\gamma)$ and
    $k_\mu(\epsilon^{n}_\mu, \phi_\mu^*$, respectively.
    \item We solve the Darcy problem
        \eqref{eq:model3}
        to get pressure and Darcy velocity in the domain, in the fracture and in
        the layer: $(\bm{q}^{n+1}_\Omega, p^{n+1}_\Omega)$, $(\bm{q}^{n+1}_\gamma,
        p^{n+1}_\gamma)$, and $(\bm{q}^{n+1}_\mu, p^{n+1}_\mu)$, respectively,
        as well as the interface flux $q_\Gamma$ on the interface $\Gamma$.
        For the discretization of the temporal derivative of porosity
        $\phi_\Omega$ and $\phi_\mu$ and fracture aperture $\epsilon_\gamma$ we
        consider both their value predicted in Step \ref{step:porosity_2} and at
        time $n-1$.
    \item \label{step:solute_ad} We solve
        the advection-diffusion part of the  solute equation,
        \eqref{eq:modelADR3}, to obtain an intermediate value of the solute:
        $u^{n+\frac{1}{2}}_\Omega$, $u^{n+\frac{1}{2}}_\gamma$, and
        $u^{n+\frac{1}{2}}_\mu$. Also the interface flux $\chi_\Gamma$ is
        computed on $\Gamma$.
        Note that we do not consider for this point the reaction term.
    \item In the previous point we have accounted for porosity changes
        using $\phi_\Omega^*$ and $\phi_\mu^*$, as well as fracture aperture
        changes using $\epsilon_\gamma^*$. The intermediate value of the solute
        $u_\Omega^{n+\frac{1}{2}}$,
        $u_\gamma^{n+\frac{1}{2}}$ and $u_\gamma^{n+\frac{1}{2}}$ already accounts for the change in pore
        volume, then also the precipitate in the porous domain, fracture, and
        layer have to be updated to account for the same variation
        \begin{gather*}
            w^{n+\frac{1}{2}}_\Omega=w^n_\Omega \frac{\phi_\Omega^n}{\phi_\Omega^*}
            \quad \text{and} \quad
            w^{n+\frac{1}{2}}_\gamma=w^n_\gamma
            \frac{\epsilon^n_\gamma}{\epsilon^*_\gamma}
            \quad \text{and} \quad
            w^{n+\frac{1}{2}}_\mu=w^n_\mu
            \frac{\phi^n_\mu}{\phi^*_\mu}.
        \end{gather*}
    \item \label{step:solute_r} We then solve the reaction step for both the solute and precipitate, by starting from the values of
    $(w^{n+\frac{1}{2}}_\Omega,
    w^{n+\frac{1}{2}}_\gamma, w^{n+\frac{1}{2}}_\mu)$ and $(u^{n+\frac{1}{2}}_\Omega,
    u^{n+\frac{1}{2}}_\gamma, u^{n+\frac{1}{2}}_\mu)$
    to get the values of $(w^{**}_\Omega, w^{**}_\gamma, w_\mu^{**})$ and
    $(u^{**}_\Omega, u^{**}_\gamma, u^{**}_\mu)$.
    \item Since the precipitate has changed in the previous step, we need
    to update the porosity of the porous matrix and layer as well as the
    fracture aperture. Considering the model \eqref{phi_w} for $\phi_\Omega$ and
    $\phi_\mu$ and \eqref{eq:reduced_aperture} for $\epsilon_\gamma$,
    we  obtain
        \begin{gather*}
            \phi_\Omega^{n+1}=\frac{\phi_\Omega^n}{1+ \eta_\Omega(w^{**}_\Omega-w^{n}_\Omega)}
            \quad \text{and} \quad
            \epsilon^{n+1}_\gamma=\frac{\epsilon^n_\gamma}{1+
            \eta_\gamma(w^{**}_\gamma-w^{n}_\gamma)}
            \quad \text{and} \quad
            \phi_\mu^{n+1}=\frac{\phi_\mu^n}{1+ \eta_\mu(w^{**}_\mu-w^{n}_\mu)}
        \end{gather*}
    We also compute the thickness of the layer $\epsilon_\mu^{n+1}$ by following
    one of the model presented in Subsection \ref{subsec:model_layer}.
    \item As the last point in the algorithm, we update the solute and precipitate concentrations to account for the variation of porosity and fracture aperture
    at the previous point. We compute
        \begin{gather*}
            w^{n+1}_\Omega=w^{**}_\Omega \frac{\phi_\Omega^*}{\phi_\Omega^{n+1}}
            \quad
            u^{n+1}_\Omega=u^{**}_\Omega
            \frac{\phi_\Omega^*}{\phi_\Omega^{n+1}}
            \quad
            w^{n+1}_\gamma=w^{**}_\gamma
            \frac{\epsilon^*_\gamma}{\epsilon^{n+1}_\gamma}
            \\
            u^{n+1}_\gamma=u^{**}_\gamma
            \frac{\epsilon^*_\gamma}{\epsilon^{n+1}_\gamma}
            \quad
            w^{n+1}_\mu=w^{**}_\mu
            \frac{\phi^*_\mu }{\phi^{n+1}_\mu}
            \quad
            u^{n+1}_\mu=u^{**}_\mu
            \frac{\phi^*_\mu}{\phi^{n+1}_\mu}.
        \end{gather*}
\end{enumerate}

The set of ordinary differential equations to be solved in Step
\ref{step:solute_r} depends on the reaction function chosen. See
\cite{Pop2017,Fumagalli2020e,Bringedal2020} for an example. For the
time discretization of Step \ref{step:solute_r} in our case we have considered a
second order Runge-Kutta scheme. For the other equations the first order
Implicit Euler scheme is used for their temporal discretization. The operator
splitting approach also introduces an error which, in our case, is of order one
in time. Globally, we obtain a first order scheme in time.

\subsection{Spatial discretization}\label{subsec:spatial_discretization}

The spatial discretization considered for the full problem is specific for each
physical phenomenon and for each spatial dimension, since the schemes for
fracture and layer are written on their tangent space.  We consider schemes with
compatible degrees of freedoms, meaning that they are associated only to cells
(primary variables) and faces (fluxes), and no interpolation operators will be
required. Since the focus of the present work is not on innovative spatial
discretizations to solve the problem, but rather on the model, and since we use
well known schemes, we only briefly mention them.

To compute a reliable Darcy velocity, which is then used as an input in the
other problems, the numerical method has to be locally mass conservative and
provide a good quality approximation of the fluxes. For this reason, our choice
is to discretize the pressure equation, in its mixed form, with the lowest-order
Raviart-Thomas finite element for the Darcy velocity and piece-wise constant
elements for the pressure fields. This scheme is also particularly suited for
strong permeability variations typical of the underground. See
\cite{Raviart1977,Roberts1991,Boffi2013} for a more detailed discussion.

For the numerical solution of the solute and temperature fields, we consider a
two-point flux approximation for the diffusion operator and a weighted upstream
for the advective part. See \cite{Flemisch2016a,Stefansson2018a,Keilegavlen2019}
for a more extensive discussion.

The coupling between the subdomains (porous media, fracture, and layer) is done via
Lagrange multipliers that represent the normal flux exchange between them. See
\cite{Boon2018,Nordbotten2018,Boon2020,Keilegavlen2019} for more details and analysis.

\section{Results}\label{sec:results}

In this section we present two groups of test cases to validate the previously
introduced model.  In the first group of test cases, in Subsection
\ref{subsec:case1}, we consider a 2D domain with one fracture, adapting the
geometry of second example of \cite{Fumagalli2020e} to our needs. In this
geometry we compare the classical fracture-matrix model described in Problem \ref{pb:fractured} with the new
multi-layer model in Problem \ref{pb:multi}, for increasing levels of complexity in the physical
parameters. In Subsection \ref{subsec:case2} instead we consider a test case in
three-dimensions, by adapting the geometry and data of Case 1 of
\cite{Berre2020a}. In all the examples, the considered numerical scheme cannot
handle the case of zero fracture aperture or layer thickness: for this reason,
at the initial time when the reactive layer has not started developing yet, we
will set a very low starting value for $\epsilon_\mu$. See \cite{Boon2018} for a
different approach that is able to handle vanishing fracture aperture. Since the
presented model for the layer thickness evolution considers mostly an advective
field as main driving force, we will set the diffusion coefficient for the
solute transport problem to a low value to obtain results that are in agreement
with the theory.

The following examples are implemented with the Python library PorePy
\cite{Keilegavlen2019} and the scripts of each test case are freely available on
GitHub.

Finally, even if the current model may be coupled with a heat equation as in
\cite{Fumagalli2020e} in these experiments we consider a given, constant
temperature field and therefore a fixed and uniform in space reaction rate
$\lambda$.

\subsection{Two-dimensional problem}\label{subsec:case1}

\begin{figure}[tb]
    \centering
    \resizebox{0.32\textwidth}{!}{\fontsize{0.75cm}{2cm}\selectfont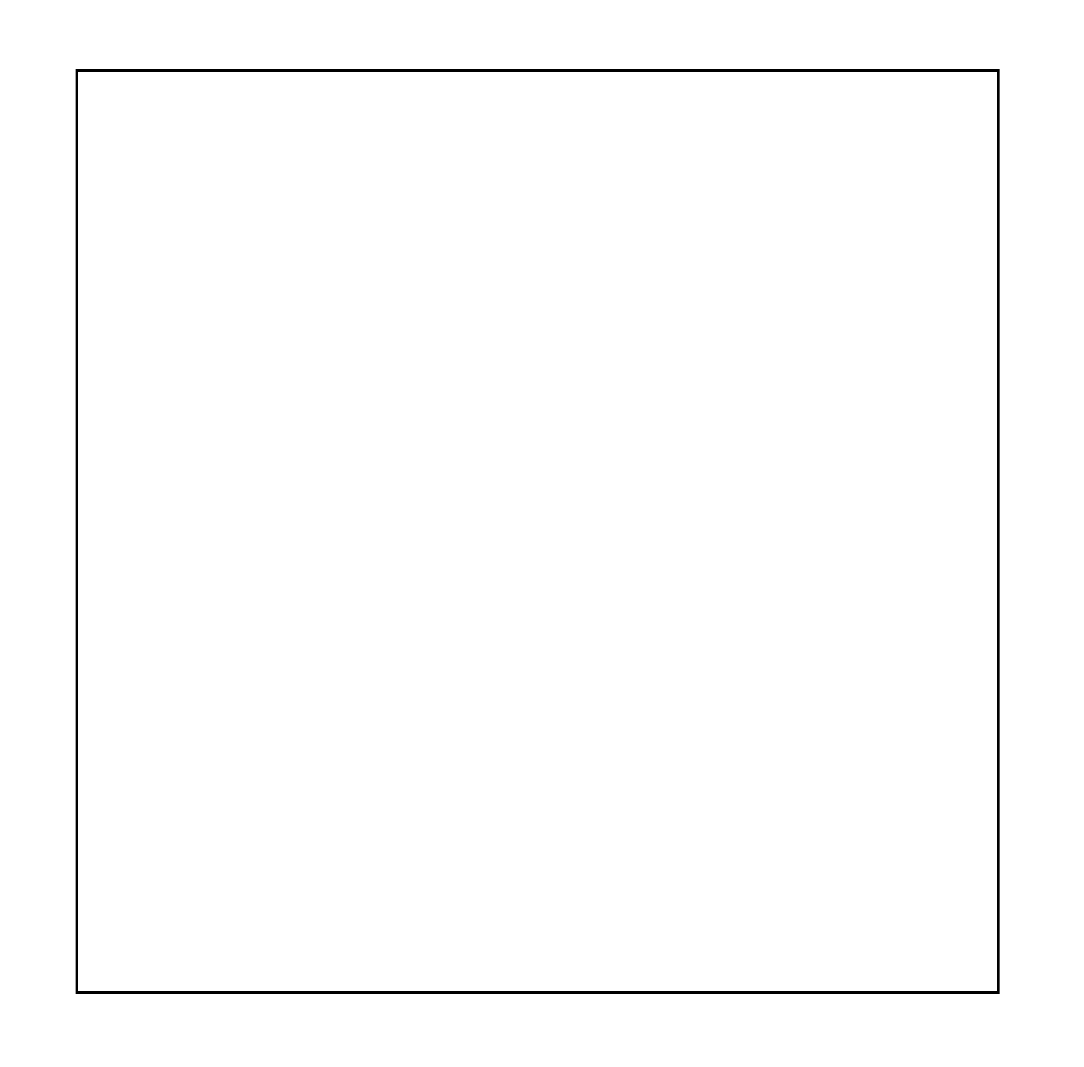}%
    \hspace*{0.01\textwidth}%
    \includegraphics[width=0.32\textwidth]{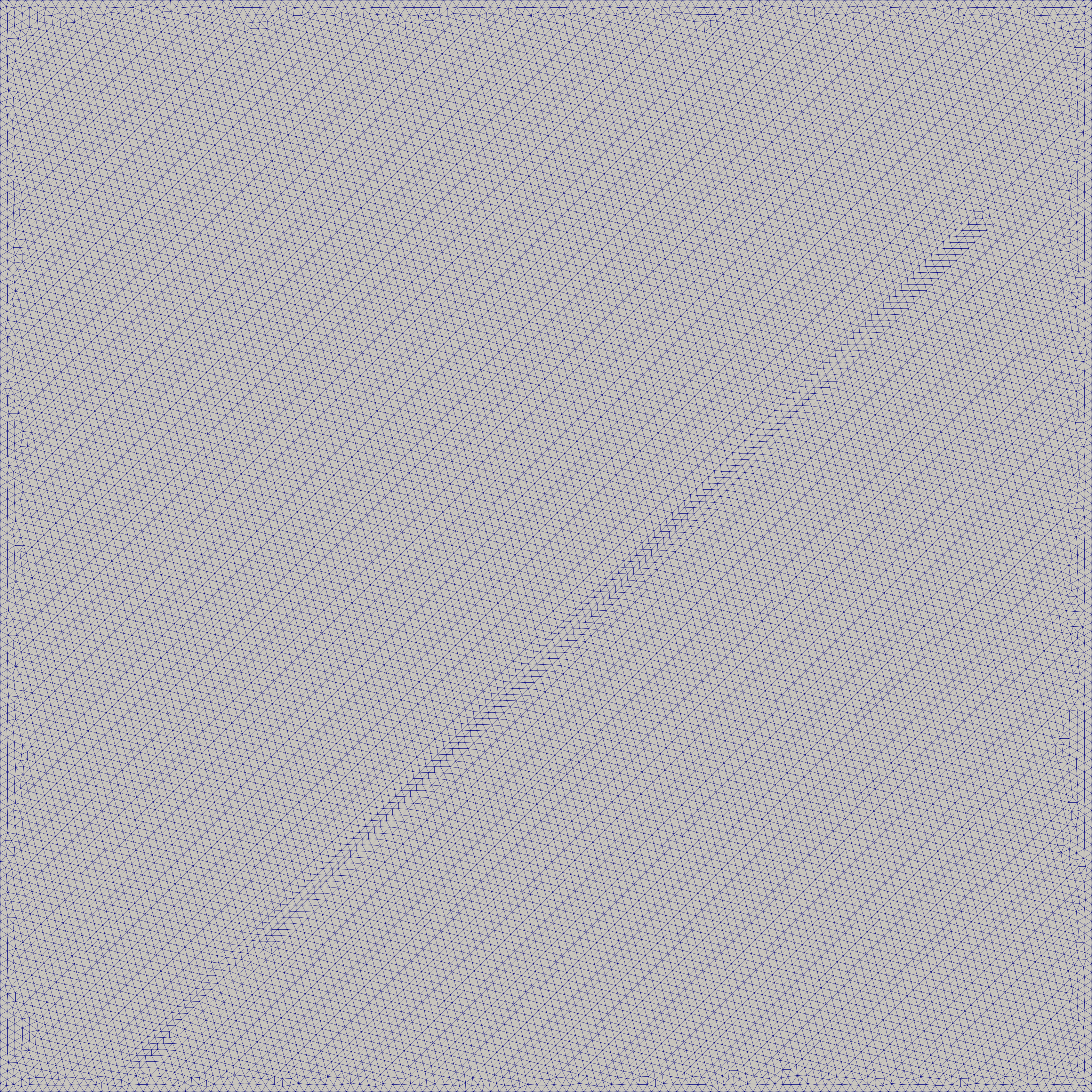}%
    \hspace*{0.01\textwidth}%
    \includegraphics[width=0.32\textwidth]{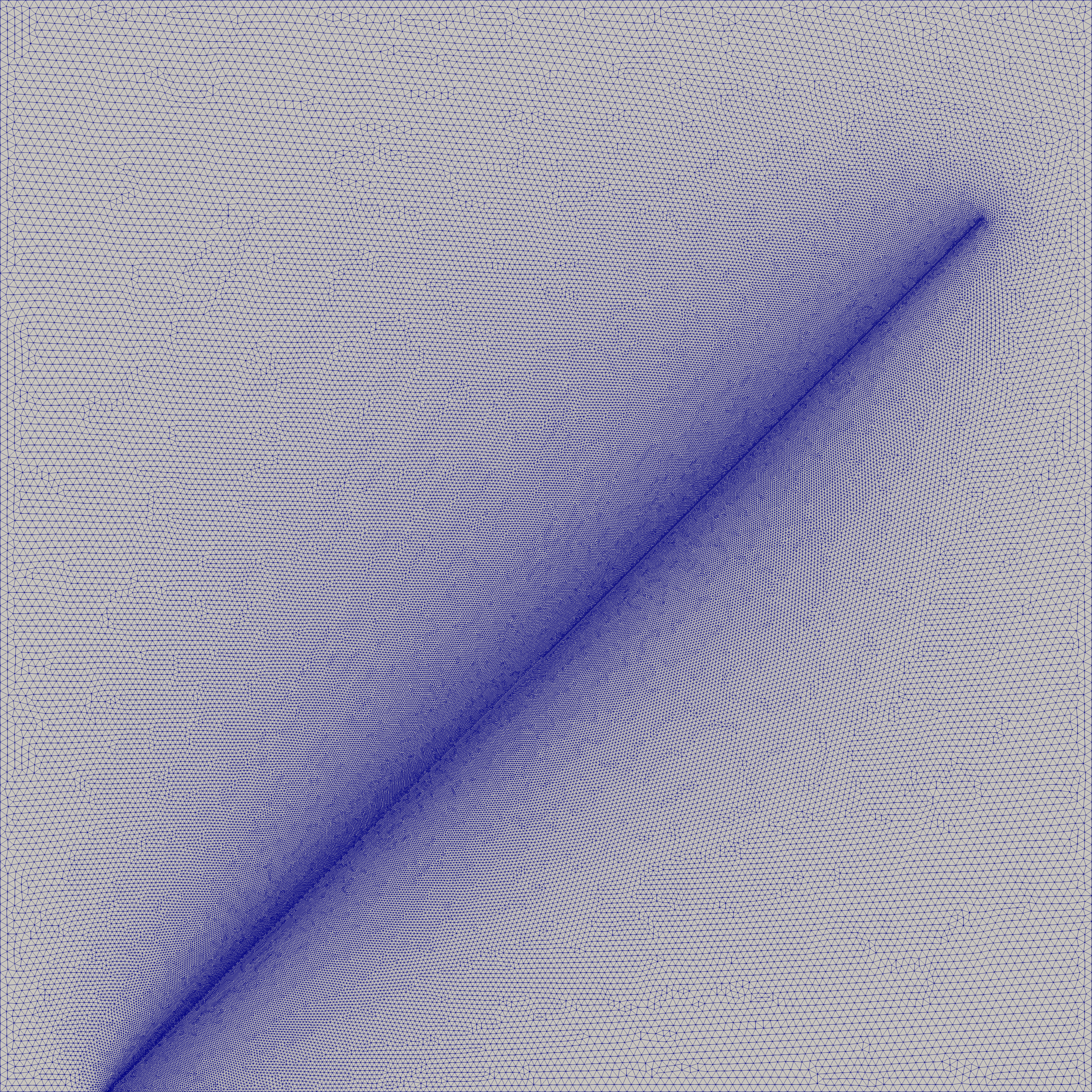}%
    \caption{On the left, domain $\Omega$, fracture $\gamma$, and the two sides of the layer
    $\mu$ for the cases of Subsection
    \ref{subsec:case1}. At the center the computational grid used for the
    multi-layer model and on the right the finer computational grid used for the simulation of the fracture-matrix model (without the layer).}%
    \label{fig:example_single_fracture_domain}
\end{figure}

In this set of tests, we consider part of the geometrical setting introduced in the
second example of \cite{Fumagalli2020e}. We refer to Table
\ref{tab:single_fracture} for a list of the data and physical parameters common
to the three cases presented in this section.
The porous medium, represented by the domain $\Omega=(0, 1)^2$, is partially
cut by a single fracture $\gamma=\{(x, y) \in \Omega:\, y = x - 0.1, x \leq
0.9\}$ with the surrounding layers $\mu$, which geometrically coincide with $
\gamma$.
See Figure \ref{fig:example_single_fracture_domain} for a graphical representation.
\begin{table}[tb]
    \centering
    \begin{tabular}{|c|c|c|c|c|c|}
        \hline
         $\delta=0.1$ & $\phi_{\Omega, 0}=0.2$ & $\phi_{\mu, 0} =0.2$ &
        $\epsilon_{\gamma, 0}=10^{-3}$ & $\epsilon_{\mu, 0} = 10^{-8}$\\ \hline
        $k_0=1$ & $k_{\gamma, 0} = 10^2$ &
        $\kappa_{\gamma, 0} = 10^2$ & $k_{\mu, 0} = 1$& $\kappa_{\mu, 0}=1$ \\ \hline
        $\mu=1$ & $f=0$ &
        $f_\gamma=0$ & $f_\mu = 0$ & $q_{\partial \Omega}^{\rm no-flow}=0$ \\\hline
        $p_{\partial \Omega}^{\rm out-flow} = 0$ & $p_{\partial\Omega}^{\rm
        in-flow}=1$
        & $p_{\partial \gamma}^{\rm in-flow}=10^{-1}$ &
        $p_{\partial \mu}^{\rm in-flow}=10^{-1}$ &
        $d=10^{-8}$
        \\\hline
        $d_\mu=10^{-6}$ &
        $\delta_\mu=10^{-6}$ &
        $d_\gamma=10^{-6}$ & $\delta_\gamma=10^{-6}$
        &$u_{\Omega, 0}=0$
        \\\hline
        $\chi_{\partial \Omega}^{\rm no-flow}=0$ &
        $u_{\partial \Omega}^{\rm in-flow}=2$ &
        $u_{\partial \Omega}^{\rm out-flow}=0$ &
        $u_{\gamma, 0}=0$ &
        $u_{\partial \gamma}^{\rm in-flow}=2$
        \\\hline
        $u_{\mu, 0}=0$
        &
        $u_{\partial \mu}^{\rm in-flow}=2$
        & $\lambda^-=100$ & $r(u) = u$ & \\\hline
    \end{tabular}
    \caption{Common data for the examples in Subsection
    \ref{subsec:case1}.}%
    \label{tab:single_fracture}
\end{table}

One of the main criteria in our evaluation, apart from a graphical observation
of the solution, is the comparison between the thickness of the layer $\mu$
estimated with model in Problem \ref{pb:multi} and the one obtained from the simulation of the
matrix-fracture Problem \ref{pb:fractured} as in \cite{Fumagalli2020e}, with a grid fine
enough to capture the concentration gradients around the fracture. This latter
high resolution simulation will numerically validate the accuracy of the proposed
model in this setting. Clearly, both test cases in this section deviate from the
assumptions at the basis of the theoretical model for the layer thickness
\eqref{growth}: the transport of solute from the fracture is not exactly
one-dimensional, there is a small diffusive effect, and, if porosity is allowed
to change due to precipitation, the Darcy velocity cannot be considered
constant. Our aim is to test the robustness of the model prediction for
different cases, to establish its usefulness in realistic situations.

The simulation has $100$ time steps of equal length, with ending time
$T_f=0.2$.  For the multi-layer model we consider a uniformly refined mesh of
38435 triangles for the porous media, 290 segments for the layer and 145
segments for the fracture, while for the model where only the fracture is a
lower dimensional object, we have considered a very fine grid around the
fracture itself which gives a non-uniform triangular grid composed of 107841
elements. The fracture is discretized with 906 equal segments. See Figure
\ref{fig:example_single_fracture_domain} for the graphical representation of the
computational grids.

\subsubsection{Case 1}\label{subsec:case1a}

For this case we consider the data and geometry describe above, and we
additionally set the following parameters: $\eta_\Omega = 0$, $\eta_\gamma = 0$,
and $\eta_\mu = 0$ thus the porosity $\phi_\Omega$ and $\phi_\mu$, as well as
the fracture aperture $\epsilon_\gamma$, are fixed for the entire simulation and
are equal to their initial value. In this case the Darcy velocity $\bm{q}$
is constant in time in the entire domain (although not necessarily exactly normal
to the fracture).

\begin{figure}[tbp]
    \centering
    \includegraphics[width=0.375\textwidth]{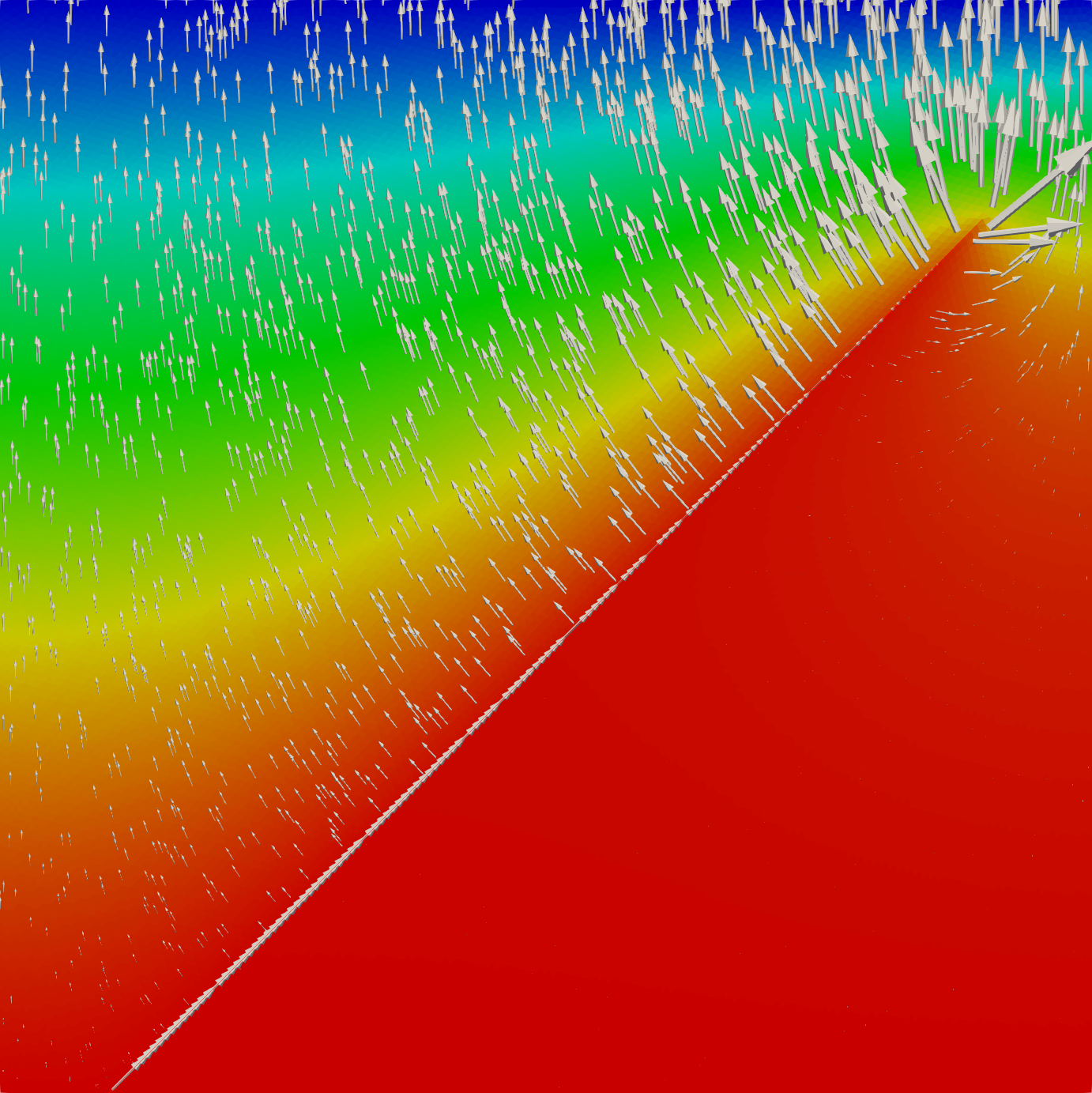}%
    \hspace*{0.05\textwidth}%
    \raisebox{0.2\height}{\includegraphics[width=0.05\textwidth]{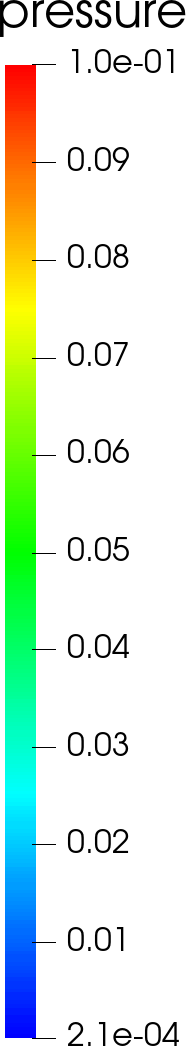}}%
    \hspace*{0.05\textwidth}%
    \includegraphics[width=0.375\textwidth]{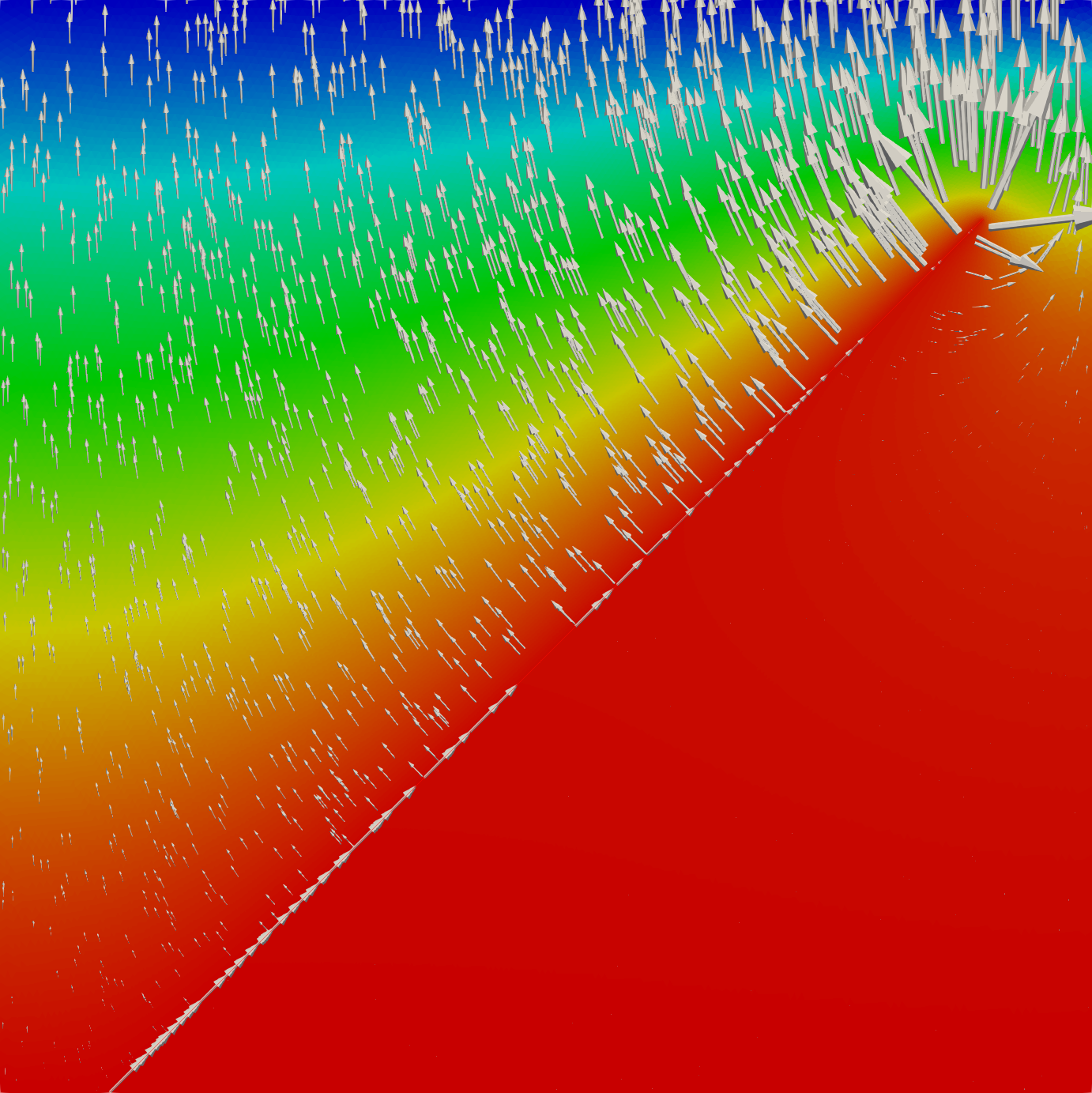}\\
    \includegraphics[width=0.375\textwidth]{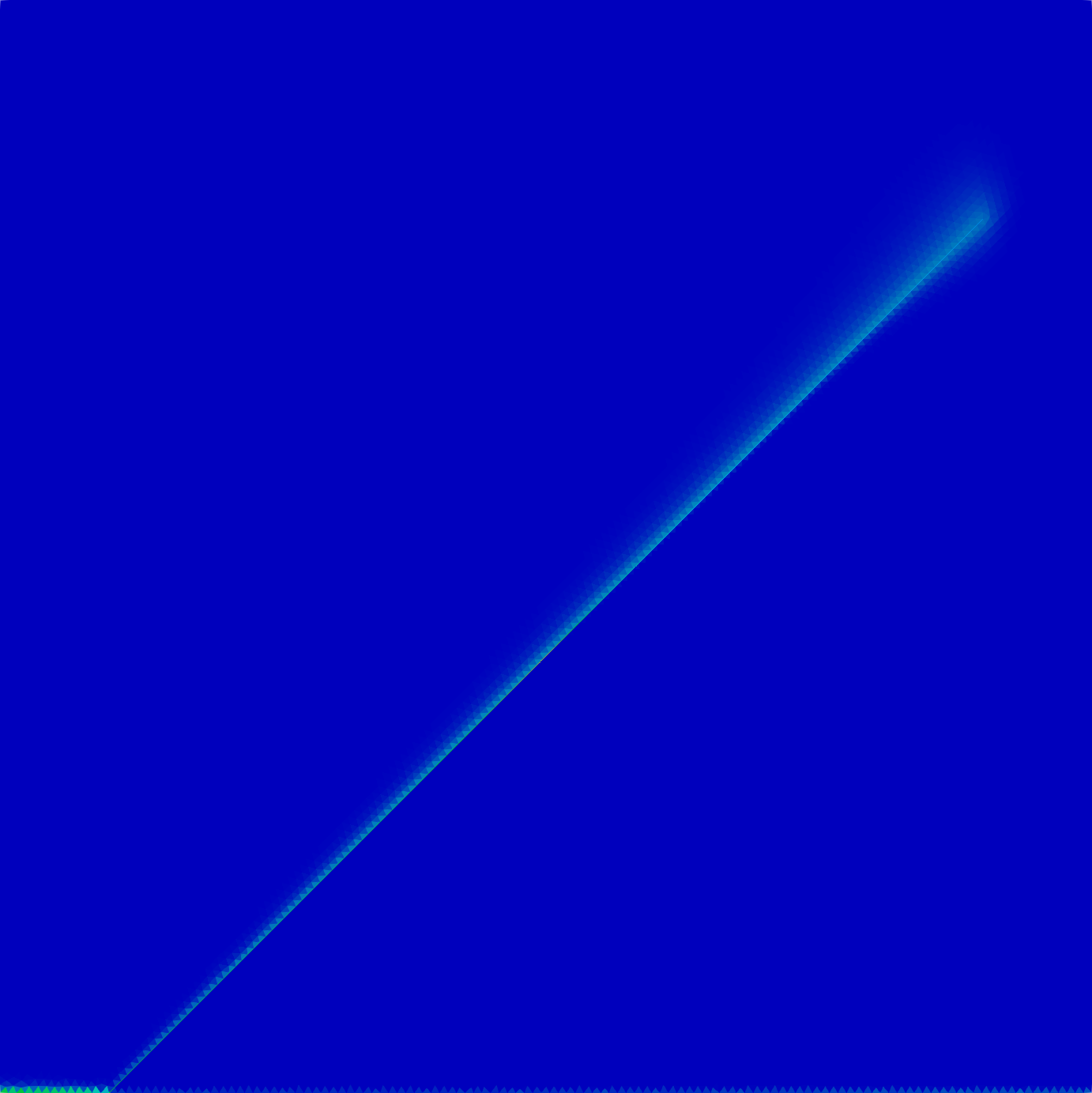}%
    \hspace*{0.05\textwidth}%
    \raisebox{0.2\height}{\includegraphics[width=0.05\textwidth]{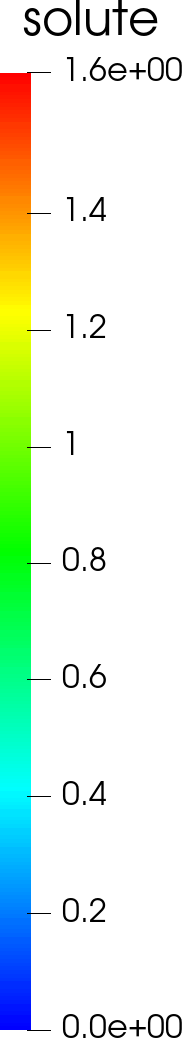}}%
    \hspace*{0.05\textwidth}%
    \includegraphics[width=0.375\textwidth]{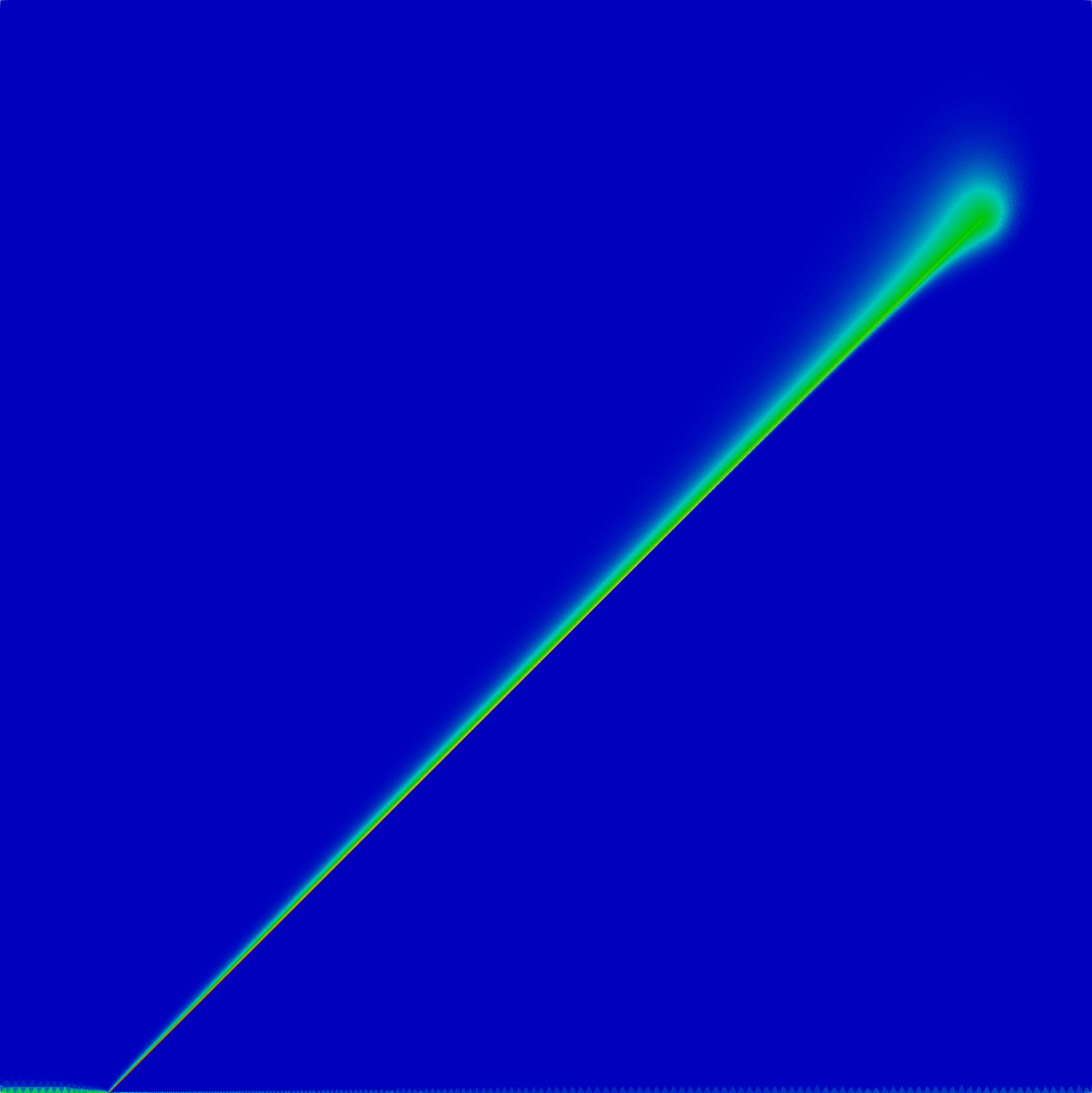}\\
    \includegraphics[width=0.375\textwidth]{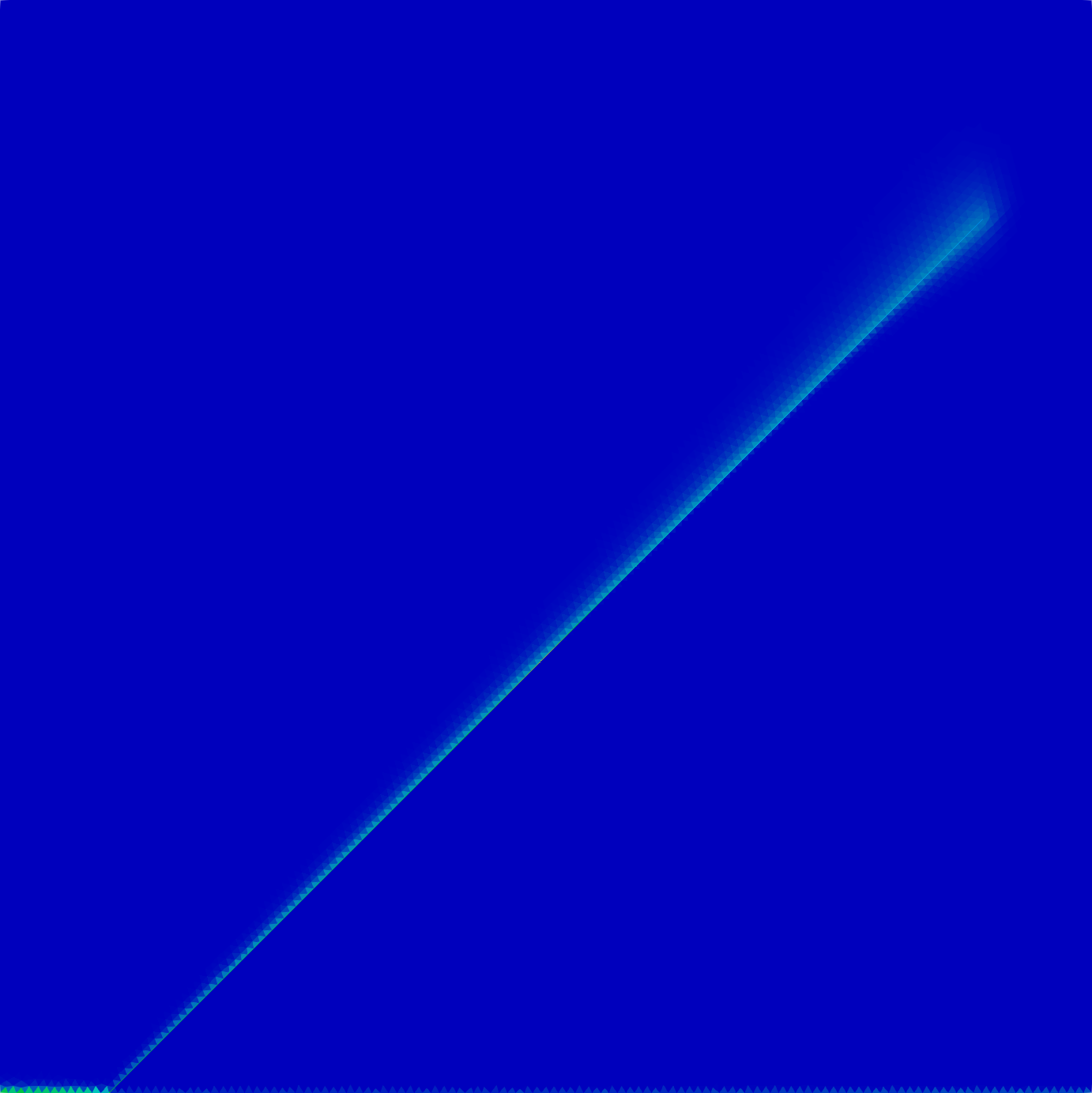}%
    \hspace*{0.05\textwidth}%
    \raisebox{0.2\height}{\includegraphics[width=0.05\textwidth]{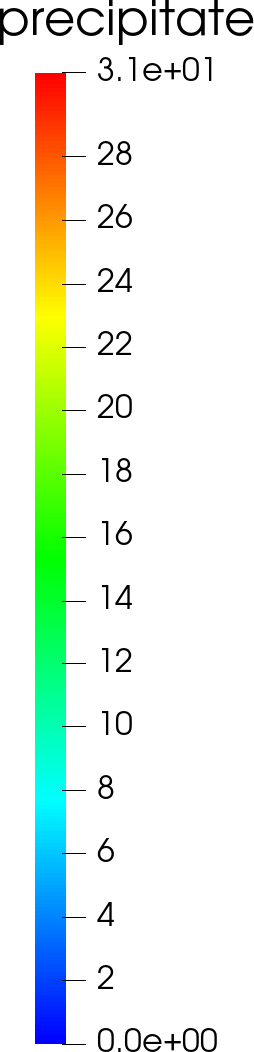}}%
    \hspace*{0.05\textwidth}%
    \includegraphics[width=0.375\textwidth]{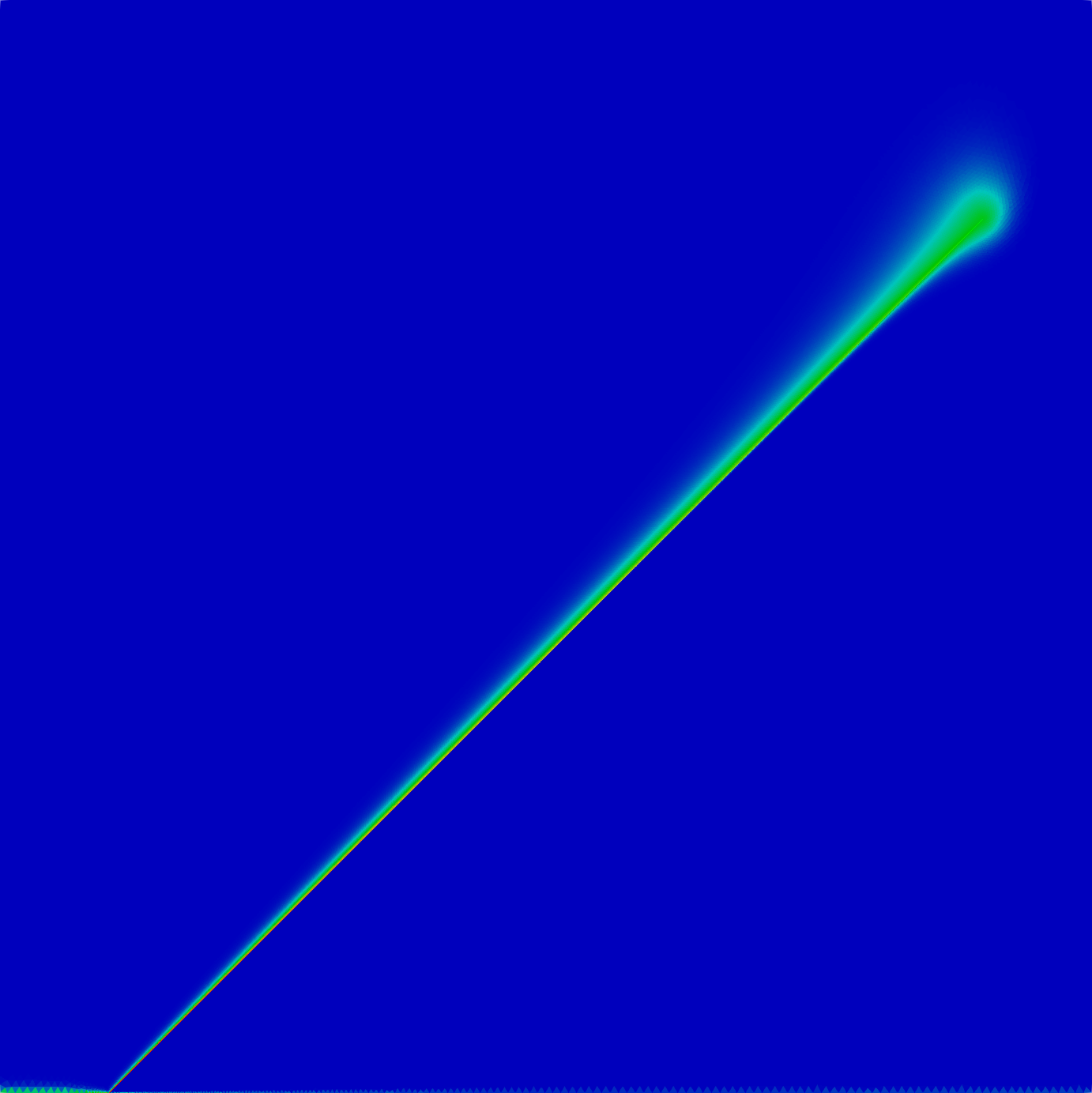}
    \caption{Graphical representation of the pressure (top), solute (centre), and
    precipitate (bottom) for the test case described in Sub-subsection
    \ref{subsec:case1a} at the end of the simulation.
    On the left for the
    multi-layer reduced model in Problem \ref{pb:multi} and on the right for the
    matrix-fracture model in Problem \ref{pb:fractured}.}
    \label{fig:2d_case1}
\end{figure}
In Figure \ref{fig:2d_case1} we compare the pressure and Darcy velocity, along
with the solute and precipitate obtained with the two models and corresponding
discretizations. First of all let us note that, since $\eta_\Omega$ is null,
pressure and Darcy velocity are fixed in time The solution shows the advantage
of adopting the introduced model, since we can observe that with a fast enough
reaction most of the dynamics for the solute and precipitate develops very close
to the fracture $\gamma$.  The graphical difference in the solute and
precipitate distribution is mainly due to the fact that, in the multi-layer
model, the part of the solution with the higher concentrations is represented by
the reduced variables $u_\mu$, $w_\mu$ in the one-dimensional layers.
\begin{figure}[tb]
    \centering
    \includegraphics[width=0.375\textwidth]{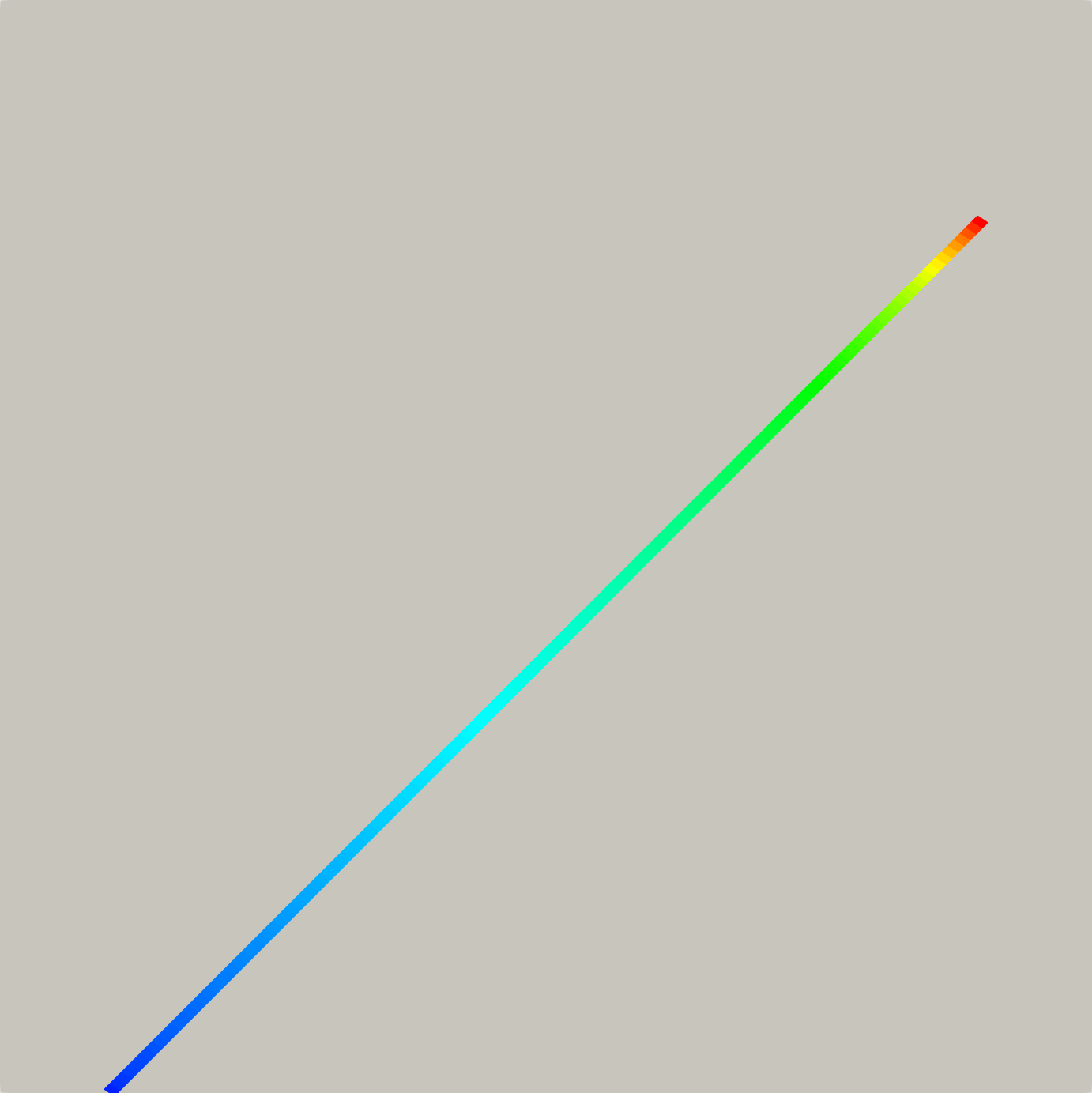}%
    \hspace*{0.05\textwidth}%
    \raisebox{0.1\height}{\includegraphics[width=0.05\textwidth]{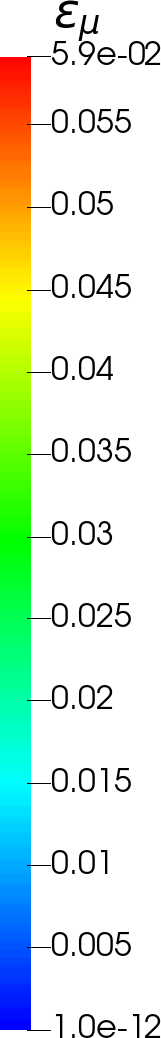}}%
    \hspace*{0.05\textwidth}%
    \includegraphics[width=0.375\textwidth]{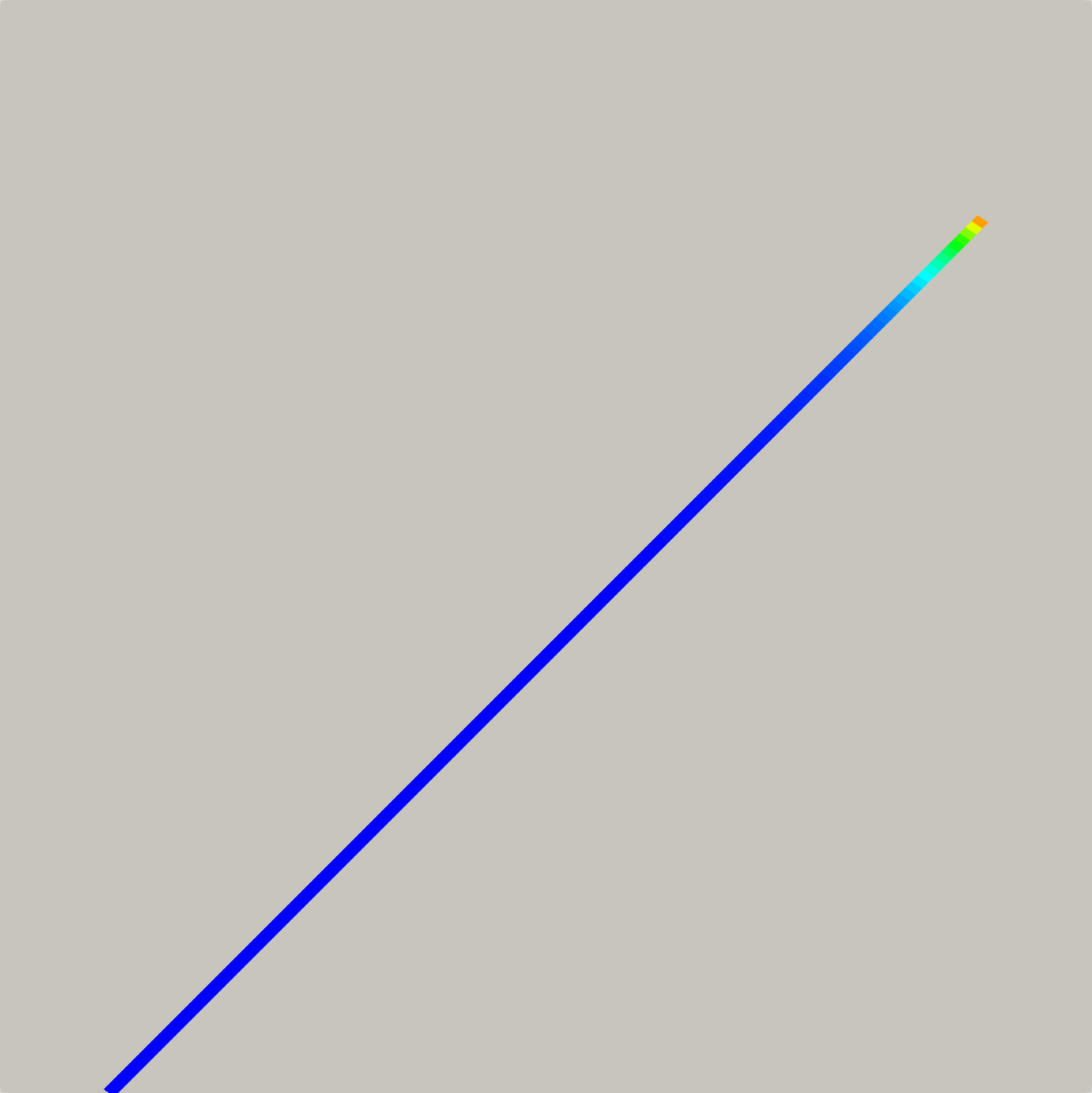}
    \caption{Layer thickness for both sides (top layer on the left, bottom layer on the right) at the end of the simulation time
    for the test case in Sub-subsection \ref{subsec:case1a}.}
    \label{fig:2d_case1_layer}
\end{figure}
In Figure \ref{fig:2d_case1_layer} we see the layer thickness at the end of the
simulation. Since it depends on the Darcy velocity at fracture-layer interface,
the top part of the layer (the one closer to the outflow) is wider than on the
bottom part. For both sides, at their tip the aperture results in a much higher
value due to the outflow from the fracture tip. Given the small layer
thickness, we can consider the proposed model to be in its range of
applicability, i.e. the layers can be reduced to their center line.
\begin{figure}[tb]
    \centering
    \includegraphics[width=0.375\textwidth]{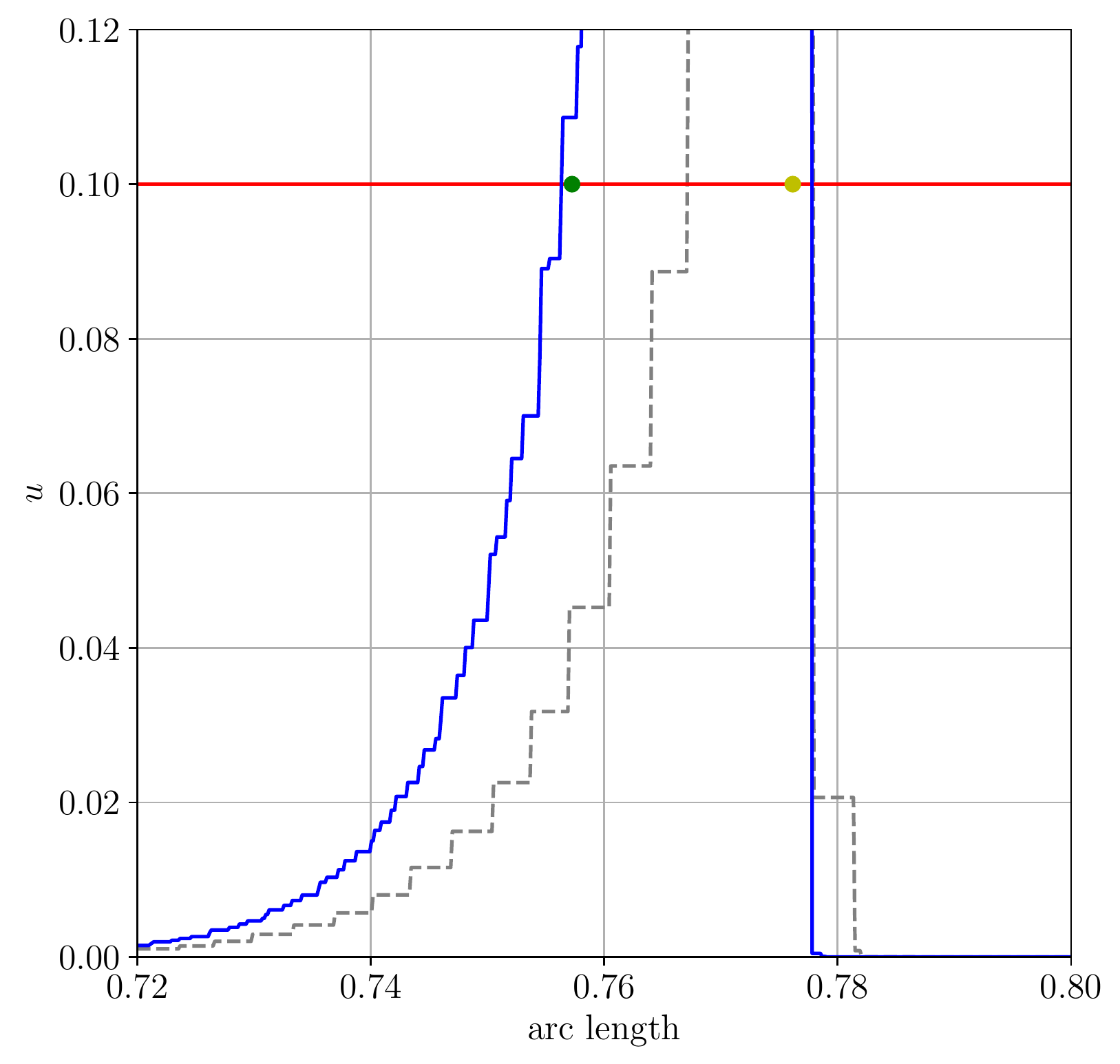}%
    \hspace*{0.05\textwidth}%
    \includegraphics[width=0.375\textwidth]{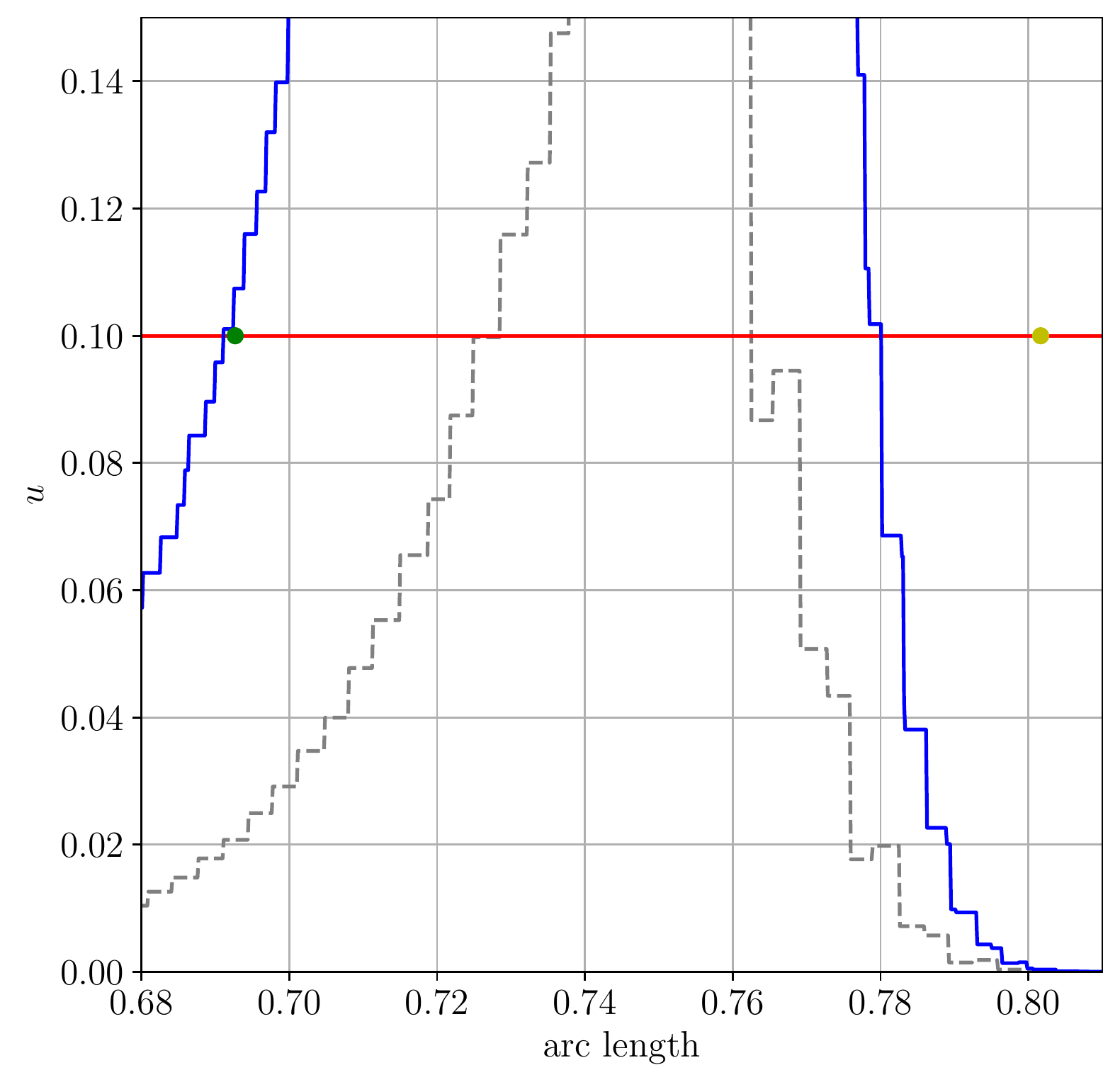}
    \caption{Profile of the solute along two lines normal to the fracture, $l_1$
    on the left and $l_2$ on the right. The blue line represent the solute
    profile obtained with the proposed multi-layer model, while the red line with
    the equi-dimensional layers $\mu$ and a fine grid. The green and yellow dot
    represent the point  $(a^\pm\pm\epsilon_\mu^\pm, \delta)$ for the
    top and bottom layer respectively. Results for the test case in
    Sub-subsection \ref{subsec:case1a}.}
    \label{fig:2d_case1_profile}
\end{figure}
Finally, Figure \ref{fig:2d_case1_profile} shows the graphical comparison
between the solute along two specifics lines normal to the fracture: $l_1$ which
connects $(0, 1)$ and $(1, 0)$, crossing the fracture at $(0.5,0.5)$ and layers,
and $l_2$ which connects $(0.3656, 1.3293)$ and $(1.3658, 0.3293)$, passing
close to the fracture tip.  The concentration profiles computed with the
matrix-fracture and the multi-layer models are plotted against the arc curve
length coordinate along $l_1$ and $l_2$. These profiles are compared with the layer
thickness predicted by \eqref{growth}, marked by the dots in position
$(a_{1,2}\pm\epsilon_\mu^\pm, \delta)$, where $a_1=0.778$ (for line $l_1$) and
$a_2=0.7521$ (line $l_2$) are the intersections with the fracture. These dots correspond to the point where solute concentration drops below the cutoff value ($\delta=0.1$), and thus mark the border of the layer $\mu^\pm$.

We can notice that for $l_1$ the results are in good agreement, while for $l_2$
we get accurate results only for the top part of the layer. For the bottom part
of $\mu$ on its tip, the model assumption that the flow is mostly normal
however, in this particular case, is not valid since the outflow from the
fracture tip creates strongly bidimensional effects in the solution.

We can conclude that, in this setting, the multi-layer reduced model is an
attractive and effective alternative which gives coherent results with the model
of \cite{Fumagalli2020e}.

\subsubsection{Case 2}\label{subsec:case1b}

In this second test case of the group, we allow for a more complex physical
interaction between the variables by setting $\eta_\Omega = 5\cdot10^{-2}$,
$\eta_\gamma = 5\cdot10^{-2}$, and $\eta_\mu =5\cdot10^{-2} $: thus, the porosity
and consequently then the permeability change in time and alter pressure and Darcy
velocity fields. The problem becomes more coupled. We would like to understand
if the presented model gives reasonably accurate outcomes even if this setting
does not satisfy the hypotheses at the basis of the derived layer evolution (even more than the
previous case).

\begin{figure}[tbp]
    \centering
    \includegraphics[width=0.375\textwidth]{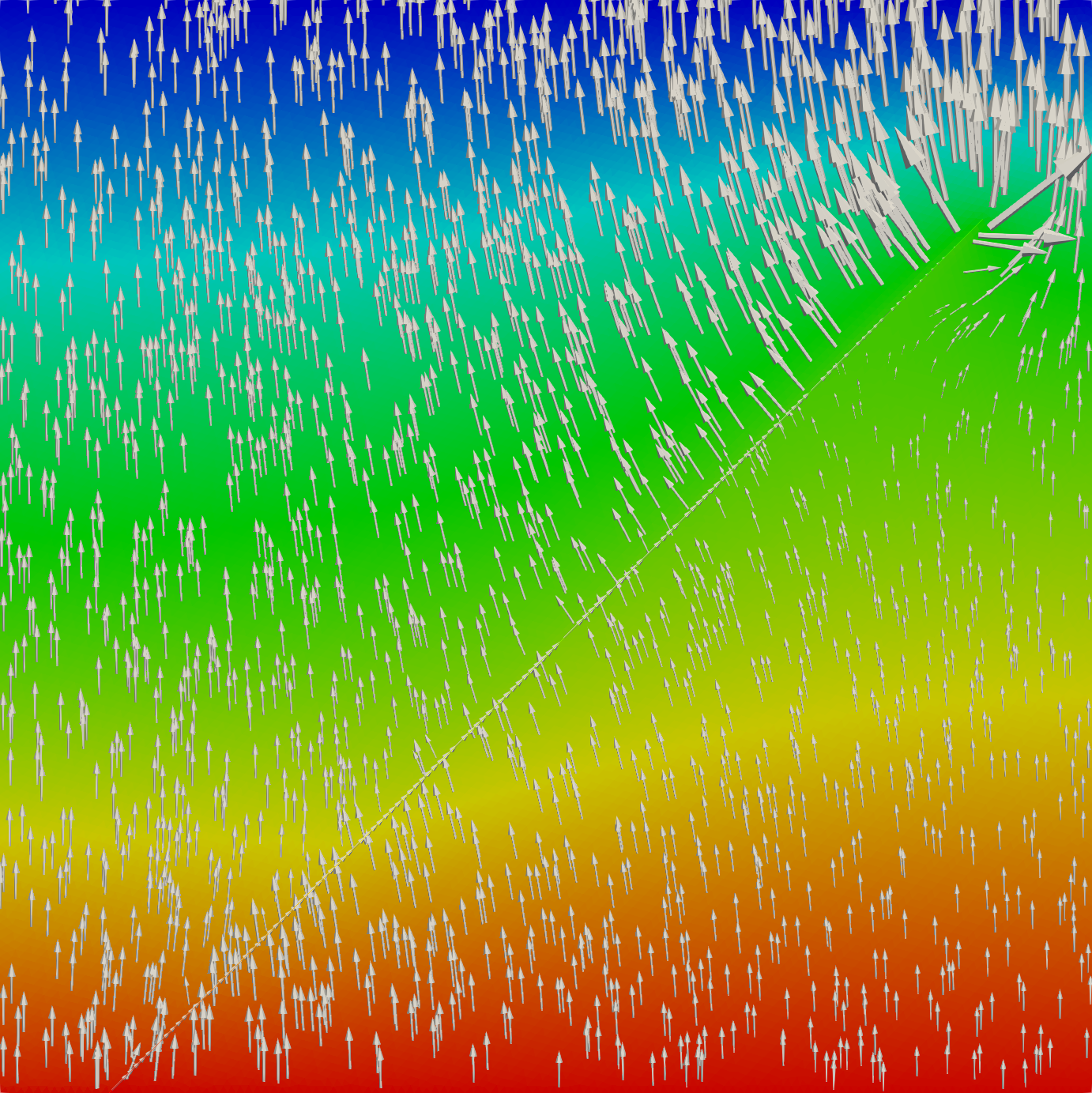}%
    \hspace*{0.05\textwidth}%
    \raisebox{0.2\height}{\includegraphics[width=0.05\textwidth]{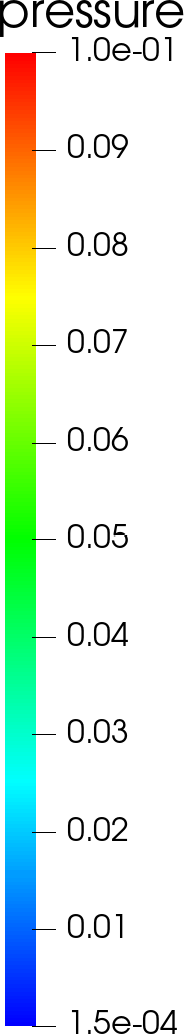}}%
    \hspace*{0.05\textwidth}%
    \includegraphics[width=0.375\textwidth]{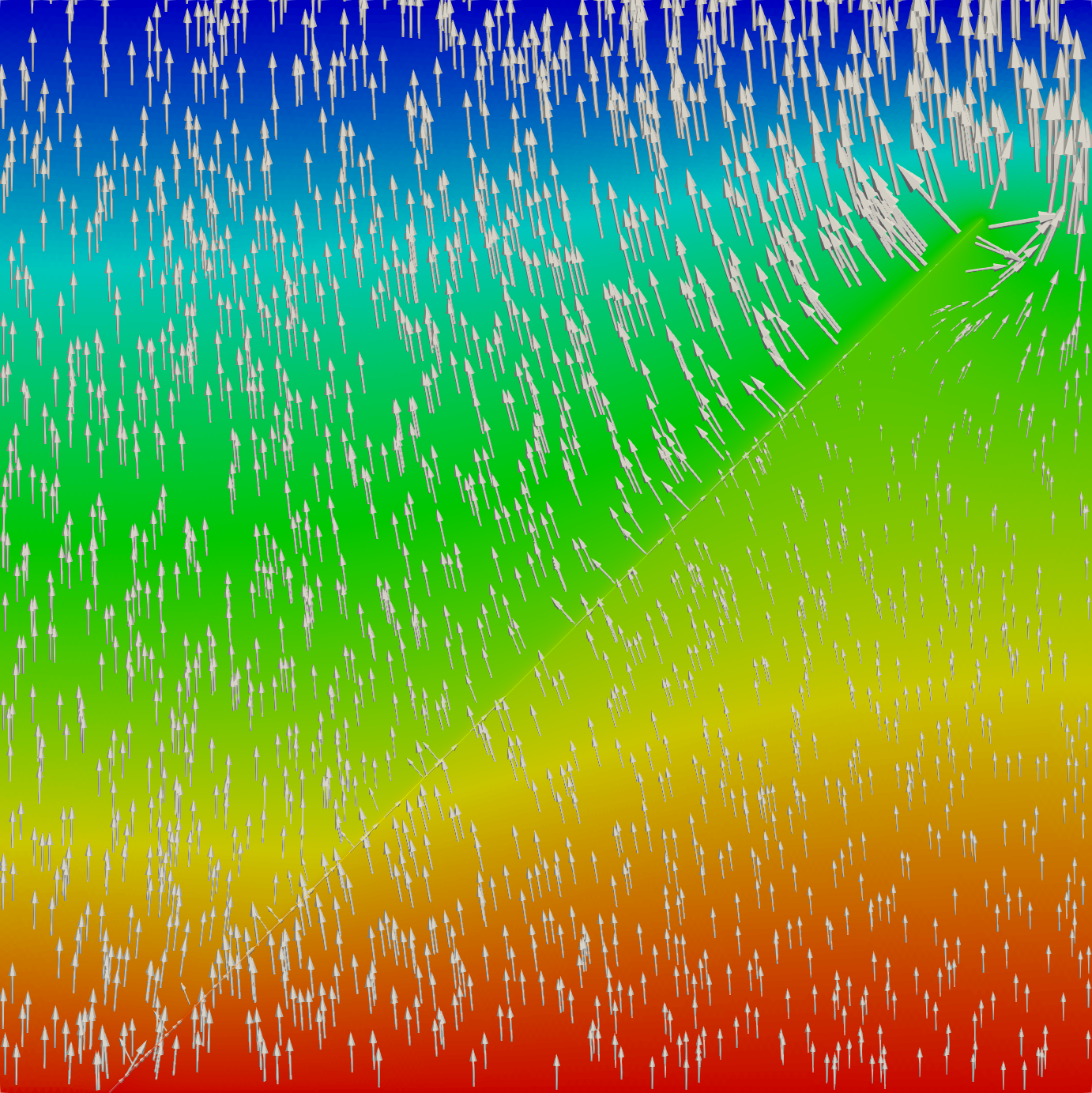}\\
    \includegraphics[width=0.375\textwidth]{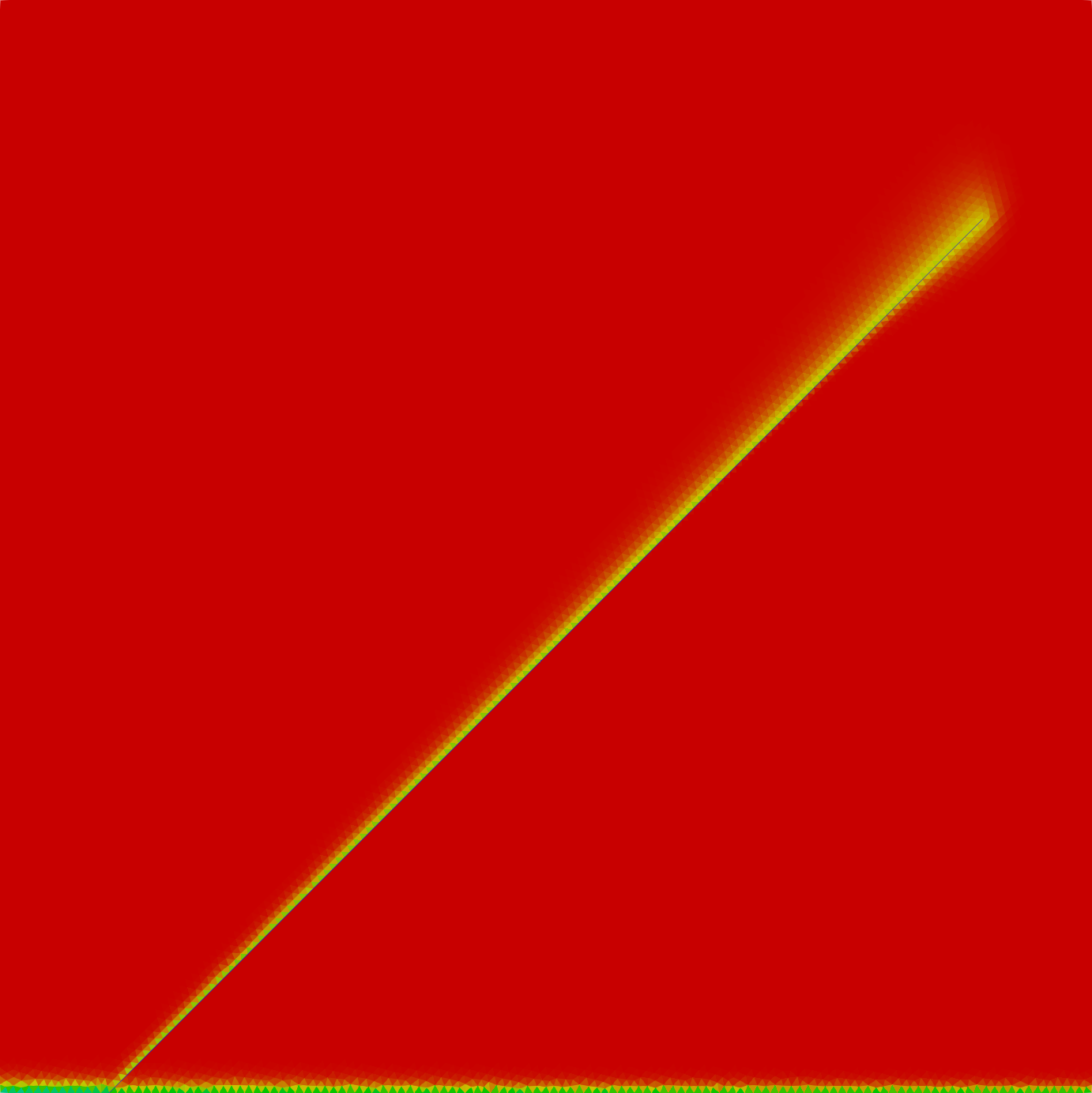}%
    \hspace*{0.05\textwidth}%
    \raisebox{0.2\height}{\includegraphics[width=0.05\textwidth]{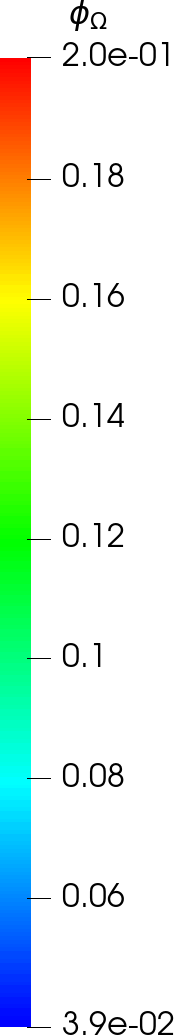}}%
    \hspace*{0.05\textwidth}%
    \includegraphics[width=0.375\textwidth]{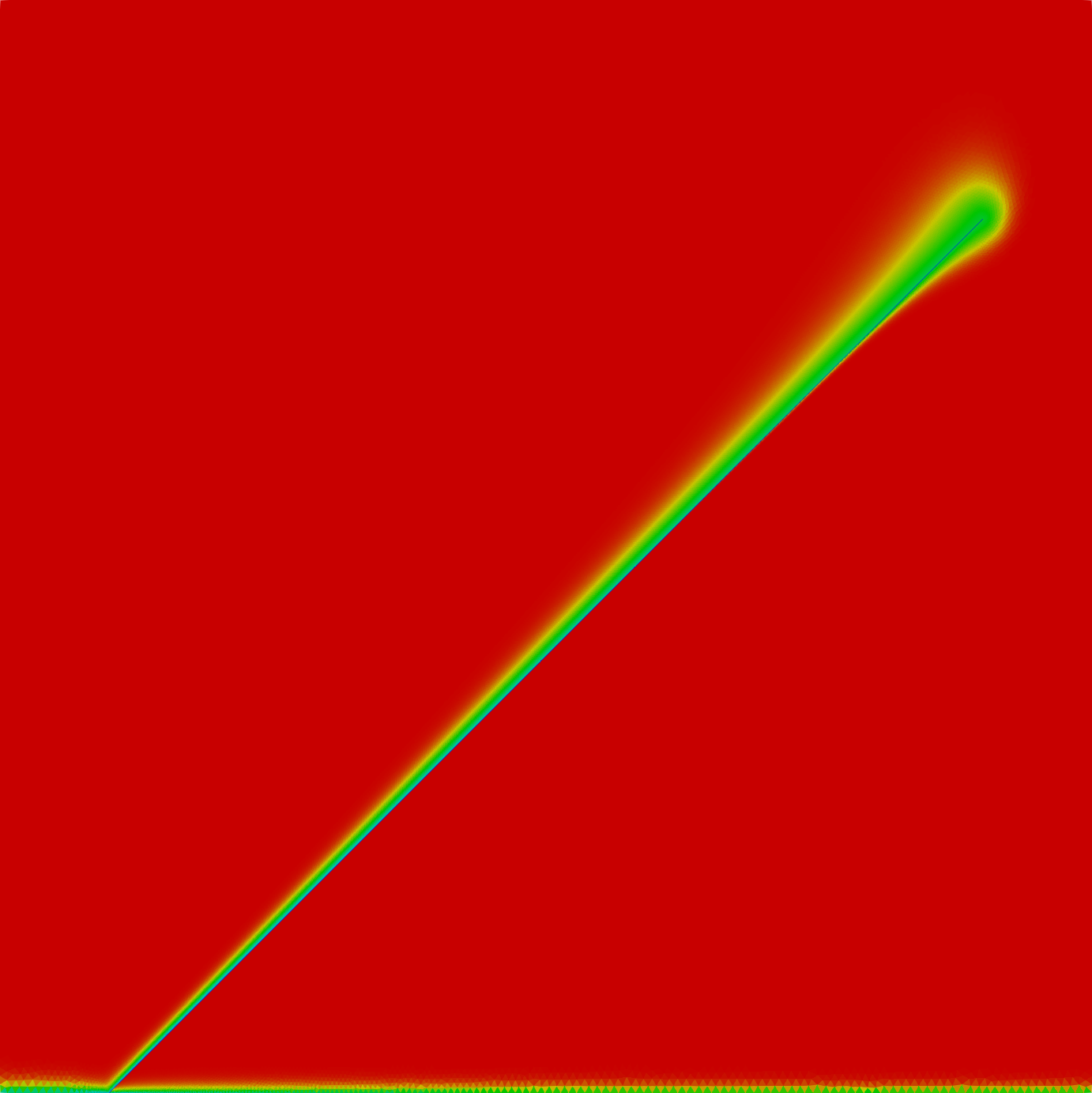}\\
    \includegraphics[width=0.375\textwidth]{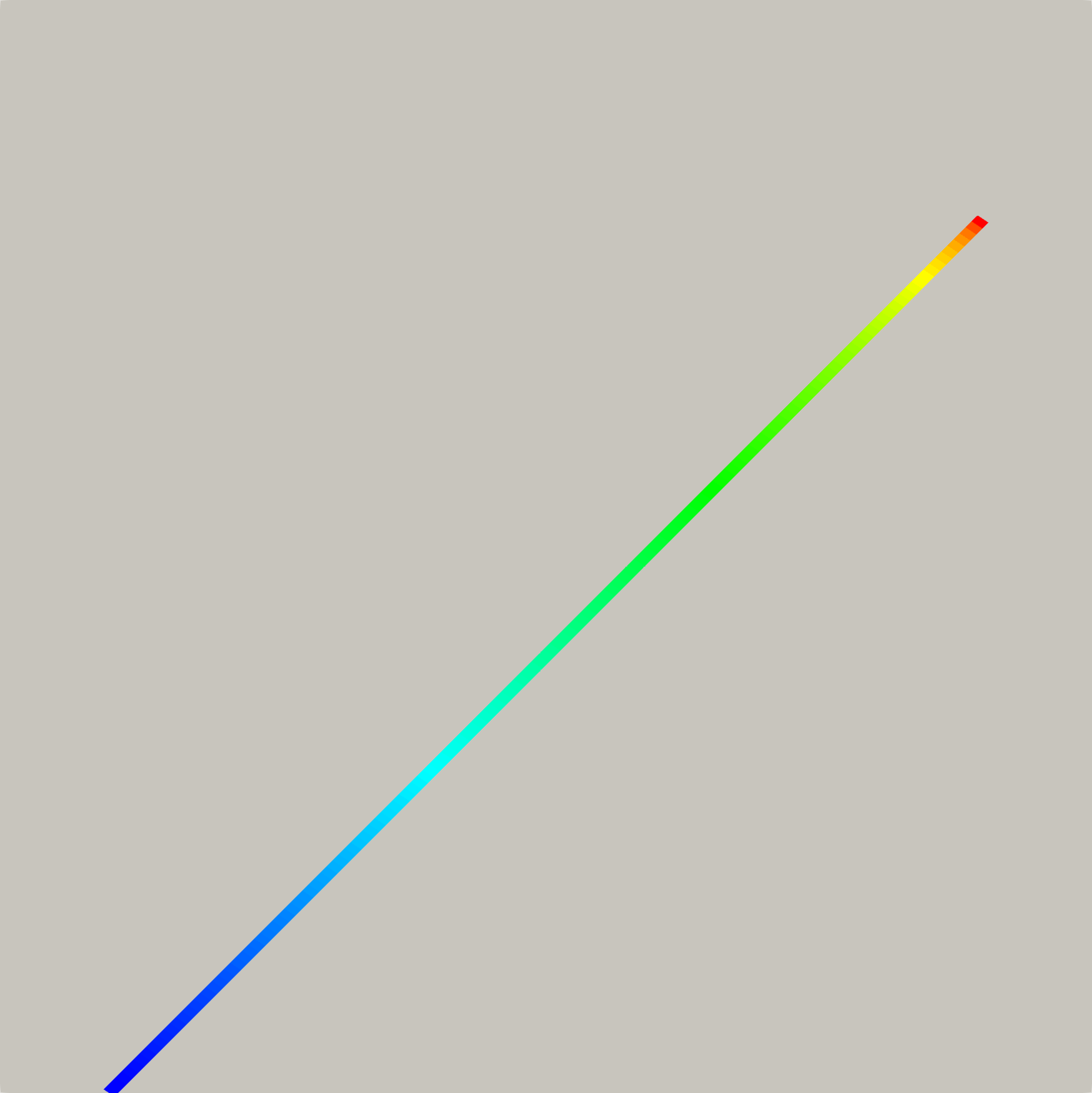}%
    \hspace*{0.05\textwidth}%
    \raisebox{0.2\height}{\includegraphics[width=0.05\textwidth]{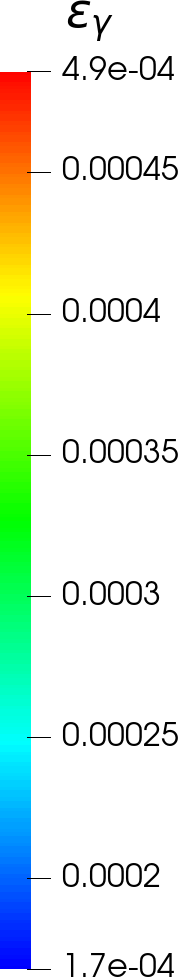}}%
    \hspace*{0.05\textwidth}%
    \includegraphics[width=0.375\textwidth]{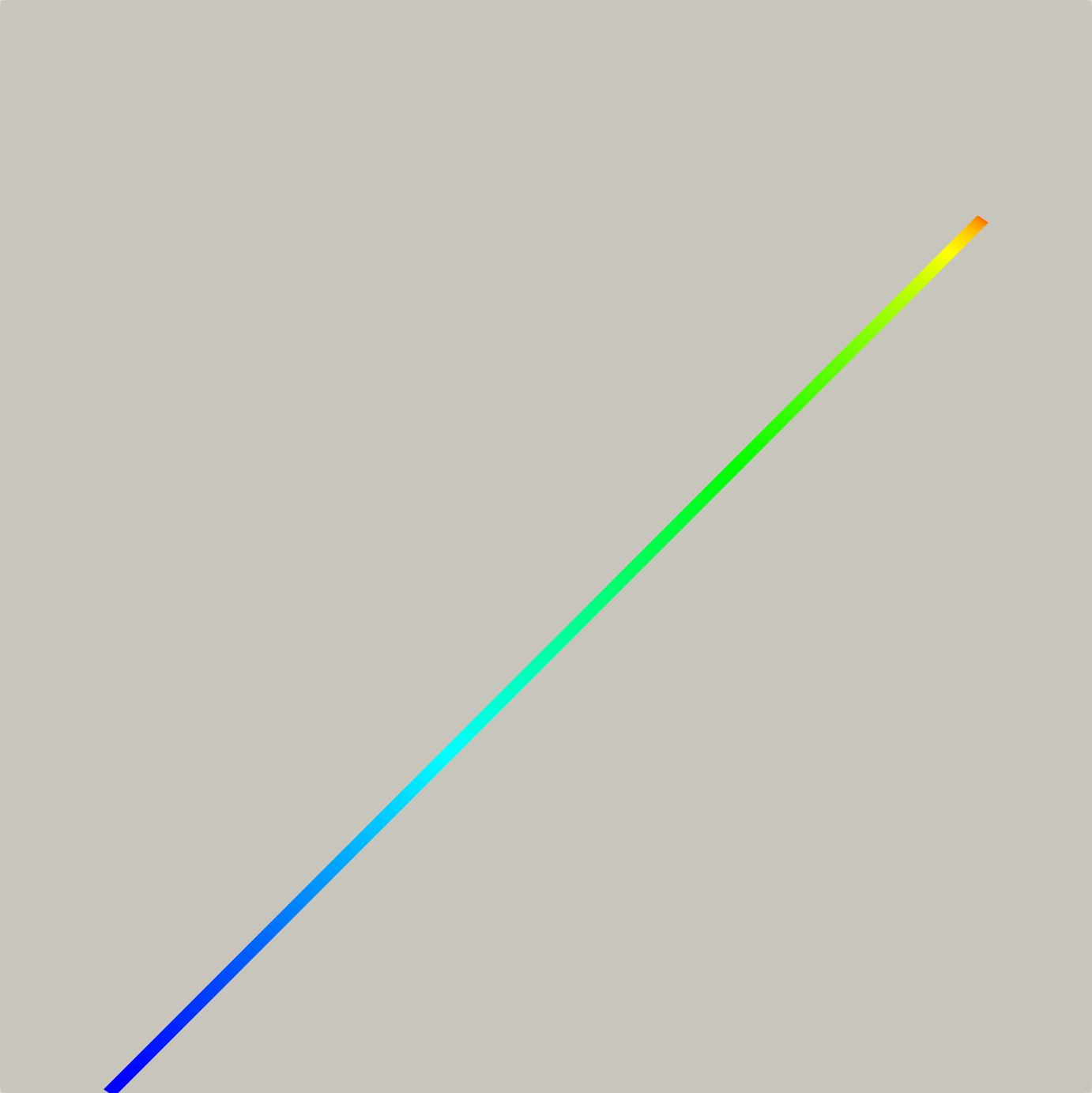}
    \caption{Graphical representation of the pressure (top), porosity (centre), and
    fracture aperture (bottom) for the test case described in Sub-subsection
    \ref{subsec:case1b} at the end of the simulation. On the left for the
    multi-layer reduced model in Problem \ref{pb:multi} and on the right for the
    matrix-fracture model in Problem \ref{pb:fractured}.}
    \label{fig:2d_case1b}
\end{figure}

Let us consider a graphical comparison of the solution obtained with the two
models, reported in Figure \ref{fig:2d_case1b}.
Comparing the pressure profile with the one obtained in the test case of
Sub-subsection \ref{subsec:case1a}, we clearly see the effect of the $\eta$
parameters on porosity and fracture aperture. The fracture indeed now becomes
less permeable (due to its shrinking aperture) as well as the layer surrounding
it (due to decreasing porosity).

We note that the predicted and observed thickness of the
layer is such that, also
in this case, it is beneficial to adopt a multi-layer approach.
We notice also that the fracture aperture is smaller closer to the inflow of the
problem: this is due to the solute that enters the domain, flows mainly into
the fracture and precipitates there, altering its aperture. This results in a
slower fracture flow which in turn affects the overall process.
\begin{figure}[tb]
    \centering
    \includegraphics[width=0.375\textwidth]{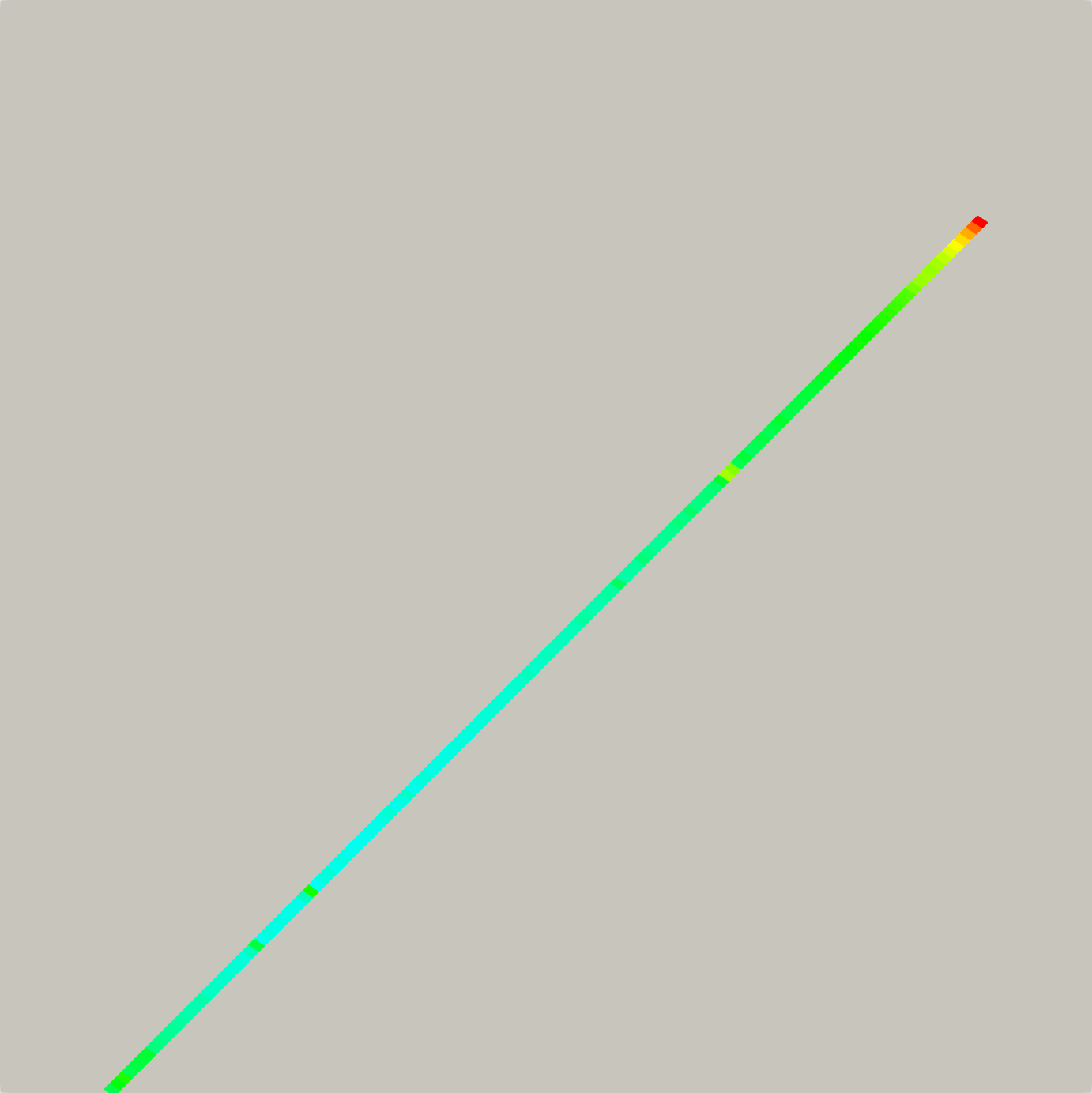}%
    \hspace*{0.05\textwidth}%
    \raisebox{0.1\height}{\includegraphics[width=0.05\textwidth]{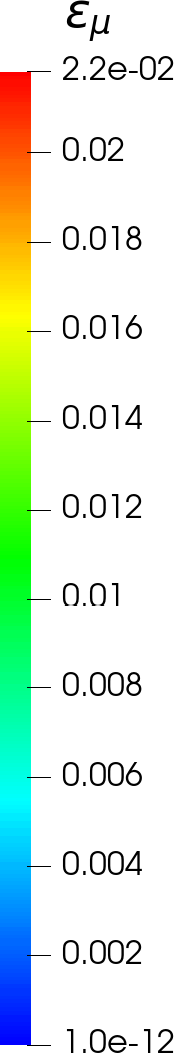}}%
    \hspace*{0.05\textwidth}%
    \includegraphics[width=0.375\textwidth]{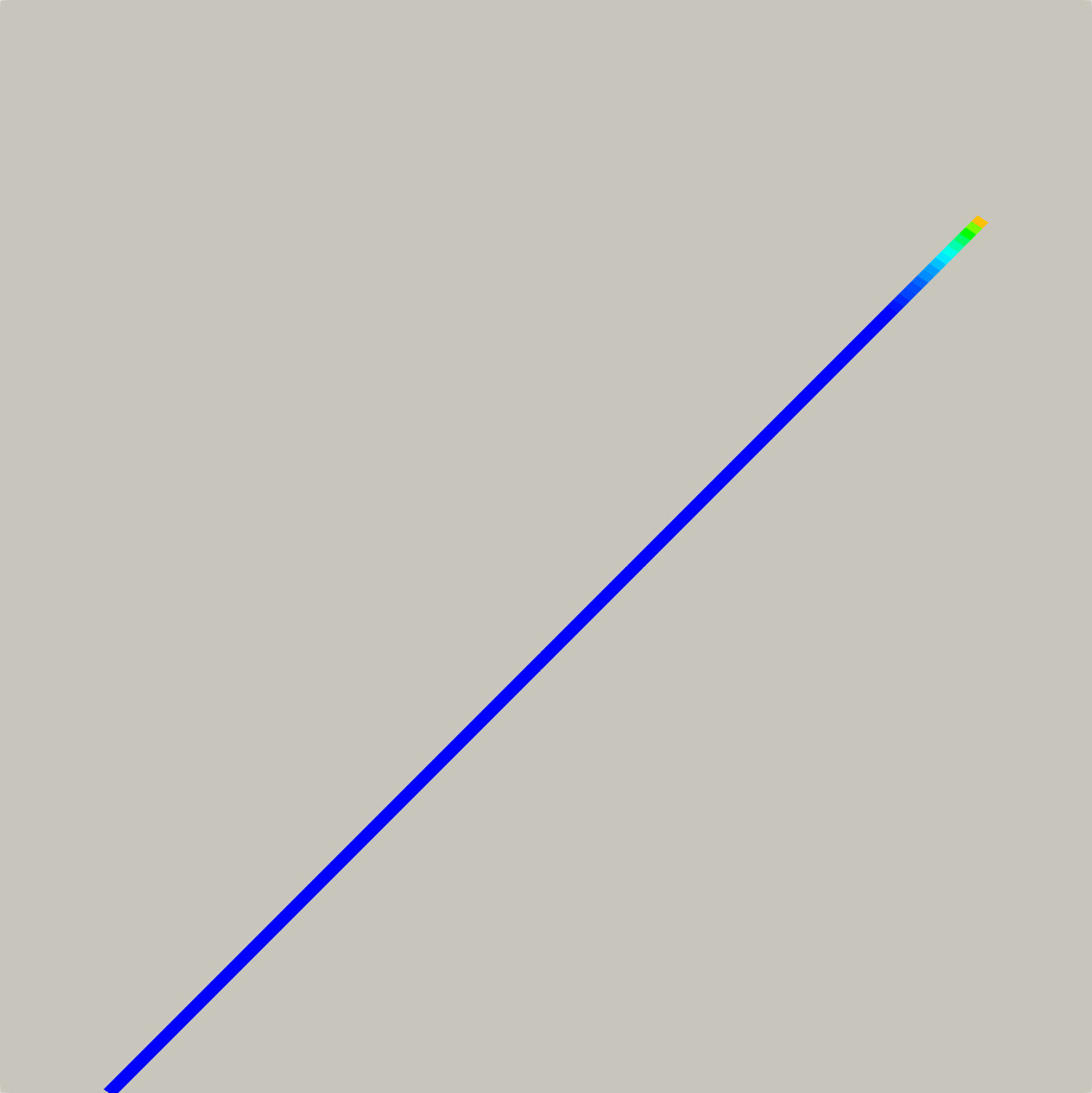}\\
    \includegraphics[width=0.375\textwidth]{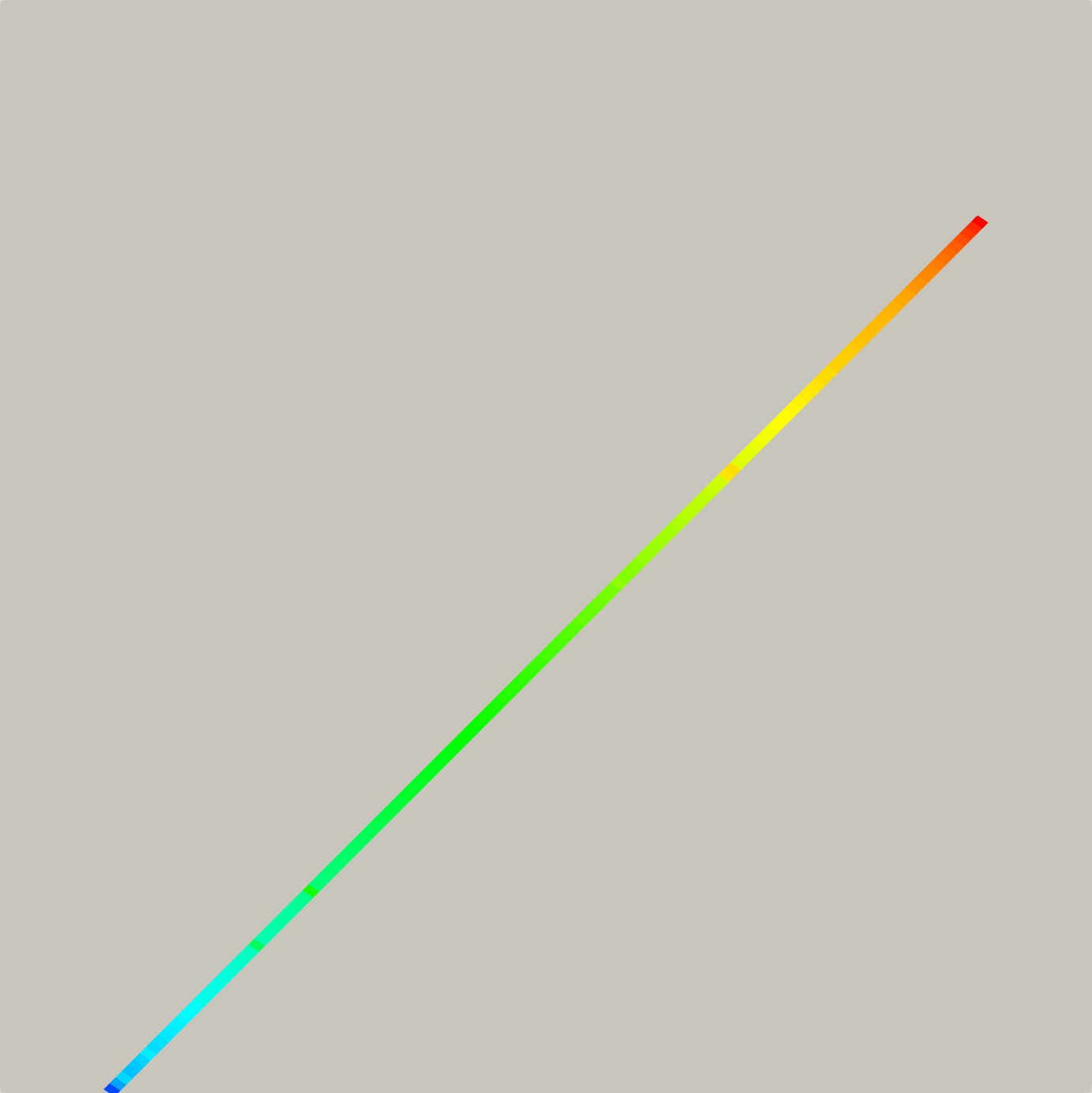}%
    \hspace*{0.05\textwidth}%
    \raisebox{0.1\height}{\includegraphics[width=0.05\textwidth]{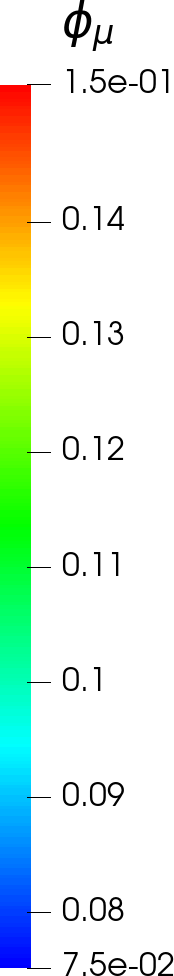}}%
    \hspace*{0.05\textwidth}%
    \includegraphics[width=0.375\textwidth]{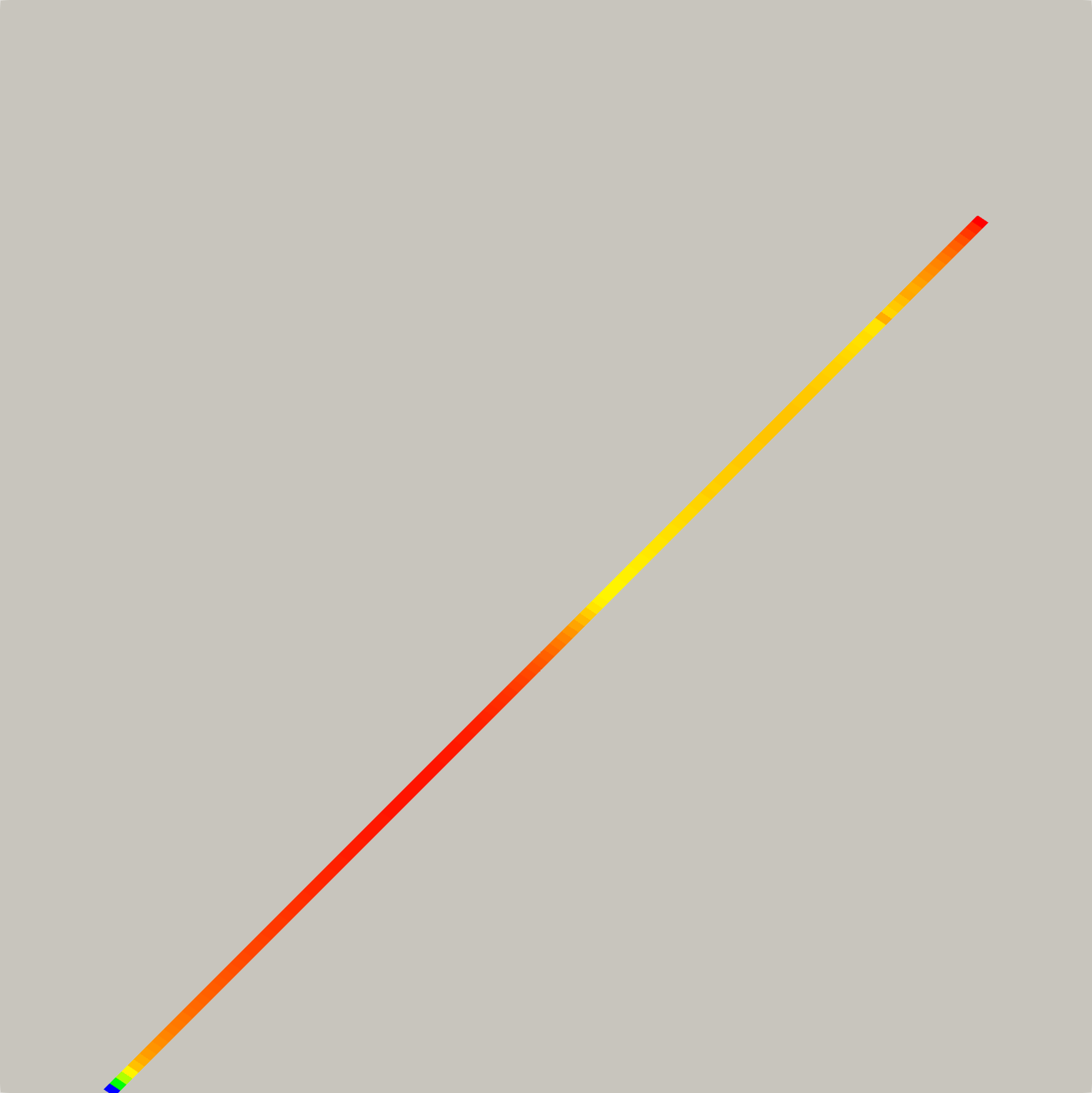}
    \caption{Layer thickness for both sides (top layer on the left, bottom layer on the right) at the end of the simulation time for the test case in Sub-subsection \ref{subsec:case1b}.}
    \label{fig:2d_case1b_layer}
\end{figure}
Figure \ref{fig:2d_case1b_layer} represents the layer thickness and porosity at
the end of the simulation. The difference between the two sides is evident,
mainly due to the difference in the flow exiting the fracture on the two sides.
\begin{figure}[tb]
    \centering
    \includegraphics[width=0.375\textwidth]{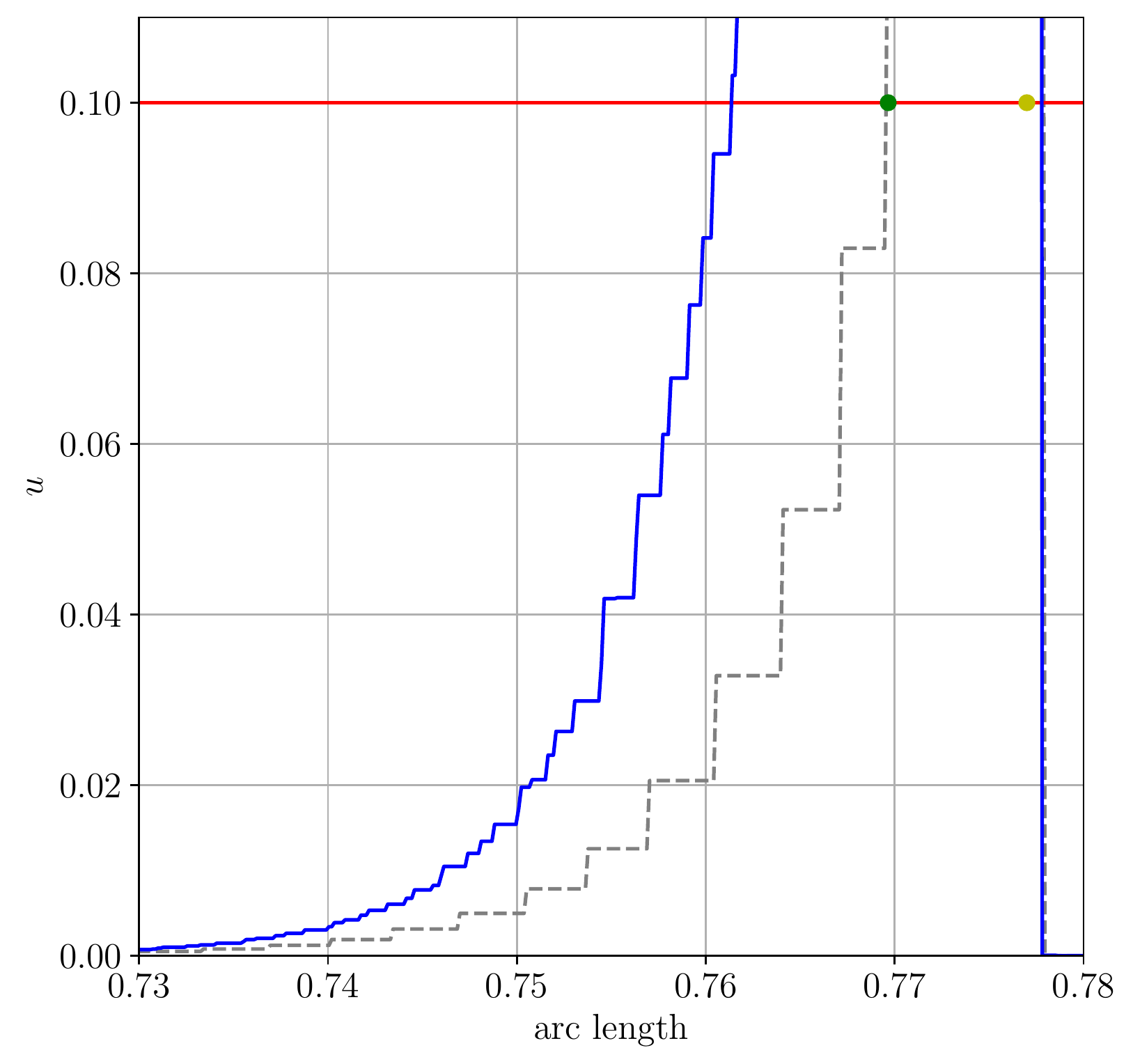}%
    \hspace*{0.05\textwidth}%
    \includegraphics[width=0.375\textwidth]{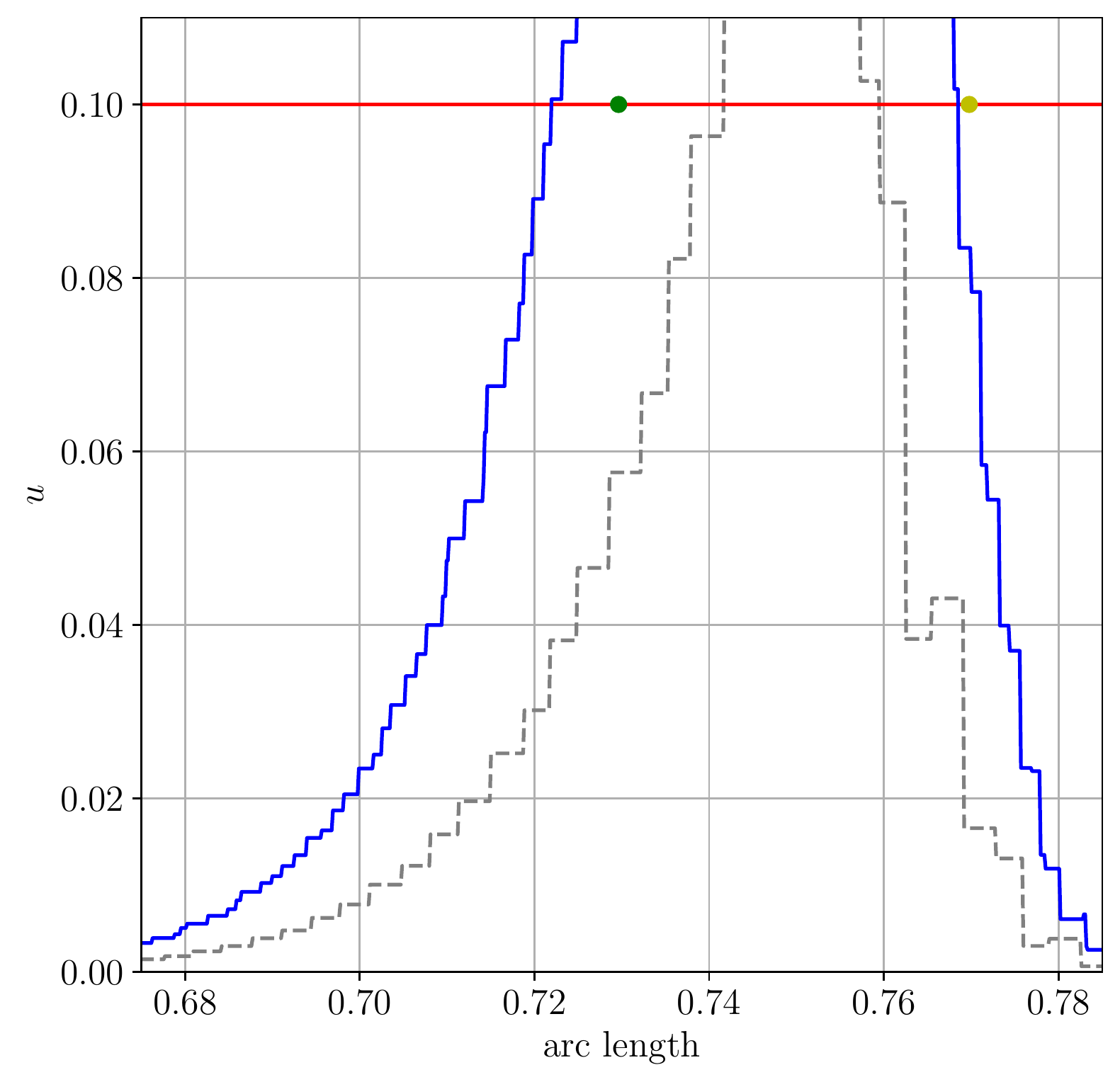}
    \caption{Profile of the solute along two lines normal to the fracture, $l_1$ on the left and $l_2$
    on the right. The blue line represent the solute profile obtained with the proposed
    multi-layer model, while the red line with the equi-dimensional layers $\mu$ and a fine grid. The green and yellow dot represent the point  $(a_{1,2}\pm\epsilon_\mu^\pm, \delta)$ for the top and bottom layer respectively. Results for the test
    case in Sub-subsection
    \ref{subsec:case1b}.}
    \label{fig:2d_case1b_profile}
\end{figure}
Finally, in Figure \ref{fig:2d_case1b_profile} we compare the solute on the same lines
specified in the previous Sub-subsection \ref{subsec:case1a}, $l_1$ and $l_2$.
The model for the thickness layer prediction is now slightly less accurate then
before, due to the effect of the not null $\eta$ parameters, however we still find
good qualitative agreement between the results.

We can conclude that, also in this setting, the multi-layer reduced model is
able to represent the effects quite accurately with a much lighter computational
cost than a refined grid, even if we are outside of the assumptions for the layer thickness evolution.

\subsubsection{Case 3}\label{subsec:case1c}

In this example we consider the same data of Case $1$ in Sub-subsection
\ref{subsec:case1a} but the reaction rate is now modeled with a non-linear
function of the solute. We set $r(u) = u^2$, $r_w(u)
=-\lambda\max\left(u^2-1,0\right)$. In this case we will note that precipitation
occurs only if $u$ exceeds $1$, the non-dimensional equilibrium value.

The aim of this test is to validate the formula \eqref{eq:thickness_layer_non_linear} for the prediction of the layer thickness.
Since this expression is derived only for the steady state and not for the
actual evolution of $\epsilon_\mu$, we cannot run the multi-layer model in
Problem \ref{pb:multi} but only the fracture-matrix Problem \ref{pb:fractured}
 and observe whether the solute/precipitate
distribution corresponds to our predictions. The use of this reaction rate in
the multi-layer model would require the derivation of an expression or an approximation of the layer thickness in time, which will be the subject of future work.

\begin{figure}[tb]
    \centering
    \includegraphics[width=0.375\textwidth]{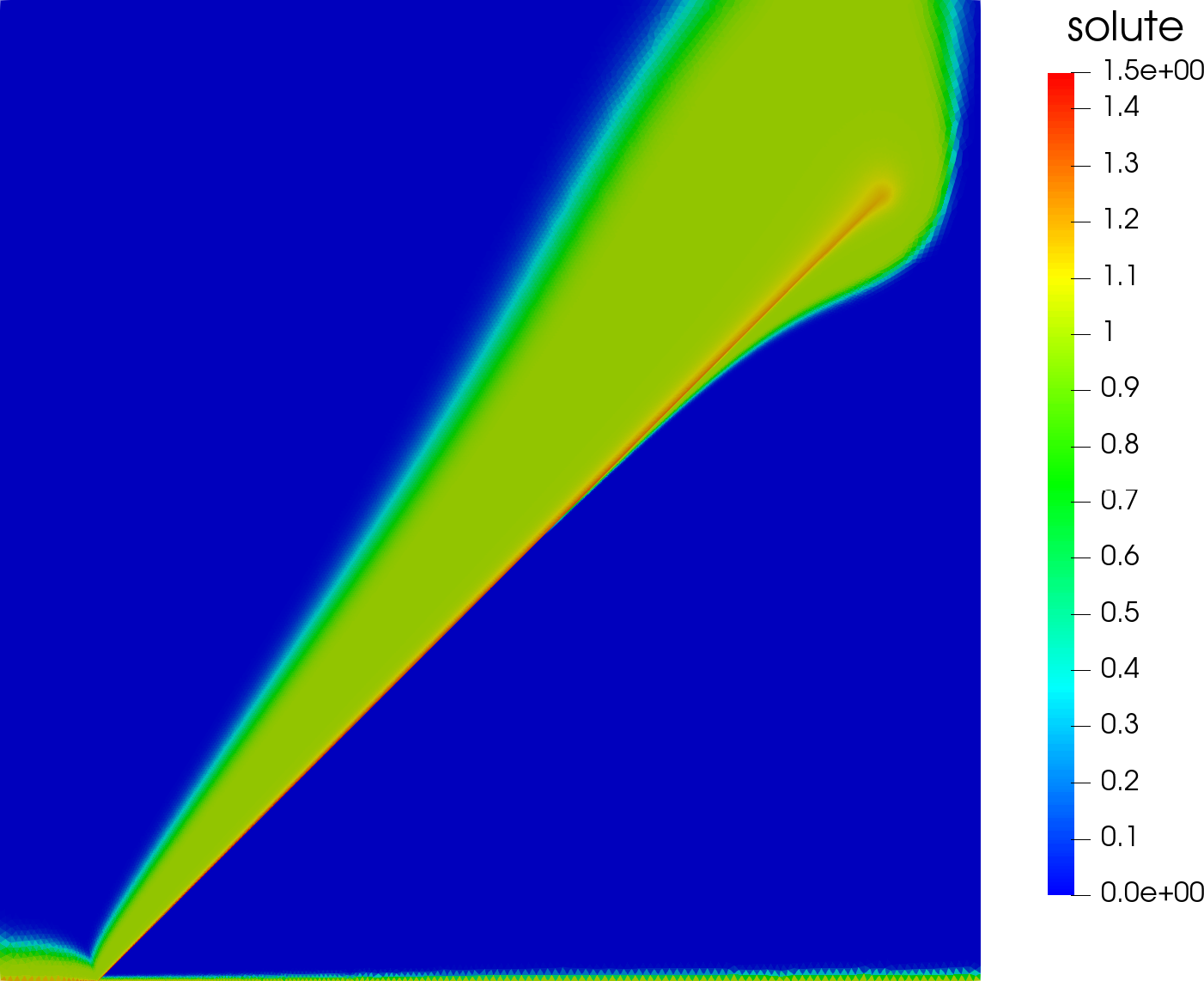}%
    \hspace*{0.05\textwidth}%
    \includegraphics[width=0.375\textwidth]{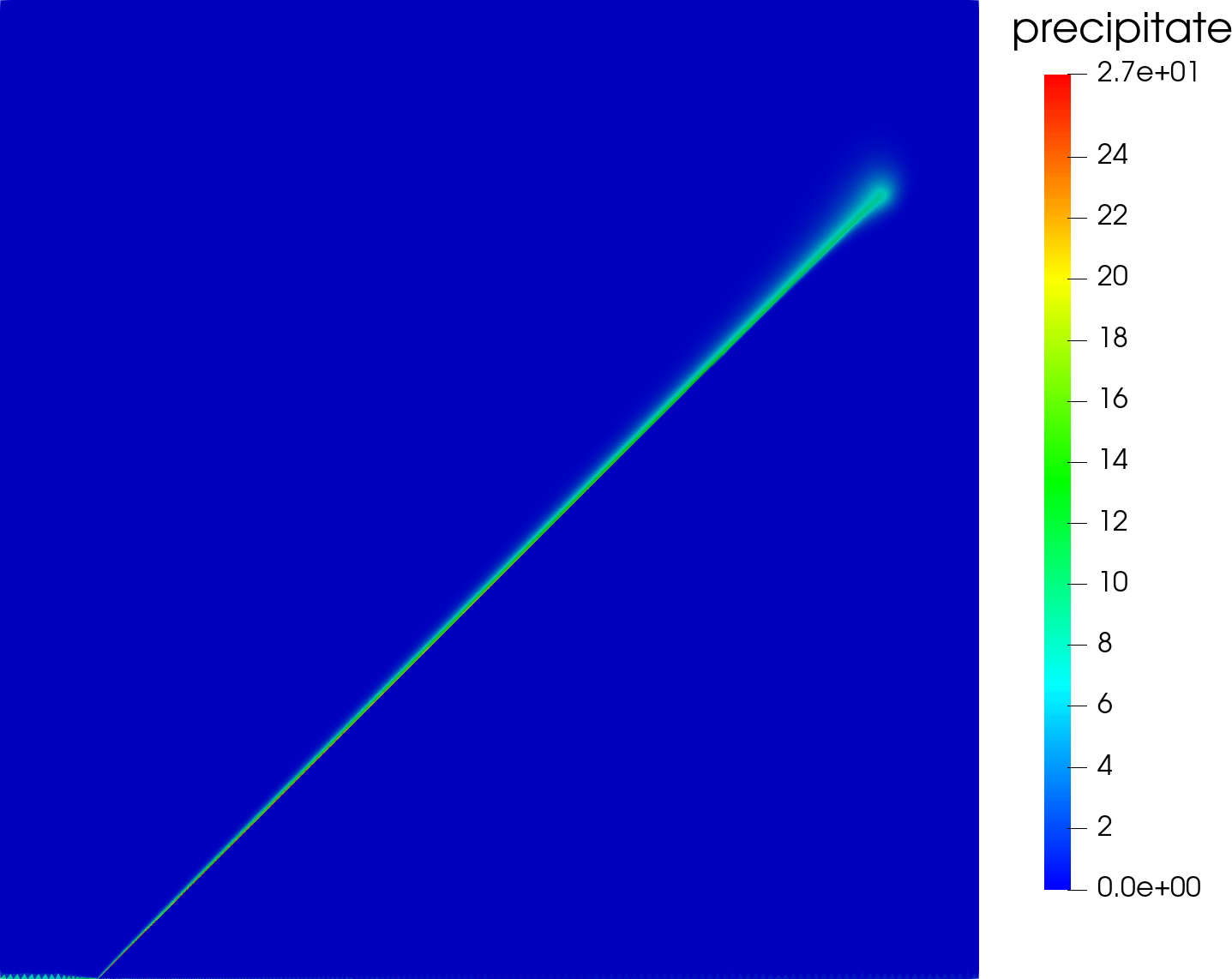}
    \caption{Solute and precipitate in the rock matrix at the end of the simulation time
    for the test case in Sub-subsection \ref{subsec:case1c}.}
    \label{fig:2d_case3_solute_precipitate}
\end{figure}

Due to the chosen data, since $\eta=0$, the porosity and fracture aperture are fixed at their initial value and thus the pressure field and Darcy velocity are the same as in Case $1$, and represented on the top of Figure \ref{fig:2d_case1}. The solute and precipitate in the rock matrix are represented in Figure \ref{fig:2d_case3_solute_precipitate}, which shows for both fields the
existence of a narrow region surrounding the fracture with very different values
than the remaining part of the rock matrix. This justifies once again the necessity
to adopt a reduced model to describe the layer around the fracture. Moreover, note that  further away form the fracture the solute concentration is below the equilibrium value ($u=1$) therefore no precipitation occurs.

\begin{figure}[tb]
    \centering
    \includegraphics[width=0.375\textwidth]{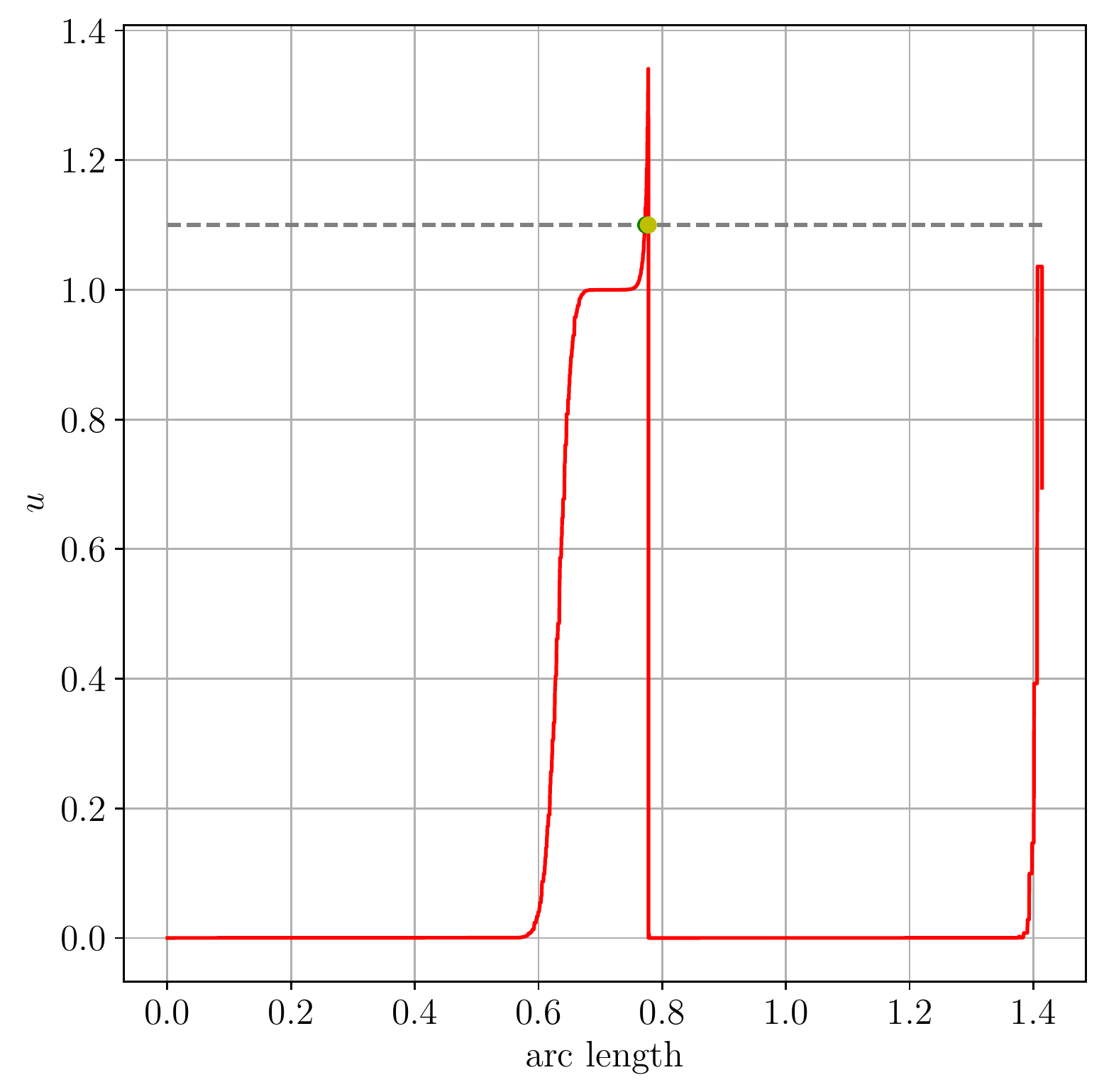}%
    \hspace*{0.05\textwidth}%
    \includegraphics[width=0.375\textwidth]{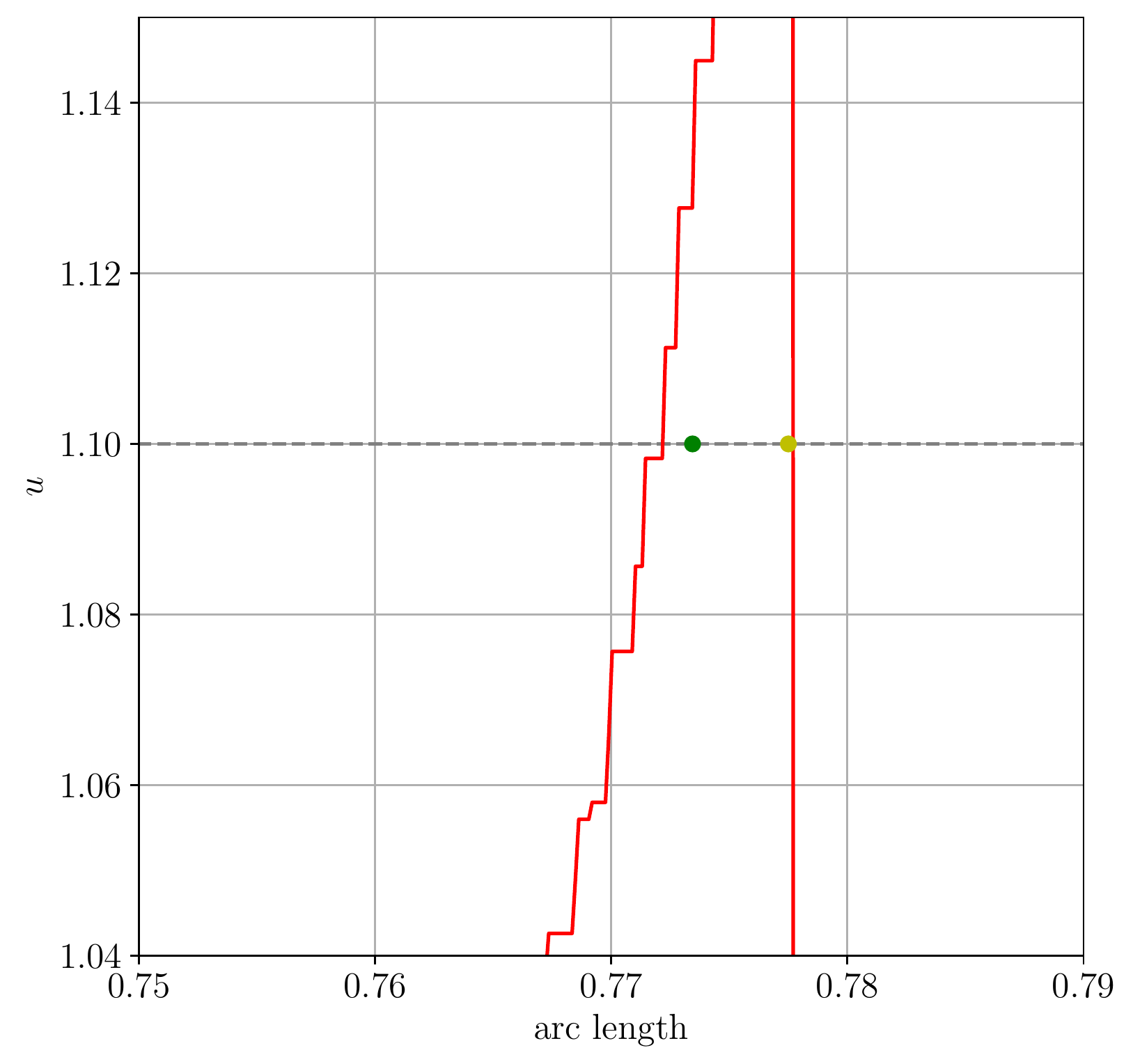}
    \includegraphics[width=0.375\textwidth]{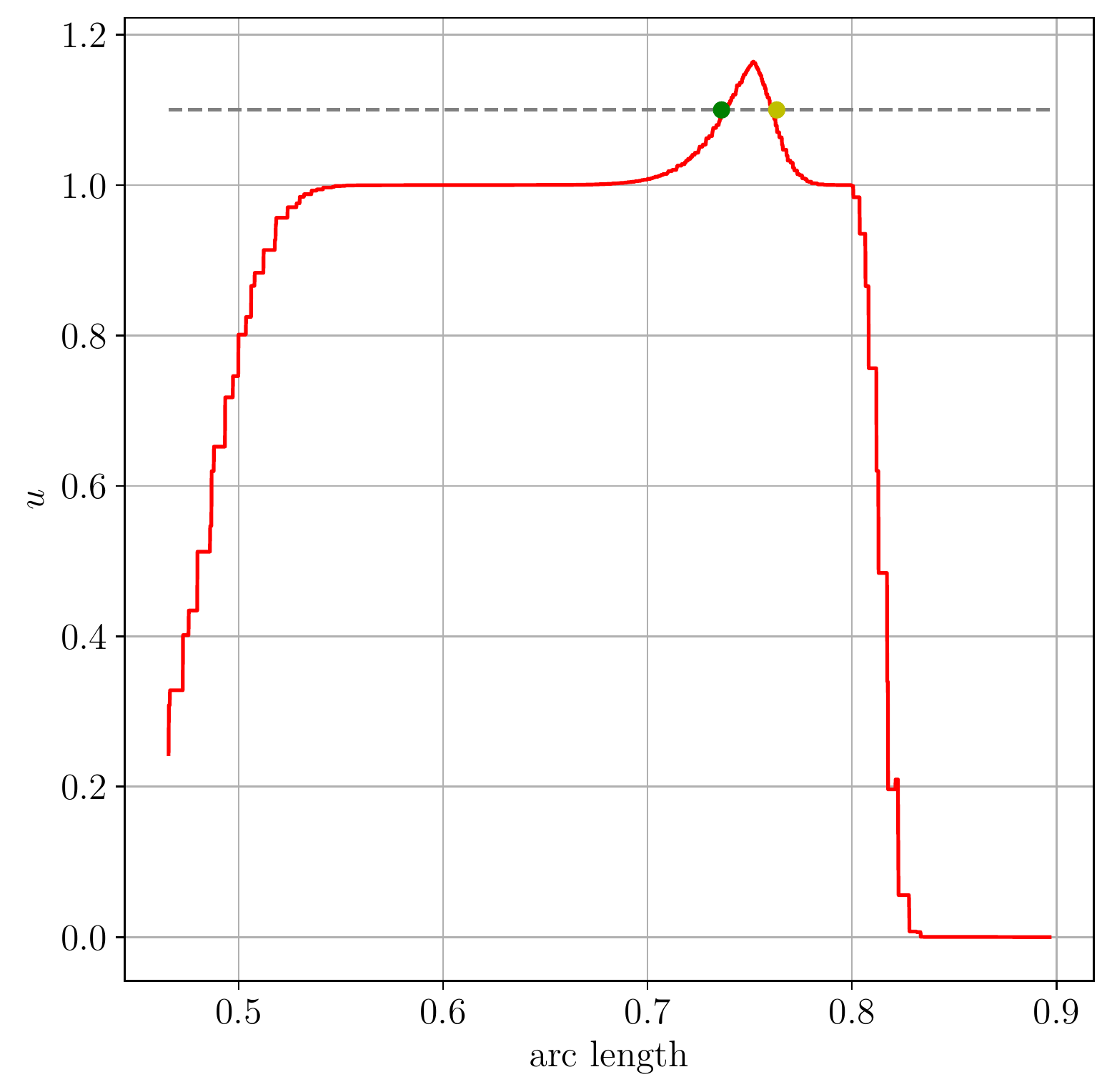}%
    \hspace*{0.05\textwidth}%
    \includegraphics[width=0.375\textwidth]{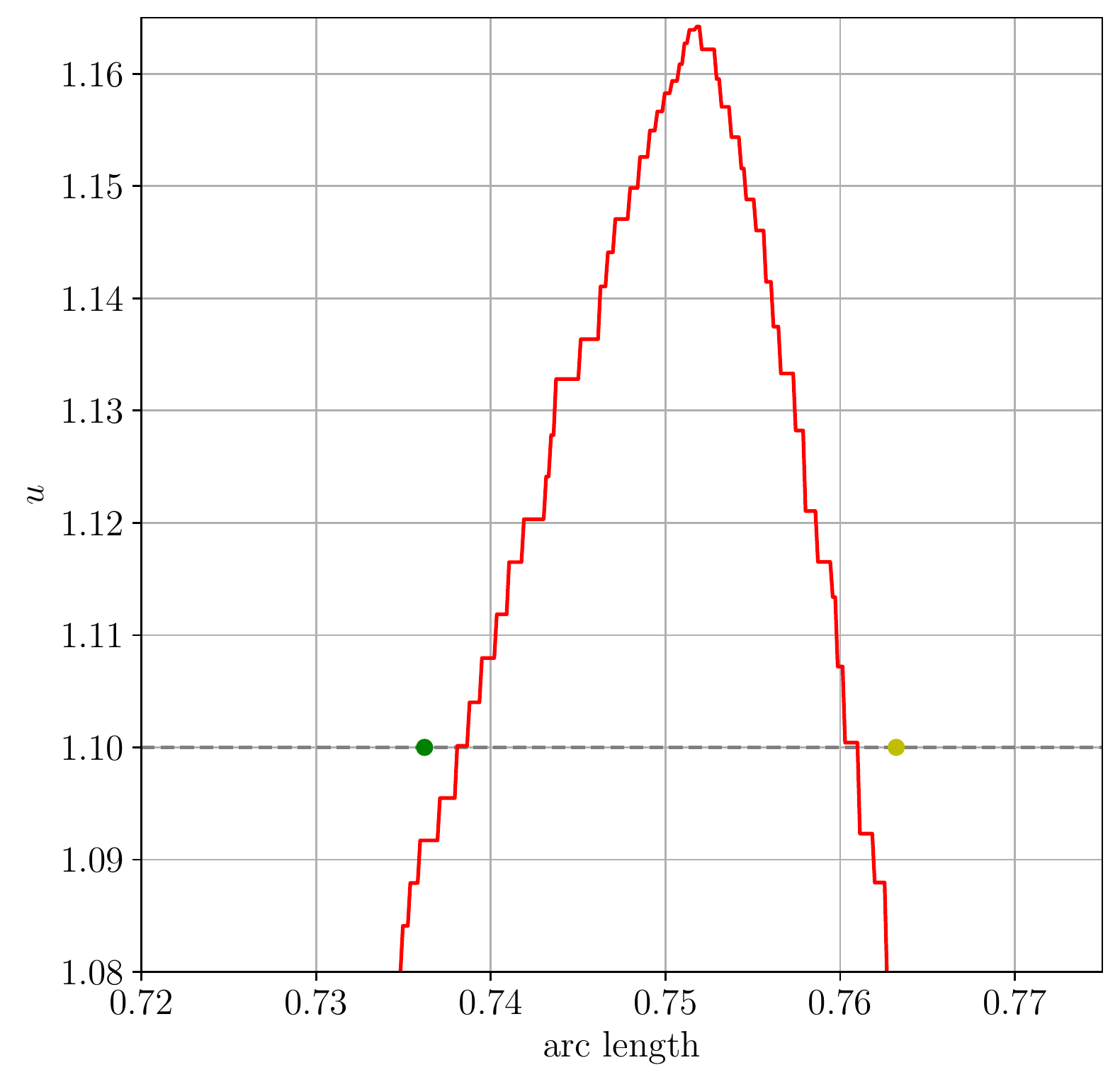}
    \caption{Profile of the solute along two lines normal to the fracture, $l_1$
    on the top and $l_2$
    on the bottom. On the right a zoom around the estimated values.
    The red line represent the solute profile obtained with the equi-dimensional layers $\mu$ and a fine grid.
    The green and yellow dot represent the point $(a_{1,2}\pm\epsilon_\mu^\pm,
    \delta+1)$ for the top and bottom layer respectively. Results for the test
    case in Sub-subsection
    \ref{subsec:case1c}.}
    \label{fig:2d_case1c_profile}
\end{figure}

Figure \ref{fig:2d_case1c_profile} shows the comparison between the layer
thickness predicted with  \eqref{eq:thickness_layer_non_linear} and the one
graphically estimated from the numerical results of model in Problem
\ref{pb:fractured}. For the comparison we consider again the
lines $l_1$ and $l_2$ introduced previously. First of all, on both lines we can
observe a peak in the solute profiles in correspondence of the fracture. The
value of solute concentration then decreases quickly reaching the plateau value
$u=1$. As done in the previous cases we can compare the predicted layer
thickness, according to \eqref{eq:thickness_layer_non_linear}, with the
numerical results: this time the dots correspond to the points in position
$(a_{1,2}\pm\epsilon_\mu^\pm, 1+\delta)$, where $a_1=0.778$ (for line $l_1$) and
$a_2=0.7521$ (line $l_2$) and $\delta=0.1$. We can observe a good agreement
between the predicted and measured layer size in both cases, even close to the
fracture tip.

Even if for the non-linear case more analysis should be done, these results can be considered promising and they confirm the feasibility of adopting a reduced model for the layer $\mu$ around the fracture.

\subsection{Three-dimensional problem}\label{subsec:case2}

\begin{figure}[tb]
    \centering
    \includegraphics[width=0.32\textwidth]{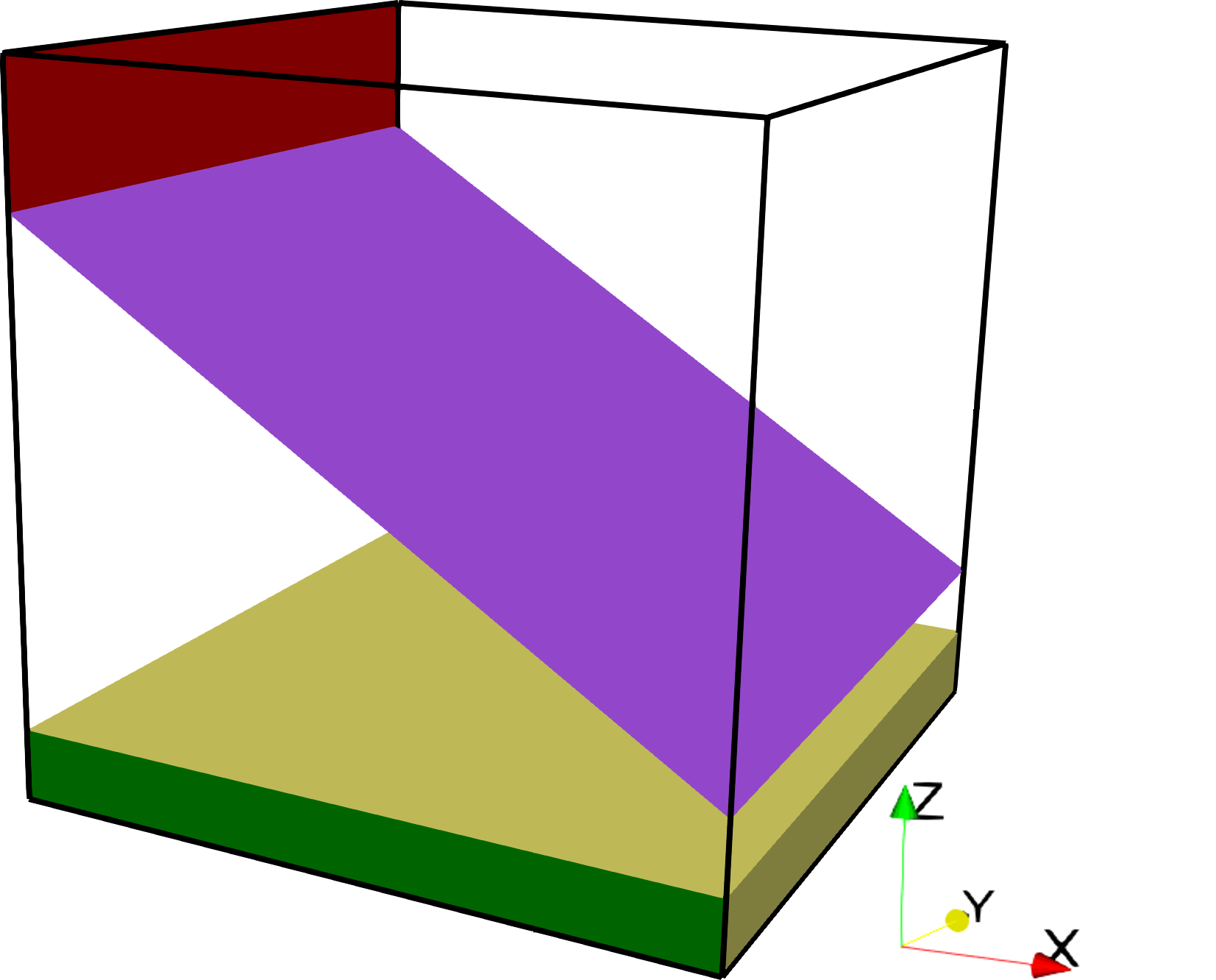}%
    \hspace*{0.01\textwidth}%
    \includegraphics[width=0.32\textwidth]{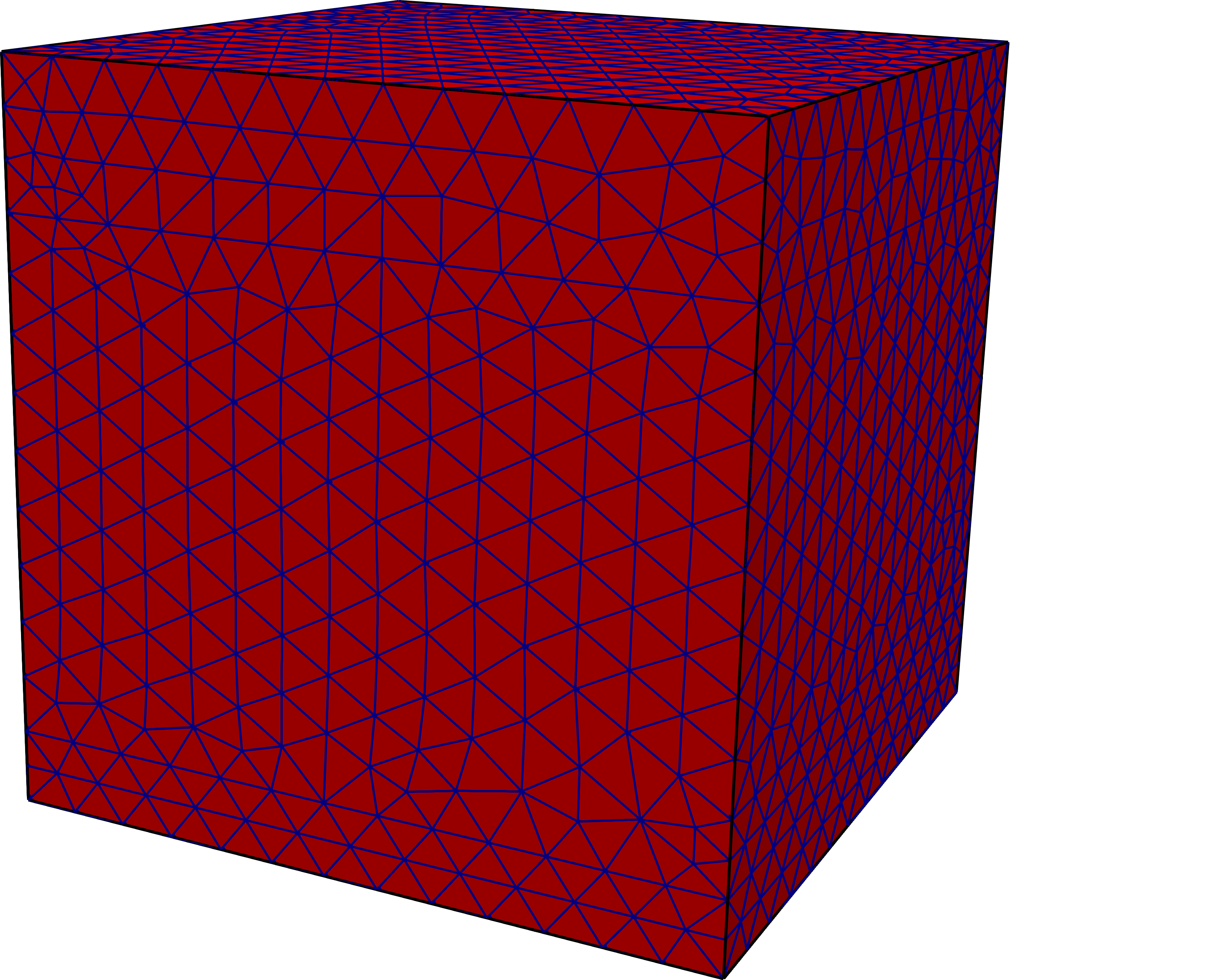}%
    \hspace*{0.01\textwidth}%
    \includegraphics[width=0.32\textwidth]{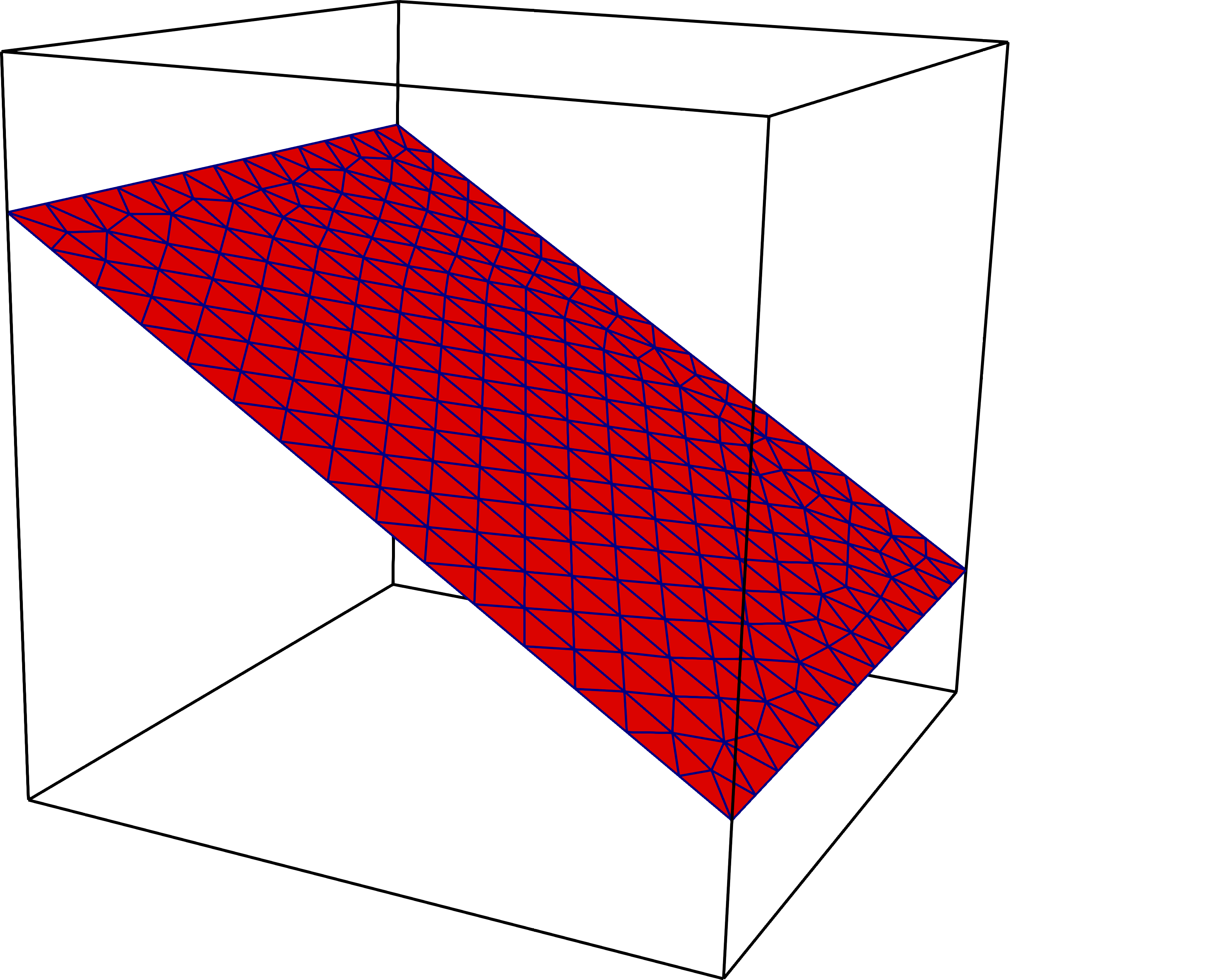}%
    \caption{Computational domain for the test case in Subsection
    \ref{subsec:case2}. The magenta represents the fracture and the layer, the
    red part is the inflow boundary region while the green is the outflow. Finally, the
    yellow block is the part with different matrix porosity and permeability.}
    \label{fig:3d_case2_domain}
\end{figure}

For this test case we consider a three dimensional setting inspired from the Case 1 of
\cite{Berre2020a}. In particular, we adopt the same geometry and part of the data for the flow problem
at the outset of the simulation.
The aim of this test case is to validate the proposed model in Problem
\ref{pb:multi} in a
three-dimensional setting.
Referring to Figure \ref{fig:3d_case2_domain}, the bottom part of the domain has
higher porosity and permeability than the remaining part. We note that the inflow
part of the boundary is slightly larger that the one in \cite{Berre2020a} to allow direct inflow into the
fracture and layer, and thus obtain a simpler flow pattern around the fracture that fits the assumptions of our model.
For the data used in the simulation see
Table \ref{tab:single_fractureb}.
\begin{table}[tb]
    \centering
    \begin{tabular}{|c|c|c|c|c|c|}
        \hline
        cutoff $=0.1$ & $\phi_{\Omega, 0}^{\rm low}=0.2$ &
        $\phi_{\Omega, 0}^{\rm high}=0.25$ & $\phi_{\mu, 0} =0.25$ &
        $\epsilon_{\gamma, 0}=4\cdot10^{-3}$
        \\ \hline
        $\epsilon_{\mu, 0} = 10^{-8}$ &
        $k_0^{\rm low}=10^{-6}$ &$k_0^{\rm high}=10^{-5}$ &  $k_{\gamma, 0} =
        10^{-1}$ &
        $\kappa_{\gamma, 0} = 10^{-1}$
        \\ \hline
        $k_{\mu, 0} = 10^{-5}$& $\kappa_{\mu,
        0}=10^{-5}$&
        $\mu=1$ & $f=0$ &
        $f_\gamma=0$
        \\ \hline
        $f_\mu = 0$ & $q_{\partial \Omega}^{\rm no-flow}=0$ &
        $p_{\partial \Omega}^{\rm out-flow} = 1$ & $p_{\partial\Omega}^{\rm
        in-flow}=4$
        & $q_{\partial \gamma}^{\rm no-flow}=0$
        \\ \hline
        $q_{\partial \mu}^{\rm no-flow}=0$ &
        $d=10^{-12}$
        &$d_\mu=10^{-12}$ &
        $\delta_\mu=10^{-12}$ &
        $d_\gamma=10^{-12}$
        \\\hline
        $\delta_\gamma=10^{-12}$&$u_{\Omega, 0}=0$
        &
        $\chi_{\partial \Omega}^{\rm no-flow}=0$ &
        $u_{\partial \Omega}^{\rm in-flow}=2$ &
        $u_{\partial \Omega}^{\rm out-flow}=0$\\\hline
        $u_{\gamma, 0}=0$ &
        $\chi_{\partial \gamma}^{\rm no-flow}=0$ &
        $u_{\partial \gamma}^{\rm in-flow}=2$ &
        $u_{\mu, 0}=0$ & $\chi_{\partial \mu}^{\rm no-flow}=0$
        \\\hline
        $u_{\partial \mu}^{\rm in-flow}=2$ &
        $\lambda^-=10^{-6}$
        &
        $r(u) = u$ & $\eta_{\Omega} = 0.5$& $\eta_{\gamma} = 0.5$
        \\\hline
        $\eta_\mu =        0.5$ &&&& \\\hline
    \end{tabular}
    \caption{Data for the examples in Subsection
    \ref{subsec:case2}.}%
    \label{tab:single_fractureb}
\end{table}
The computational grid is composed of 11436 tetrahedra for the porous media, 470
triangles for the fracture and 940 triangles for the layer.
The final simulation time is $5\cdot10^{5}$ divided uniformly in 100 time
steps. We note that the final time is shorter than in \cite{Berre2020a}
since most of the dynamic of our interest happens at an early stage.

\begin{figure}[tb]
    \centering
    \includegraphics[width=0.375\textwidth]{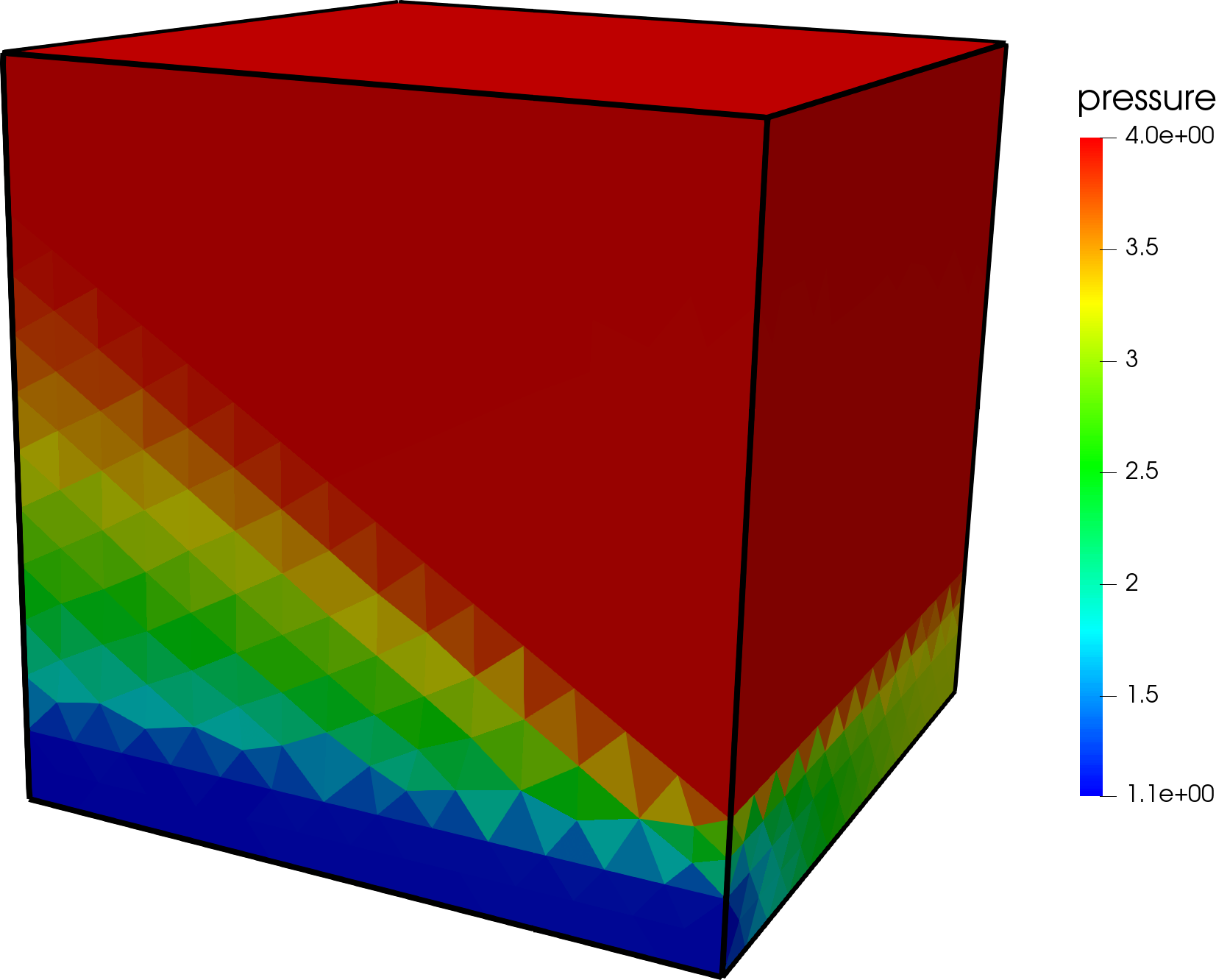}%
    \hspace*{0.05\textwidth}%
    \includegraphics[width=0.375\textwidth]{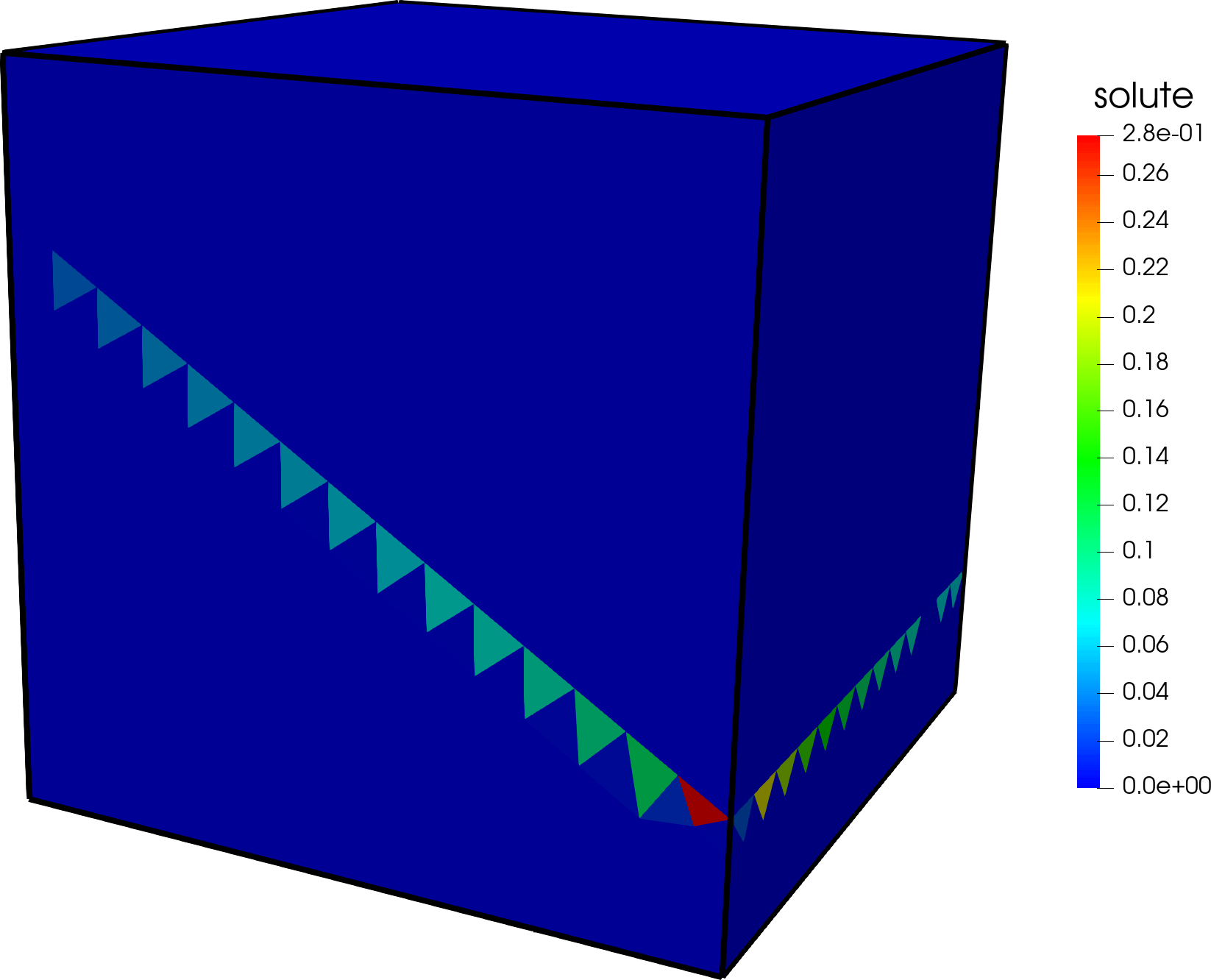}
    \caption{On the left,
    the pressure field and, on the right, the solute for the test case in Subsection
    \ref{subsec:case2}.}
    \label{fig:3d_pressure_solute_rock_matrix}
\end{figure}

\begin{figure}[tb]
    \centering
    \includegraphics[width=0.375\textwidth]{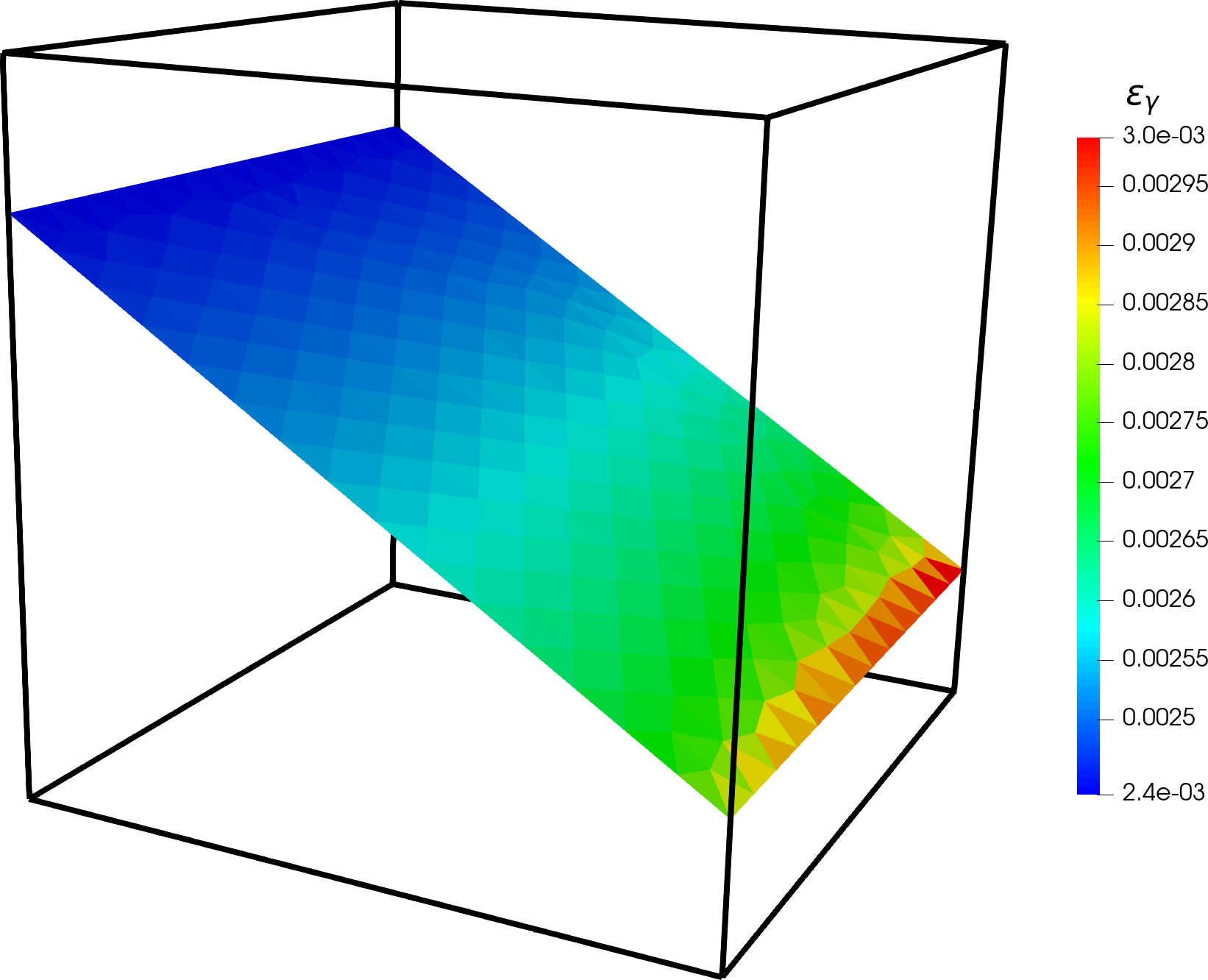}%
    \hspace*{0.05\textwidth}%
    \includegraphics[width=0.375\textwidth]{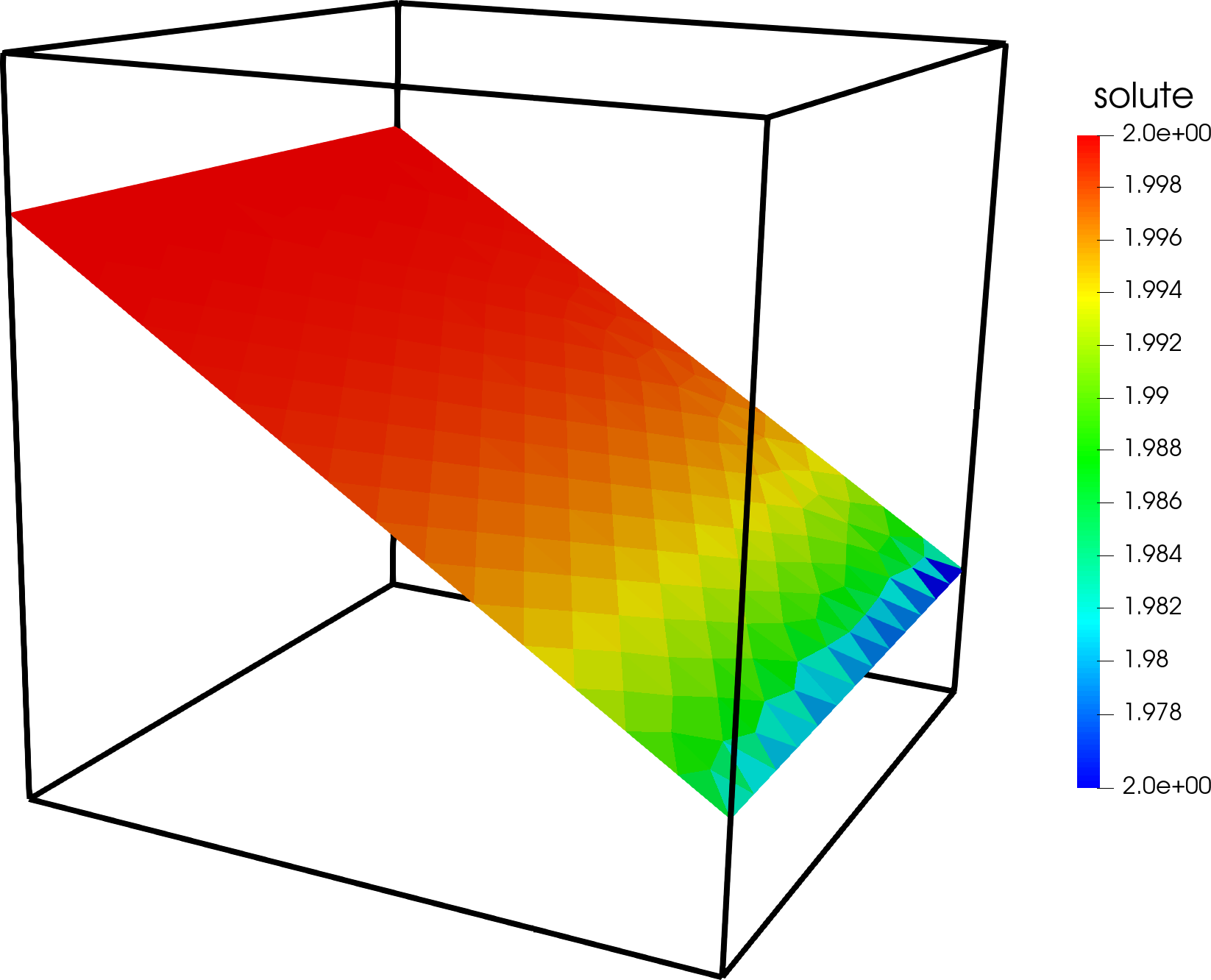}
    \caption{On the left,
    the fracture aperture and, on the right, the fracture solute for the test case in Subsection
    \ref{subsec:case2}.}
    \label{fig:3d_fracture_aperture_solute}
\end{figure}

Figure \ref{fig:3d_pressure_solute_rock_matrix} shows the pressure and solute in the
 rock matrix at the end of the simulation time. We notice that the
fracture remains highly conductive and also that the solute in the rock matrix
is quite low. Indeed, at the end of the simulation time, most of the dynamics
happened only in the fracture and surrounding layer.
In Figure \ref{fig:3d_fracture_aperture_solute} we represent the fracture
aperture and solute at the end of the simulation time. As noted before, the
fracture remains highly conductive and the inflow concentration of the solute is
transported quickly in the whole fracture. This also implies also precipitation inside the fracture and
thus fracture aperture variation, as well as a strong influence on the layer
thickness evolution.

\begin{figure}[tb]
    \centering
    \includegraphics[width=0.375\textwidth]{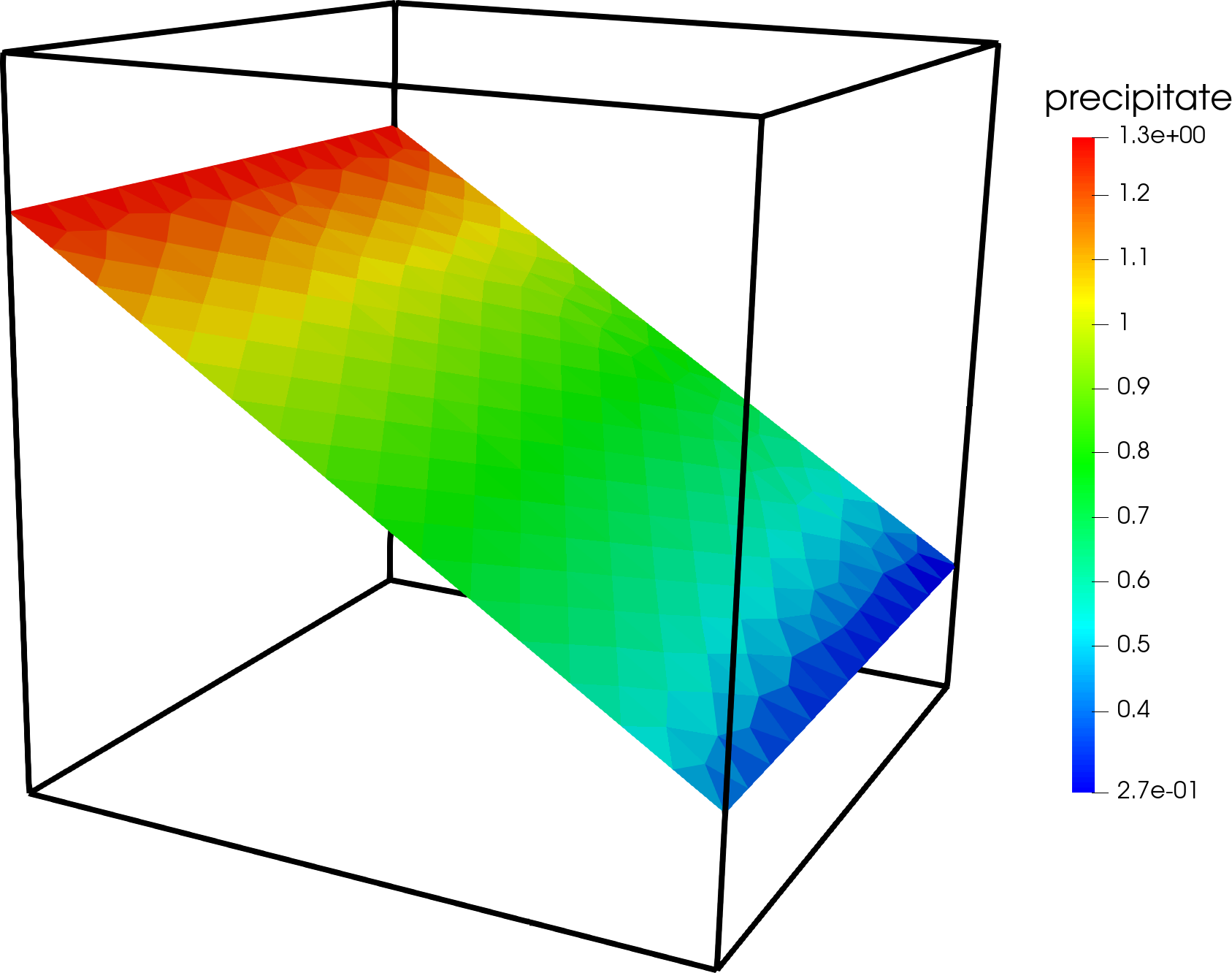}%
    \hspace*{0.05\textwidth}%
    \includegraphics[width=0.375\textwidth]{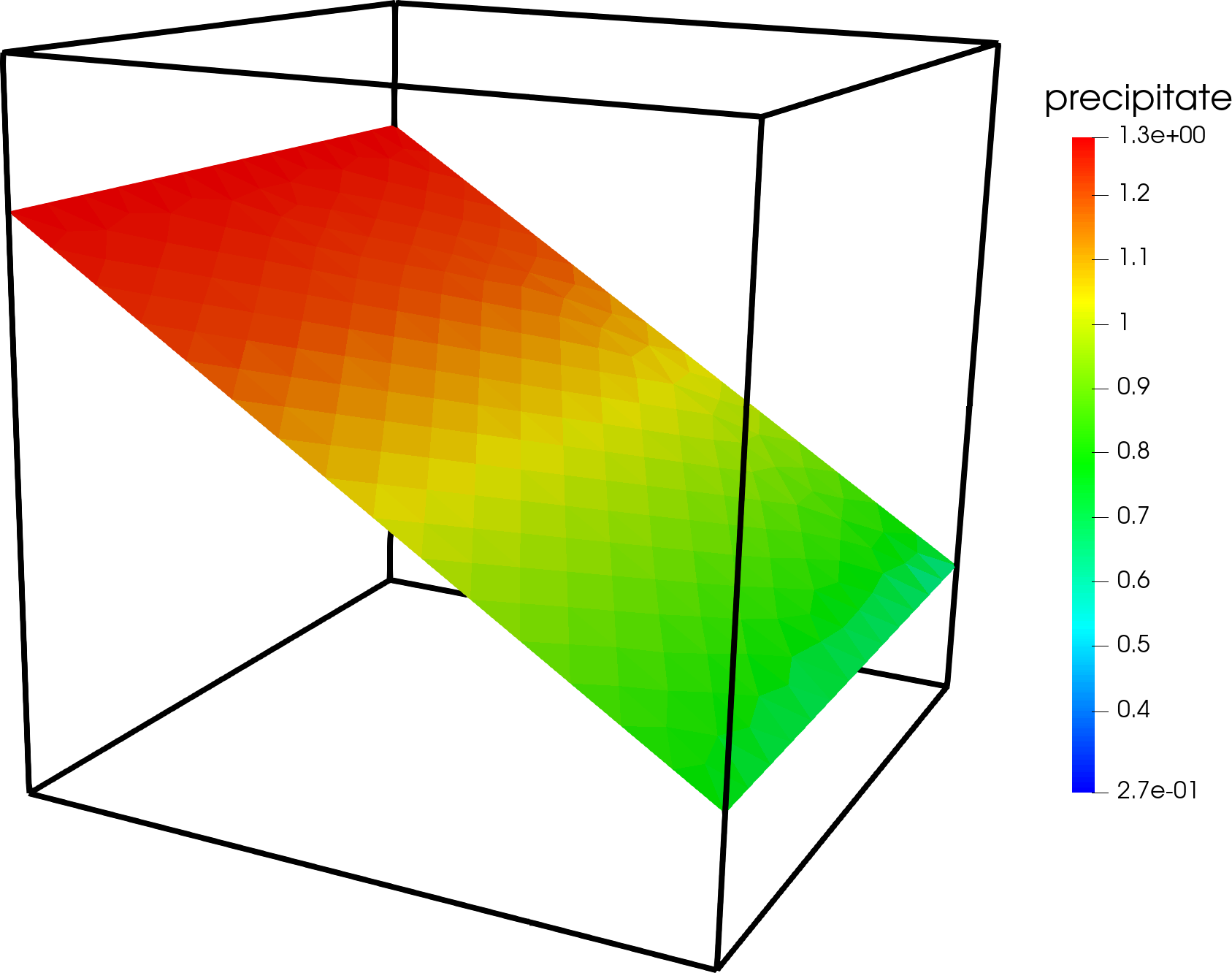}
    \includegraphics[width=0.375\textwidth]{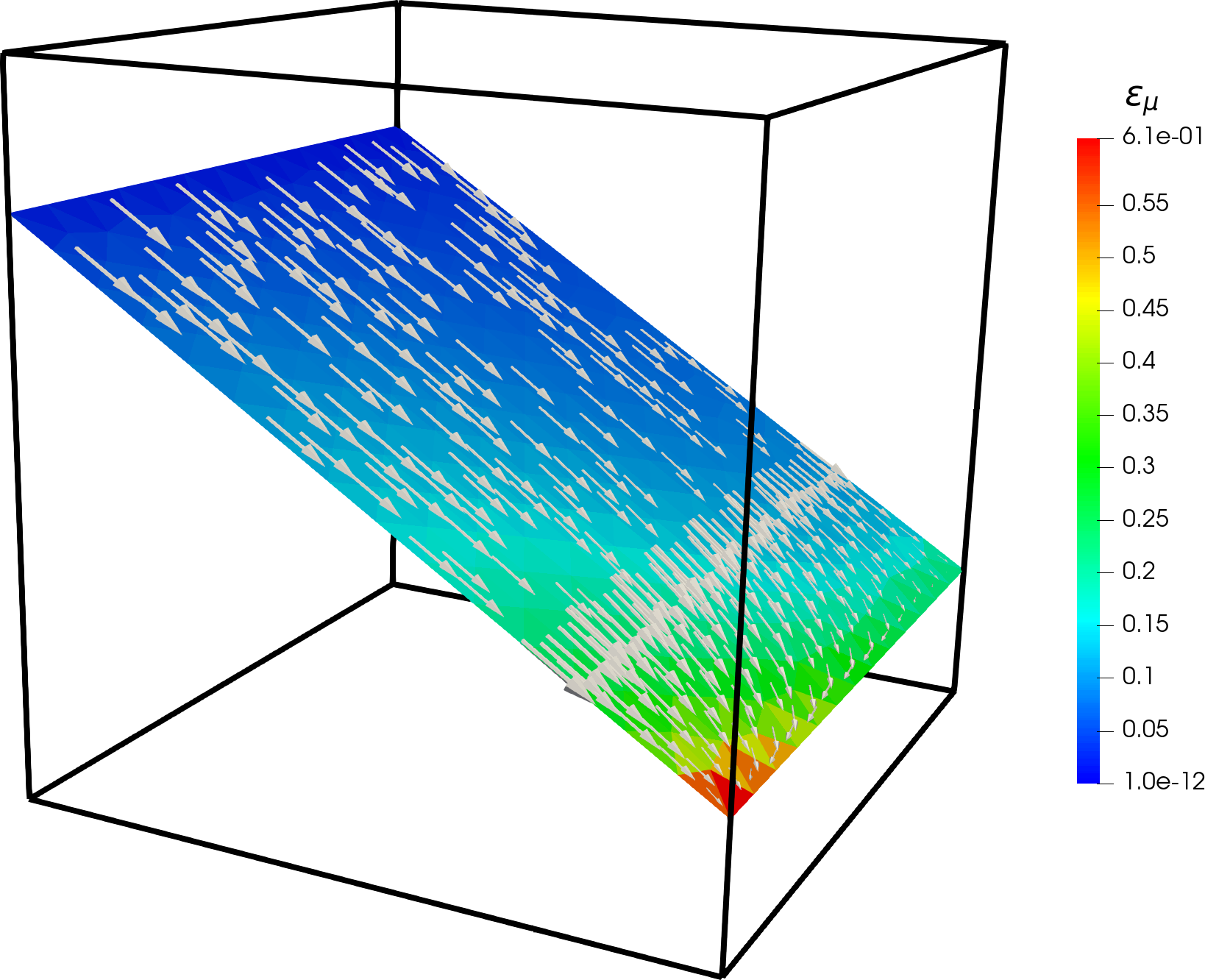}%
    \hspace*{0.05\textwidth}%
    \includegraphics[width=0.375\textwidth]{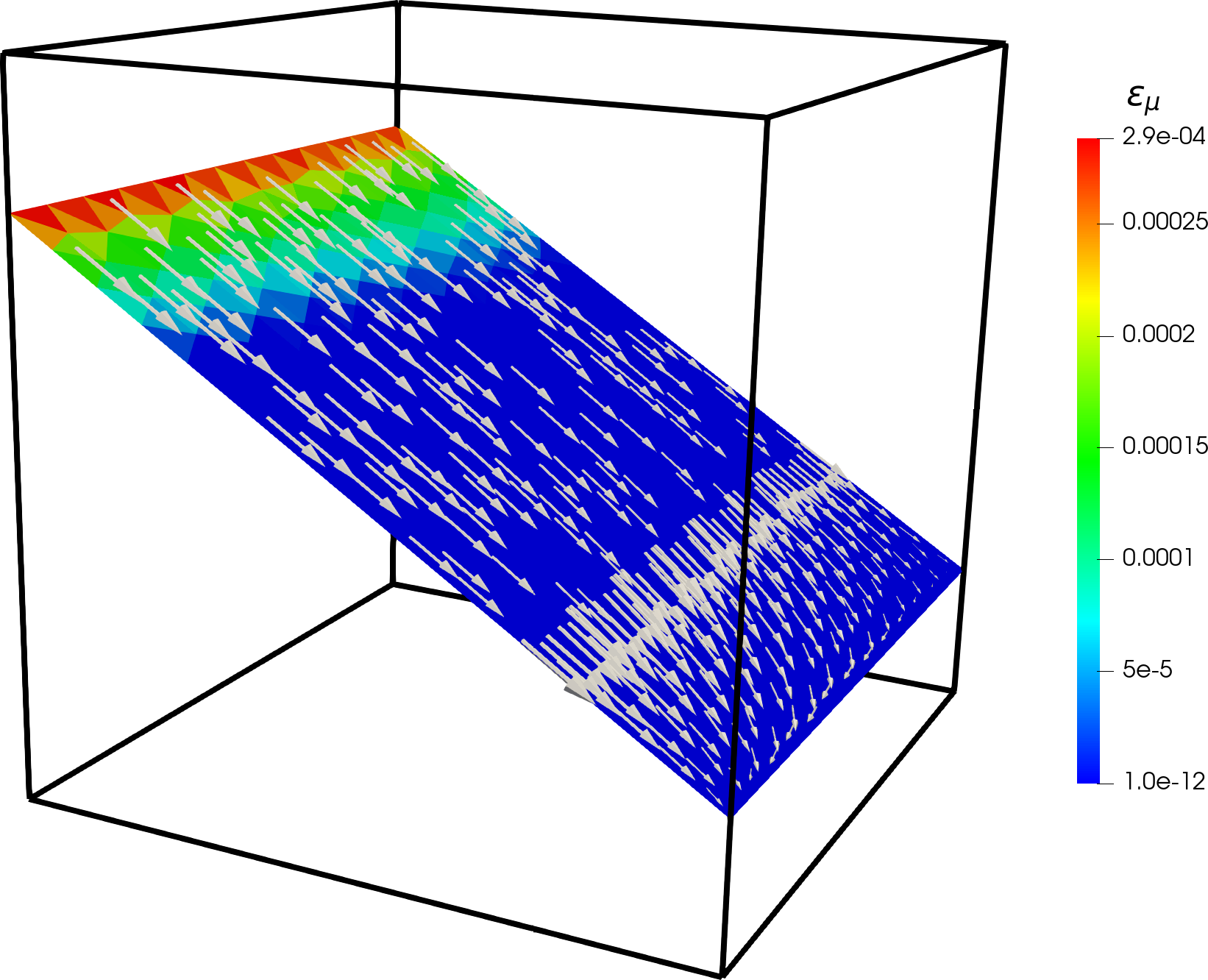}
    \caption{On the top layer precipitate and on the bottom the layer aperture.
    On the left we observe the portion of the layer on the top (close
    to the outflow) and on the right the bottom part of the layer. We overlap to the aperture
    layer the Darcy velocity in the fracture. Note that the last two figures are
    scaled differently.  Solution for the test case in Subsection
    \ref{subsec:case2}.}
    \label{fig:3d_layer_precipitate_aperture}
\end{figure}

Figure \ref{fig:3d_layer_precipitate_aperture} show the dynamics inside the
layer. We obtain more precipitate in the top part of the layer $\mu$ due to
the inflow into the layer itself from the top part of the rock
matrix, and also because, in the bottom part of $\mu$, the solute tends to flow towards the outflow boundary at the bottom, resulting in a smaller concentration of precipitate in the bottom part of $\mu$.
The layer thickness is also represented, with two different scales, overlapped
with the Darcy velocity in the fracture, as a proxy for the flow
exchange between the fracture and the layer. We see that the top part of the
layer is rather thin and in principle might be neglected, however
on the bottom part a higher value of the thickness reveals the importance of
having the layer explicitly represented. The aperture in this case is not
uniform, but rather, larger near the outflow of the problem, as one could
expect.

Considering the size of the computational domain, the values of the layer
thickness obtained is in the limit of a reduced model. To be able to capture this small
layer around the fracture, and thus use the model in Problem
\ref{pb:fractured}, we should refine the grid obtaining a problem that is too
computational expensive to solve, even for such simple test case. This test case, with the considered data, shows the importance of the presented multi-layer reduced model, which can be considered an attractive alternative.

\section{Conclusion}\label{sec:conclusion}

In this work we have introduced a mathematical model that is able to simulate in
an accurate, yet affordable way simple reactive transport flow problems in the
presence of a fracture. In particular, when the reaction rate is high enough
compared to transport, we observe that a narrow region, denoted as reactive
layer, forms just around the fracture: here the porous medium has different
physical properties from the surrounding porous matrix due to mineral
precipitation or dissolution. These changes in porosity and permeability might
substantially alter the flow field, resulting in a fully coupled and non-linear
mathematical model. Moreover, in this layer we expect steep gradients of the
variables, in particular the solute and precipitate. For large Damk\"ohler
numbers, as proven by numerical simulations and experimental observations, these
reactive layers can be extremely thin, to the point that it is difficult to
capture their geometry and solution dynamics with a refined computational grid.
For this reason in this work we have proposed and tested a reduced model where
not only the fractures, but the reactive layers as well are represented as
co-dimension 1 objects coupled with the porous matrix, and among themselves. We
have derived, under suitable assumptions, a model for the evolution in time of
the layer thickness which provided reliable results compared to a very refined
numerical simulation of the corresponding equi-dimensional model.  The model has
been derived and tested for a simple linear reaction rate model and, at the
steady state only, for a more complex reaction rate that accounts for
equilibrium solubility and supersaturation. By increasing further the complexity
of the reaction rate we expect that the model for the layer evolution might
become more involved and will require a numerical approximated solution (as
opposed to a closed form expression) to estimate, at each point and each time
step, the layer thickness: this will be part of a future study. In the numerical
study we have also shown a three-dimensional model where the proposed approach
might be even more attractive to substantially lighten the computational burden
associated with mesh refinement.

%\section*{Acknowledgments}

\section*{Conflict of interest}

There are no conflict of interest.

\bibliographystyle{plain}
\bibliography{biblio}

%%\begin{thebibliography}{999}
%%
%%\bibitem{authour1}
%%     \newblock  \textbf{Journal article style:} Benoist Y, Foulon P, Labourie F, et al. (Year) Anosov flows with stable and unstable differentiable distributions.
%%     \newblock \emph{ J Amer Math Soc}
%%     \newblock  Volume: StaringPage{--\hspace*{-2mm}--}Ending Page.
%%
%%
%%\bibitem{authour2}
%%    \newblock \textbf{Book style:}
%%    \newblock Serrin J, (1971) Gradient estimates for solutions of nonlinear elliptic and parabolic equations, In: Zarantonello, E.Z. Author,
%%    \newblock \emph{Contributions to Nonlinear Functional Analysis},
%%    \newblock  2 Eds., New York: Academic Press, 35{--\hspace*{-2mm}--}75.
%%
%%
%%\bibitem{authour3}
%%    \newblock \textbf{Online content:}
%%    \newblock  SARS Expert Committee, SARS in Hong Kong: From Experience to Action. Hong Kong SARS Expert Committee, 2003. Available from: \\ \url{http://www.sars-expertcom.gov.hk/english/reports/reports.html.}
%%
%%\end{thebibliography}

\end{document}

%% file: fig/domain_equi.pdf_tex
%% Creator: Inkscape inkscape 0.92.5, www.inkscape.org
%% Accompanies image file 'domain_equi.pdf' (pdf, eps, ps)
\begingroup%
  \makeatletter%
  \providecommand\color[2][]{%
    \errmessage{(Inkscape) Color is used for the text in Inkscape, but the package 'color.sty' is not loaded}%
    \renewcommand\color[2][]{}%
  }%
  \providecommand\transparent[1]{%
    \errmessage{(Inkscape) Transparency is used (non-zero) for the text in Inkscape, but the package 'transparent.sty' is not loaded}%
    \renewcommand\transparent[1]{}%
  }%
  \providecommand\rotatebox[2]{#2}%
  \newcommand*\fsize{\dimexpr\f@size pt\relax}%
  \newcommand*\lineheight[1]{\fontsize{\fsize}{#1\fsize}\selectfont}%
  \ifx\svgwidth\undefined%
    \setlength{\unitlength}{168.909983bp}%
    \ifx\svgscale\undefined%
      \relax%
    \else%
      \setlength{\unitlength}{\unitlength * \real{\svgscale}}%
    \fi%
  \else%
    \setlength{\unitlength}{\svgwidth}%
  \fi%
  \global\let\svgwidth\undefined%
  \global\let\svgscale\undefined%
  \makeatother%
  \begin{picture}(1,0.62459253)%
    \lineheight{1}%
    \setlength\tabcolsep{0pt}%
    \put(0,0){\includegraphics[width=\unitlength,page=1]{domain_equi.pdf}}%
    \put(0.11360182,0.43425379){\color[rgb]{0,0,0}\makebox(0,0)[lt]{\lineheight{0}\smash{\begin{tabular}[t]{l}$\Omega^-$\end{tabular}}}}%
    \put(0,0){\includegraphics[width=\unitlength,page=2]{domain_equi.pdf}}%
    \put(0.78066455,0.43425375){\color[rgb]{0,0,0}\makebox(0,0)[lt]{\lineheight{0}\smash{\begin{tabular}[t]{l}$\Omega^+$\end{tabular}}}}%
    \put(0,0){\includegraphics[width=\unitlength,page=3]{domain_equi.pdf}}%
    \put(0.38128896,0.3484496){\color[rgb]{0,0,0}\makebox(0,0)[lt]{\lineheight{0}\smash{\begin{tabular}[t]{l}$\mu^-$\end{tabular}}}}%
    \put(0.5742286,0.34540041){\color[rgb]{0,0,0}\makebox(0,0)[lt]{\lineheight{0}\smash{\begin{tabular}[t]{l}$\mu^+$\end{tabular}}}}%
    \put(0.48172507,0.22307793){\color[rgb]{0,0,0}\makebox(0,0)[lt]{\lineheight{0}\smash{\begin{tabular}[t]{l}$\gamma$\end{tabular}}}}%
  \end{picture}%
\endgroup%

%% file: fig/domain_mono.pdf_tex
%% Creator: Inkscape inkscape 0.92.5, www.inkscape.org
%% Accompanies image file 'domain_mono.pdf' (pdf, eps, ps)
\begingroup%
  \makeatletter%
  \providecommand\color[2][]{%
    \errmessage{(Inkscape) Color is used for the text in Inkscape, but the package 'color.sty' is not loaded}%
    \renewcommand\color[2][]{}%
  }%
  \providecommand\transparent[1]{%
    \errmessage{(Inkscape) Transparency is used (non-zero) for the text in Inkscape, but the package 'transparent.sty' is not loaded}%
    \renewcommand\transparent[1]{}%
  }%
  \providecommand\rotatebox[2]{#2}%
  \newcommand*\fsize{\dimexpr\f@size pt\relax}%
  \newcommand*\lineheight[1]{\fontsize{\fsize}{#1\fsize}\selectfont}%
  \ifx\svgwidth\undefined%
    \setlength{\unitlength}{231.8740244bp}%
    \ifx\svgscale\undefined%
      \relax%
    \else%
      \setlength{\unitlength}{\unitlength * \real{\svgscale}}%
    \fi%
  \else%
    \setlength{\unitlength}{\svgwidth}%
  \fi%
  \global\let\svgwidth\undefined%
  \global\let\svgscale\undefined%
  \makeatother%
  \begin{picture}(1,0.54716416)%
    \lineheight{1}%
    \setlength\tabcolsep{0pt}%
    \put(0,0){\includegraphics[width=\unitlength,page=1]{domain_mono.pdf}}%
    \put(0.47950963,0.5157085){\color[rgb]{0,0,0}\makebox(0,0)[lt]{\lineheight{0}\smash{\begin{tabular}[t]{l}$\gamma$\end{tabular}}}}%
    \put(0.42870626,0.00608494){\color[rgb]{0,0,0}\makebox(0,0)[lt]{\lineheight{0}\smash{\begin{tabular}[t]{l}$\Gamma^-$\end{tabular}}}}%
    \put(0.5293886,0.00689657){\color[rgb]{0,0,0}\makebox(0,0)[lt]{\lineheight{0}\smash{\begin{tabular}[t]{l}$\Gamma^+$\end{tabular}}}}%
    \put(0,0){\includegraphics[width=\unitlength,page=2]{domain_mono.pdf}}%
    \put(0.10109607,0.36145703){\color[rgb]{0,0,0}\makebox(0,0)[lt]{\lineheight{0}\smash{\begin{tabular}[t]{l}$\Omega^-$\end{tabular}}}}%
    \put(0,0){\includegraphics[width=\unitlength,page=3]{domain_mono.pdf}}%
    \put(0.84562085,0.38098076){\color[rgb]{0,0,0}\makebox(0,0)[lt]{\lineheight{0}\smash{\begin{tabular}[t]{l}$\Omega^+$\end{tabular}}}}%
    \put(0,0){\includegraphics[width=\unitlength,page=4]{domain_mono.pdf}}%
    \put(0.27038654,0.20357747){\color[rgb]{0,0,0}\makebox(0,0)[lt]{\lineheight{0}\smash{\begin{tabular}[t]{l}$\mu^-$\end{tabular}}}}%
    \put(0.68835571,0.20808294){\color[rgb]{0,0,0}\makebox(0,0)[lt]{\lineheight{0}\smash{\begin{tabular}[t]{l}$\mu^+$\end{tabular}}}}%
    \put(0,0){\includegraphics[width=\unitlength,page=5]{domain_mono.pdf}}%
    \put(0.51766801,0.29043068){\color[rgb]{0,0,0}\makebox(0,0)[lt]{\lineheight{1.25}\smash{\begin{tabular}[t]{l}\textbf{n}\end{tabular}}}}%
  \end{picture}%
\endgroup%

%% file: fig/domain.pdf_tex
%% Creator: Inkscape inkscape 0.92.5, www.inkscape.org
%% Accompanies image file 'domain.pdf' (pdf, eps, ps)
\begingroup%
  \makeatletter%
  \providecommand\color[2][]{%
    \errmessage{(Inkscape) Color is used for the text in Inkscape, but the package 'color.sty' is not loaded}%
    \renewcommand\color[2][]{}%
  }%
  \providecommand\transparent[1]{%
    \errmessage{(Inkscape) Transparency is used (non-zero) for the text in Inkscape, but the package 'transparent.sty' is not loaded}%
    \renewcommand\transparent[1]{}%
  }%
  \providecommand\rotatebox[2]{#2}%
  \newcommand*\fsize{\dimexpr\f@size pt\relax}%
  \newcommand*\lineheight[1]{\fontsize{\fsize}{#1\fsize}\selectfont}%
  \ifx\svgwidth\undefined%
    \setlength{\unitlength}{231.8740244bp}%
    \ifx\svgscale\undefined%
      \relax%
    \else%
      \setlength{\unitlength}{\unitlength * \real{\svgscale}}%
    \fi%
  \else%
    \setlength{\unitlength}{\svgwidth}%
  \fi%
  \global\let\svgwidth\undefined%
  \global\let\svgscale\undefined%
  \makeatother%
  \begin{picture}(1,0.54716416)%
    \lineheight{1}%
    \setlength\tabcolsep{0pt}%
    \put(0,0){\includegraphics[width=\unitlength,page=1]{domain.pdf}}%
    \put(0.37004535,0.5157085){\color[rgb]{0,0,0}\makebox(0,0)[lt]{\lineheight{0}\smash{\begin{tabular}[t]{l}$\mu^-$\end{tabular}}}}%
    \put(0.59681472,0.5157085){\color[rgb]{0,0,0}\makebox(0,0)[lt]{\lineheight{0}\smash{\begin{tabular}[t]{l}$\mu^+$\end{tabular}}}}%
    \put(0.47950963,0.5157085){\color[rgb]{0,0,0}\makebox(0,0)[lt]{\lineheight{0}\smash{\begin{tabular}[t]{l}$\gamma$\end{tabular}}}}%
    \put(0.29940362,0.00608494){\color[rgb]{0,0,0}\makebox(0,0)[lt]{\lineheight{0}\smash{\begin{tabular}[t]{l}$M^-$\end{tabular}}}}%
    \put(0.64347807,0.00608494){\color[rgb]{0,0,0}\makebox(0,0)[lt]{\lineheight{0}\smash{\begin{tabular}[t]{l}$M^+$\end{tabular}}}}%
    \put(0.42870626,0.00608494){\color[rgb]{0,0,0}\makebox(0,0)[lt]{\lineheight{0}\smash{\begin{tabular}[t]{l}$\Gamma^-$\end{tabular}}}}%
    \put(0.54290502,0.00608494){\color[rgb]{0,0,0}\makebox(0,0)[lt]{\lineheight{0}\smash{\begin{tabular}[t]{l}$\Gamma^+$\end{tabular}}}}%
    \put(0,0){\includegraphics[width=\unitlength,page=2]{domain.pdf}}%
    \put(0.15342667,0.27134751){\color[rgb]{0,0,0}\makebox(0,0)[lt]{\lineheight{0}\smash{\begin{tabular}[t]{l}$\Omega^-$\end{tabular}}}}%
    \put(0,0){\includegraphics[width=\unitlength,page=3]{domain.pdf}}%
    \put(0.80830852,0.27134751){\color[rgb]{0,0,0}\makebox(0,0)[lt]{\lineheight{0}\smash{\begin{tabular}[t]{l}$\Omega^+$\end{tabular}}}}%
    \put(0,0){\includegraphics[width=\unitlength,page=4]{domain.pdf}}%
    \put(0.51451884,0.28555574){\color[rgb]{0,0,0}\makebox(0,0)[lt]{\lineheight{1.25}\smash{\begin{tabular}[t]{l}$\bm{n}$\end{tabular}}}}%
  \end{picture}%
\endgroup%

%% file: fig/examples/ex1/domain_single_frac.pdf_tex
%% Accompanies image file 'domain_single_frac.pdf' (pdf, eps, ps)
\begingroup%
  \makeatletter%
  \providecommand\color[2][]{%
    \errmessage{(Inkscape) Color is used for the text in Inkscape, but the package 'color.sty' is not loaded}%
    \renewcommand\color[2][]{}%
  }%
  \providecommand\transparent[1]{%
    \errmessage{(Inkscape) Transparency is used (non-zero) for the text in Inkscape, but the package 'transparent.sty' is not loaded}%
    \renewcommand\transparent[1]{}%
  }%
  \providecommand\rotatebox[2]{#2}%
  \newcommand*\fsize{\dimexpr\f@size pt\relax}%
  \newcommand*\lineheight[1]{\fontsize{\fsize}{#1\fsize}\selectfont}%
  \ifx\svgwidth\undefined%
    \setlength{\unitlength}{311.96932983bp}%
    \ifx\svgscale\undefined%
      \relax%
    \else%
      \setlength{\unitlength}{\unitlength * \real{\svgscale}}%
    \fi%
  \else%
    \setlength{\unitlength}{\svgwidth}%
  \fi%
  \global\let\svgwidth\undefined%
  \global\let\svgscale\undefined%
  \makeatother%
  \begin{picture}(1,0.98417801)%
    \lineheight{1}%
    \setlength\tabcolsep{0pt}%
    \put(0,0){\includegraphics[width=\unitlength,page=1]{domain_single_frac.pdf}}%
    \put(0.40318924,0.93547031){\color[rgb]{0,0,0}\makebox(0,0)[lt]{\lineheight{1.25}\smash{\begin{tabular}[t]{l}out-flow\end{tabular}}}}%
    \put(0.40068498,0.0009078){\color[rgb]{0,0,0}\makebox(0,0)[lt]{\lineheight{1.25}\smash{\begin{tabular}[t]{l}in-flow\end{tabular}}}}%
    \put(0.9990922,0.40283447){\color[rgb]{0,0,0}\rotatebox{90}{\makebox(0,0)[lt]{\lineheight{1.25}\smash{\begin{tabular}[t]{l}no-flow\end{tabular}}}}}%
    \put(0.00090779,0.64233492){\color[rgb]{0,0,0}\rotatebox{-90}{\makebox(0,0)[lt]{\lineheight{1.25}\smash{\begin{tabular}[t]{l}no-flow\end{tabular}}}}}%
    \put(0,0){\includegraphics[width=\unitlength,page=2]{domain_single_frac.pdf}}%
    \put(0.51866719,0.33870657){\color[rgb]{0,0,1}\makebox(0,0)[lt]{\lineheight{1.25}\smash{\begin{tabular}[t]{l}$\gamma$\end{tabular}}}}%
    \put(0.52389125,0.57606337){\color[rgb]{0.85098039,0,0}\makebox(0,0)[lt]{\lineheight{1.25}\smash{\begin{tabular}[t]{l}$\mu^{\rm top}$\end{tabular}}}}%
    \put(0.10899575,0.82610803){\color[rgb]{0,0,0}\makebox(0,0)[lt]{\lineheight{1.25}\smash{\begin{tabular}[t]{l}$\Omega$\end{tabular}}}}%
    \put(0,0){\includegraphics[width=\unitlength,page=3]{domain_single_frac.pdf}}%
    \put(0.68189648,0.51267539){\color[rgb]{0.85098039,0,0}\makebox(0,0)[lt]{\lineheight{1.25}\smash{\begin{tabular}[t]{l}$\mu^{\rm bottom}$\end{tabular}}}}%
  \end{picture}%
\endgroup%